\documentclass[11pt,a4paper]{article}
\usepackage[T1]{fontenc}
\usepackage[utf8]{inputenc}
\usepackage[english]{babel}
\usepackage{lmodern}
\usepackage[margin=25mm,headheight=15pt]{geometry}
\usepackage{amsmath,amssymb,amsthm,mathtools,booktabs,array,tabularx}
\usepackage{enumitem,fancyhdr,needspace,xcolor,microtype}
\usepackage[unicode,hypertexnames=false,colorlinks=true,linkcolor=blue!40!black,citecolor=blue!40!black,urlcolor=blue!40!black]{hyperref}
\newtheorem{theorem}{Theorem}[section]
\newtheorem{lemma}[theorem]{Lemma}
\newtheorem{proposition}[theorem]{Proposition}
\newtheorem{corollary}[theorem]{Corollary}
\theoremstyle{remark}\newtheorem{remark}[theorem]{Remark}

\DeclareMathOperator{\supp}{supp}
\numberwithin{equation}{section}
\setlist[enumerate]{itemsep=3pt,topsep=4pt}
\hypersetup{pdftitle={Horizontal estimates for Kloosterman sums under Hecke eigenvalue constraints},pdfsubject={English manuscript},pdflang={en-US}}
\newcommand{\papertitle}[1]{\begin{center}\LARGE\bfseries #1\par\end{center}\vspace{5pt}}
\begin{document}
\thispagestyle{plain}
\papertitle{Horizontal estimates for Kloosterman sums\\under Hecke eigenvalue constraints}

\begin{center}
\normalsize
Hongze Li and Mingxuan Zhong
\strut
\end{center}
\begin{abstract}
We prove weighted cancellation under small-value conditions
on Hecke eigenvalues of finitely many non-CM holomorphic newforms. We also obtain
many squarefree integers of each sign with fixed normalized amplitude,
simultaneously small Hecke eigenvalues, and a fixed bound on the number
of prime factors, uniformly in growing reduced residue classes.
Our questions are motivated by Xi's noncorrelation conjectures and
the mixed-moment estimates of Jiang and Lau. Their estimates imply that
small Hecke eigenvalues can occur together with non-small Kloosterman
values. Our result retains this property with a fixed prime-factor
bound, both signs, and the residue conditions.
Using OpenAI's zero-free half-plane theorem for Dirichlet $L$-functions,
we give an example with at most $83$ prime factors, normalized amplitude
at least $10^{-21}$, and moduli $M\le C^{11/200}$.
\end{abstract}
\begingroup
\renewcommand{\thefootnote}{}
\makeatletter
\renewcommand{\@makefntext}[1]{\noindent#1}
\makeatother
\footnotetext[0]{%
\textbf{Keywords:} Kloosterman sums; Hecke eigenvalues; sieve methods;
sign changes; Sato--Tate distribution.\par
\noindent\textbf{2020 MSC:} 11L05(11F30, 11N36).}
\endgroup
\section{Introduction}
\subsection{Background}
\label{en:intro:background}

Kloosterman sums first appeared in Poincar\'e's work on modular forms
and were later used by Kloosterman to study representations by
quadratic forms \cite{en:poincare,en:kloosterman}.
Their connections with modular forms, spectral theory, and
equidistribution have led to many applications in analytic number
theory. Michel surveys these developments under the name
\emph{Kloostermania}, the term proposed by Huxley
\cite{en:michelKloostermania}.
We use the notation
\[
 \begin{gathered}
 S(a,b;c)=\sum_{v\bmod c}^{*}
       e\!\left(\frac{av+b\bar v}{c}\right),\\
 K(n)=\frac{S(1,1;n)}{\sqrt n},\qquad
 x_n=2^{-\omega(n)}K(n).
 \end{gathered}
\]
Here $e(t)=\exp(2\pi it)$, the star restricts the sum to reduced
residue classes, $\bar v$ is the inverse of $v\bmod c$, and
$\omega(n)$ is the number of distinct prime factors of $n$.
These sums are real. Weil's bound \cite{en:weil} gives $|K(p)|\le2$
for primes $p$, and the Chinese remainder theorem gives $|x_n|\le1$
for squarefree $n$.

For a fixed nonzero integer $a$, Katz asked about the signs and
distribution of $p^{-1/2}S(1,a;p)$ as $p$ varies over the primes
\cite{en:katz80}. Write
$p^{-1/2}S(1,a;p)=2\cos\theta_p(a)$ for $p\nmid a$.
The horizontal Sato--Tate conjecture predicts that these angles
are distributed on $[0,\pi]$ with measure
$(2/\pi)\sin^2\theta\,d\theta$. In particular, it predicts that each
sign has density $1/2$ among the primes. This remains open.
Katz proved the vertical version, in which $a$ varies modulo $p$
and then $p$ tends to infinity \cite{en:katz}.

Over function fields, Fu--Lau--Li--Xi \cite{en:fllx} prove
horizontal Sato--Tate laws with explicit error terms as the places
vary in arithmetic progressions or short intervals. They also
establish a joint version for suitable finite families of
Kloosterman sums.

One direction is to study signs on integers with few prime factors.
Fouvry--Michel \cite{en:fmchanges} obtained
$\gg X/\log X$ squarefree moduli of each sign in $(X,2X]$ with at
most $23$ prime factors. The bound was reduced in subsequent work
\cite{en:sivak,en:matomakiSigns,en:xiSigns,en:xiSignsII,
en:xiSignsCorrection,en:zhangzhong}.
Recent work of Xiao \cite{en:xiaoFive} and
Xiao--Li \cite{en:xiaoli26I,en:xiaoli26} gives the same counting
order with at most five and four prime factors, respectively.
Drappeau--Maynard \cite{en:dm} obtain sign changes at products of
two primes under an exceptional-zero hypothesis.

Another direction concerns signed sums, absolute moments, and
arithmetic restrictions. Spectral methods give cancellation in
progressions and over squarefree moduli
\cite{en:bm,en:eass2024,en:steinerPeriodic}.
Absolute moments were studied by Fouvry--Michel
\cite{en:fmabsolute,en:fmabsoluteCorr}; recent work of
Xi \cite{en:xidensity26} and Jiang--Lau \cite{en:jl26} determines
the order of the absolute first moment:
\[
 \sum_{n\le X}|K(n)|
 \asymp X(\log X)^{8/(3\pi)-1}.
\]
These estimates distinguish cancellation of signs from small
absolute values.

The connection with Hecke eigenvalues goes back to Katz's third
problem \cite[Chapter~1]{en:katz80}; see also
\cite[Section~1]{en:xihecke}. For a fixed nonzero integer $a$, he
asked whether the Euler product formed from the prime values
$p^{-1/2}S(1,a;p)$ could be the $L$-function of a Maass form whose
level is a power of $2$. This would identify these Kloosterman
values with the Hecke eigenvalues at the unramified primes.
For a fixed holomorphic newform, such an identification can already
be ruled out by algebraic properties \cite[Section~1]{en:xihecke}.

For a primitive Hecke--Maass cusp form $f$ with trivial nebentypus,
Xi \cite[Theorems~1.1--1.2]{en:xihecke} proves
$\lambda_f(n)\ne\pm K(n)$ for infinitely many squarefree integers
with a bounded number of prime factors, where $\lambda_f(n)$ is
the normalized Hecke eigenvalue. He also obtains both signs of
$\lambda_f(n)-\eta K(n)$ for each fixed nonzero real $\eta$.
For a fixed primitive holomorphic or Maass cusp form $f$, he further
conjectures \cite[Conjectures~2.1--2.2]{en:xihecke}
\[
 \sum_{p\le X}\lambda_f(p)K(p)=o(X/\log X),
 \qquad
 \sum_{n\le X}\lambda_f(n)K(n)\ll_f X(\log X)^{-2018}.
\]
Jiang--Lau \cite[Section~4.3]{en:jl26} study mixed absolute
moments of $|\lambda_f(n)|$ and $|K(n)|$ for full-level
holomorphic eigenforms. Related results on small Hecke coefficients
and on exact equality with multiplicative functions can be found in
\cite{en:lrsTypical,en:yangzhangDistinct}.

\subsection{Our results and the plan of the paper}
\label{en:intro:results}

Motivated by Xi's noncorrelation conjectures and the mixed-moment
results of Jiang--Lau, we ask whether Kloosterman sums still have
cancellation and both signs under small-value conditions on
several Hecke eigenvalues.
Fix normalized non-CM holomorphic newforms $f_1,\ldots,f_s$ with
trivial nebentypus, and positive numbers $A_1,\ldots,A_s$.
Put $A=A_1+\cdots+A_s$ and $L=\log C$.
We impose small-eigenvalue conditions at selected prime divisors
and construct a nonnegative arithmetic weight $W_C$.
For the three weighted sums defined in
Section~\ref{en:intro:notation}, Theorem~\ref{en:thm:main-hecke}
gives
\[
 T_C=o(D_C),\qquad
 B_C\asymp D_C\asymp\frac{C}{\varphi(M)L^{1+A}},
\]
uniformly in each reduced class $a\pmod M$ with $M\le C^\delta$
for a fixed $\delta>0$. In particular, the signed sum is
negligible compared with the absolute sum for the same weight.

Corollary~\ref{en:app:finite-forms} gives fixed constants
$b,c_1,\ldots,c_s>0$ and an integer $R_0$ such that each sign
$\sigma\in\{-1,1\}$ occurs on
\[
 \gg\frac{C}{\varphi(M)L^{1+A}}
\]
integers satisfying
\[
 \begin{gathered}
 C<n<2C,\qquad \mu^2(n)=1,\qquad n\equiv a\pmod M,\\
 \omega(n)\le R_0,\qquad \sigma x_n\ge b,\qquad
 |\lambda_{f_i}(n)|\le c_iL^{-A_i}\quad(1\le i\le s).
 \end{gathered}
\]
Thus all the small-value conditions hold on the same integers,
and the bound on the number of prime factors is independent of $C$.
The signed estimate above concerns the full weight; the conclusion
after imposing $\omega(n)\le R_0$ is a separate count for each sign.

These conclusions concern Kloosterman sums under Hecke restrictions;
they do not establish Xi's mixed-sum conjectures. The estimates of Jiang--Lau
\cite[Theorem~3.8 and Corollary~4.5]{en:jl26} already imply, for
full-level holomorphic forms, that simultaneously small Hecke
eigenvalues can occur together with $|K(n)|\ge1$.
This deduction from high moments does not give a fixed bound for
$\omega(n)$. Here we retain such a bound, both signs with fixed
normalized amplitude, and the stated count in every permitted
reduced residue class.

The proof uses sieve weights, spectral estimates, local
equidistribution, and effective Sato--Tate estimates. Different
marked primes are assigned to different forms, so single-form
Sato--Tate estimates suffice.
Using the zero-free half-plane theorem for Dirichlet $L$-functions
in the OpenAI preprint \cite{en:oaizero26},
Theorem~\ref{en:num:theorem} gives a two-marker example with
\[
 R_0=83,\qquad b=10^{-21},\qquad M\le C^{11/200}.
\]
The general theorem does not require this additional input.

Section~\ref{en:unif:definition} defines the weights and states
the main results. Section~\ref{en:unif:signed-proof} proves the
signed estimate and the positive upper bounds.
Section~\ref{en:unif:absolute-proof} gives the absolute lower bound.
Section~\ref{en:app:main} derives the small Hecke eigenvalues and
the counts for both signs.
Section~\ref{en:num:main} gives the explicit example.

\subsection{Notation}
\label{en:intro:notation}

All variables in summation are positive integers. We write
$P^+(n)$ and $P^-(n)$ for the largest and smallest prime factors
of $n$, with $P^+(1)=1$ and $P^-(1)=\infty$, and use $*$ for
Dirichlet convolution.
Fix a nonzero function $0\le g\in C_c^\infty((1,2))$.
For a nonnegative weight $W_C$, which may depend on $C$, set
\begin{equation}\label{en:eq:three-masses}
 \begin{aligned}
 D_C&=\sum_n\mu^2(n)g(n/C)W_C(n),\\
 T_C&=\sum_n\mu^2(n)g(n/C)x_nW_C(n),\\
 B_C&=\sum_n\mu^2(n)g(n/C)|x_n|W_C(n).
 \end{aligned}
\end{equation}
The notation $\asymp$ means two-sided bounds by positive constants.

For a normalized holomorphic newform $f$ of weight $w_f\ge2$,
level $N_f$, and trivial nebentypus, write
\[
 f(z)=\sum_{n\ge1}a_f(n)e(nz),\qquad
 \lambda_f(n)=a_f(n)n^{-(w_f-1)/2}.
\]
At good primes, meaning $p\nmid N_f$, write
$\lambda_f(p)=2\cos\vartheta_f(p)$ with
$0\le\vartheta_f(p)\le\pi$.
For squarefree $n$, Hecke multiplicativity gives
$\lambda_f(n)=\prod_{p\mid n}\lambda_f(p)$.
The bound $|\lambda_f(p)|\le2$ also holds at bad primes
($p\mid N_f$). We use the equivalent Sato--Tate measures
\[
 d\mu_{\rm ST}(\vartheta)=\frac2\pi\sin^2\vartheta\,d\vartheta,
 \qquad
 d\nu_{\rm ST}(t)=\frac1{2\pi}\sqrt{4-t^2}\,
                   \mathbf1_{[-2,2]}(t)\,dt.
\]
Constants concerning the marked primes may depend on the fixed
forms and exponents, but not on $C$ or on the admissible modulus
and residue classes.

\section{Arithmetic weights and main results}
\label{en:unif:definition}

Fix $0<\eta<1/2$, $0<\beta<1/10$ is sufficiently small relative
to $\eta$. Put $y=C^\eta$ and $z=C^\beta$. Use upper-bound beta-sieve
coefficients of dimension $4$, as in \cite[\S3.2]{en:dm}, satisfying
\begin{equation}\label{en:unif:sieve}
 \begin{gathered}
 \xi_1=1,\quad |\xi_e|\le1,\quad
 \theta_C(n)=\sum_{e\mid n}\xi_e\ge0,\\
 \xi_e=0\quad\hbox{unless }e\le y,\quad P^+(e)\le z,\quad\mu^2(e)=1.
 \end{gathered}
\end{equation}
The fixed ratio $\eta/\beta$ is taken sufficiently large 
so that the estimations in the following text are valid. Define
\begin{equation}\label{en:unif:positive-weight}
 \Theta_C=\Lambda*\theta_C=\xi_C*\log,\qquad
 A_C(n)=\frac{\Theta_C(n)}{\log n}\ (n\ge2),\qquad A_C(1)=0.
\end{equation}
Both factors in the first convolution are nonnegative, so $A_C\ge0$.
Moreover, $A_C(n)=1$ whenever $P^-(n)>z$. 

Let $k\ge0$ be fixed. When $k=0$, use the empty marking convention
$H_C=F=d(C)=1$, $P=1$, and $\sigma=M_0=0$; the parameters and assumptions
on prime restrictions below are then unnecessary. For $k\ge1$, fix
$\alpha>\beta$ and $0<\sigma<1/2-\eta$, and a marked compact set
\[
 \mathcal K\subset\{(t_1,\ldots,t_k):
 \alpha<t_1<\cdots<t_k,\ \textstyle\sum_jt_j<\sigma\}.
\]
Take a fixed function $0\le F\le1$ supported on $\mathcal K$, with
$F\ge c_F>0$ on a fixed open box 
\[
\prod_j(a_j,b_j)=\left\{(t_1,\ldots,t_k)\in\mathbb R^k:
a_j<t_j<b_j\ \text{for all }j\right\}=Q\subset\mathcal K.
\]

Neither smoothness nor separation of variables is required of $F$.
For each $j$, specify a set of primes $\mathcal P_{j,C}$ and a density
$0<d_j(C)\le1$. For every fixed $\lambda>1$ sufficiently close to $1$, assume
uniformly in the relevant fixed power ranges that
\begin{equation}\label{en:unif:prime-colors}
 \#(\mathcal P_{j,C}\cap[Y,\lambda Y))\asymp d_j(C)Y/L.
\end{equation}
It suffices to have upper bounds near the coordinate projections of the
support and lower bounds in the specified box. For instance, one may impose
this estimate throughout $C^\alpha\le Y\le C^\sigma$. Put
\begin{equation}\label{en:unif:marker}
 \begin{split}
 H_C(n)&=\sum_{\substack{p_1<\cdots<p_k\\p_1\cdots p_k\mid n}}
 F\!\left(\frac{\log p_1}L,\ldots,\frac{\log p_k}L\right)
       \prod_{j=1}^k\mathbf1_{p_j\in\mathcal P_{j,C}},\\
 d(C)&=\prod_{j=1}^k d_j(C)\ge L^{-M_0},
 \end{split}
\end{equation}
where $M_0$ is fixed. Since every marked prime exceeds $C^\alpha$, we have
$H_C(n)\ll1$ for $n\asymp C$. All given rules are prescribed independently
of the Kloosterman values.

Let $M=M(C)\le C^\delta$ and choose any nonempty set
$\mathcal R_C\subseteq(\mathbb Z/M\mathbb Z)^\times$. Define
\begin{equation}\label{en:unif:weight-density}
 \rho_C=\frac{|\mathcal R_C|}{\varphi(M)},\qquad
 \Delta_C=d(C)\rho_C\frac CL,\qquad
 W_C(n)=2^{\omega(n)}A_C(n)H_C(n)
               \mathbf1_{n\bmod M\in\mathcal R_C}.
\end{equation}
For $M=1$, the unit group is interpreted as a singleton and $\rho_C=1$.

\subsection{Weighted estimates retaining the conditional density}
\label{en:unif:main-statement}

In \eqref{en:eq:three-masses}, we use the weight defined in \eqref{en:unif:weight-density}. 
For $\kappa=2,4$, also write
\[
 D_{\kappa,C}=\sum_n\mu^2(n)g(n/C)\kappa^{\omega(n)}A_C(n)H_C(n)
                  \mathbf1_{n\bmod M\in\mathcal R_C};\qquad D_{2,C}=D_C.
\]

\begin{proposition}[Cancellation for nonnegative arithmetic weights]\label{en:thm:unified}
Under the preceding fixed-parameter and prime-restriction assumptions, there exist
$\delta_0>0$ and a sufficiently small fixed $\delta>0$ such that, uniformly
for $M\le C^\delta$, nonempty $\mathcal R_C$, and the admissible restricted sets,
\begin{equation}\label{en:unif:main}
 \begin{gathered}
 D_C\asymp B_C\asymp\Delta_C\longrightarrow\infty,\qquad
 |T_C|\ll |\mathcal R_C|C^{1-\delta_0}=o(\Delta_C),\\
 D_{4,C}\ll\Delta_C.
 \end{gathered}
\end{equation}
Constants may depend on the fixed parameters and the uniform comparison
constants in \eqref{en:unif:prime-colors}. 
\end{proposition}

\subsection{The main theorem with Hecke conditions at marked primes}

Fix normalized non-CM holomorphic newforms $f_1,\ldots,f_s$, all with trivial
nebentypus, and positive real numbers $A_1,\ldots,A_s$. Choose integers
$k_i>2A_i$ and put $k=\sum_i k_i$. Partition the marking indices into disjoint
sets $J_i$ with $|J_i|=k_i$. For $j\in J_i$, set
\begin{equation}\label{en:def:hecke-colors}
 \mathcal P_{j,C}=\{p:p\nmid N_{f_i},\
          |\lambda_{f_i}(p)|\le L^{-A_i/k_i}\}.
\end{equation}
Use $H_C$ above, with a function $F$ satisfying
the stated support conditions. Such parameters always exist: take the
centers of the position intervals to be $t_j=j/(2k(k+1))$, whose sum is
$1/4$, choose sufficiently narrow disjoint intervals, set $\eta=1/10$
and $0.3<\sigma<0.4$, and finally take $\beta<\alpha$ sufficiently small.
All these choices are fixed before $C$ tends to infinity.

Write
\[
 m_i(C)=\nu_{\rm ST}([-L^{-A_i/k_i},L^{-A_i/k_i}]),\quad
 d(C)=\prod_{i=1}^s m_i(C)^{k_i},\quad A=\sum_{i=1}^s A_i.
\]
Direct integration gives
\begin{equation}\label{en:eq:hecke-density}
 m_i(C)=\frac2\pi L^{-A_i/k_i}+O(L^{-3A_i/k_i}),
 \qquad d(C)\asymp L^{-A}.
\end{equation}

\begin{theorem}[Cancellation and absolute mass under Hecke prime-factor conditions]\label{en:thm:main-hecke}
With these choices above, there exist fixed $\delta,\delta_0>0$ such that, for all
sufficiently large $C$, all $M\le C^\delta$, and all nonempty
$\mathcal R_C\subseteq(\mathbb Z/M\mathbb Z)^\times$, the weight
\[
 W_C(n)=2^{\omega(n)}A_C(n)H_C(n)
                \mathbf1_{n\bmod M\in\mathcal R_C}
\]
satisfies
\begin{equation}\label{en:eq:hecke-main}
 \begin{gathered}
 D_C\asymp B_C\asymp
 \frac{|\mathcal R_C|C}{\varphi(M)L^{1+A}}\longrightarrow\infty,\\
 |T_C|\ll |\mathcal R_C|C^{1-\delta_0}=o(D_C),\qquad
 \sum_n\mu^2(n)g(n/C)4^{\omega(n)}A_C(n)H_C(n)
     \mathbf1_{n\bmod M\in\mathcal R_C}\ll D_C.
 \end{gathered}
\end{equation}
Furthermore, for every squarefree $n$ with $H_C(n)>0$, we have
\begin{equation}\label{en:eq:normalized-hecke-small}
 2^{-\omega(n)}|\lambda_{f_i}(n)|\le2^{-k_i}L^{-A_i}
 \qquad(1\le i\le s).
\end{equation}
The weights are specified by sieve coefficients, prime-factor positions,
Hecke conditions, and residue classes, independently of the magnitudes and
signs of the Kloosterman sums.
\end{theorem}

We prove Prop.~\ref{en:thm:unified} in the next two Sections and Thm.~\ref{en:thm:main-hecke} 
in Section 5.
The two first-moment estimates together give
\[
 \frac{|T_C|}{B_C}\longrightarrow0,\qquad \frac{B_C}{D_C}\asymp1.
\]
The first ratio controls the difference between positive and negative
contributions; the second excludes degeneration of the absolute average.
Both are measured against the actual mass arising from the shrinking Hecke
conditions and the congruence restriction. Cancellation therefore remains
detectable within the constructed thin arithmetic family.

The main theorem and its application serve distinct purposes. The normalization
in \eqref{en:eq:normalized-hecke-small} cannot be dropped on the full weighted
support; a fixed bound on $\omega(n)$ turns it into smallness of the global
Hecke eigenvalues. Corollary~\ref{en:app:finite-forms} combines cancellation,
the absolute lower bound, and the prime-factor tail bound to extract fixed-amplitude
witnesses of both signs while retaining these small-value and congruence
conditions. 
The conclusion after a fixed truncation of $\omega(n)$
does not assert a little-$o$ signed first moment under that truncation.

For a concrete example, use two marked primes for Ramanujan's form $\Delta$:
\begin{equation}\label{en:def:double-marker}
 \begin{split}
 H_C(n)=\sum_{\substack{p<q,\ pq\mid n}}
 &\mathbf1_{C^{0.36}<pq<C^{0.38}}
  \mathbf1_{C^{0.02}<q/p<C^{0.04}}\\[-2pt]
 &{}\cdot\mathbf1_{|\lambda_\Delta(p)|\le L^{-1/4}}
         \mathbf1_{|\lambda_\Delta(q)|\le L^{-1/4}},
 \qquad \lambda_\Delta(n)=\tau(n)n^{-11/2}.
 \end{split}
\end{equation}
With the choice $\eta=1/2000$ from Section~\ref{en:num:main}, the total divisor
length is less than $C^{761/2000}<\sqrt{C}$. The scale
in a single reduced class is $C/(\varphi(M)L^{3/2})$. In particular, for
$M=2^{\lfloor\log_2 L\rfloor}$ and $\mathcal R_C=\{1\}$,
\[
 D_C\asymp B_C\asymp C/L^{5/2},\qquad
 T_C\ll C^{1-\delta_0}=o(C/L^{5/2}).
\]
Section~\ref{en:num:main} supplies explicit factor-count and amplitude constants,
together with an explicit range of arbitrary moduli, for this example,
with the additional zero-free input stated there.

\section{Proof of Proposition~\ref{en:thm:unified}: Positive sieve mass and signed cancellation}
\label{en:unif:signed-proof}
\subsection{The sieve input and positivity of the convolution}

We use the fundamental lemma of the upper-bound beta sieve in dimension $4$.
The following is the precise consequence required here. If $y=z^u$ and the
fixed number $u$ is sufficiently large, the same coefficients in
\eqref{en:unif:sieve} can be chosen so that
\begin{equation}\label{en:eq:fundamental-input}
 \sum_e\xi_e h(e)\ll\prod_{p\le z}(1-h(p)),\qquad
 h(p)=\frac{\kappa}{p+\kappa},\quad\kappa=2,4.
\end{equation}
Here $h$ is multiplicative and supported on squarefree integers. The general
fundamental lemma assumes
$\prod_{w\le p<v}(1-h(p))^{-1}\le K(\log v/\log w)^4$.
Mertens' formula verifies this condition for the two local densities above,
with fixed $K$. The coefficients depend on the dimension and sieve parameters,
so they can be used for both densities. This is a known sieve input; see \cite{en:opera},
\cite[\S3.2]{en:dm} and \cite[Lemma 9.1]{en:matomakiTeravainen}.
Section~\ref{en:r3:presieve} gives its pre-sieved form when explicit constants
are needed.

The identity $\Lambda*1=\log$ and associativity give
\begin{equation}\label{en:eq:convolution-identity}
 \Lambda*(\xi_C*1)=\xi_C*\log,
 \qquad A_C(n)=\sum_{e\mid n}\xi_e\frac{\log(n/e)}{\log n}.
\end{equation}
The convolution on the left proves positivity of $A_C$, whereas the expression on
the right involves only divisors $e\le C^\eta$. If $P^-(n)>z$, only $e=1$
occurs and $A_C(n)=1$. For squarefree $n$, each logarithmic ratio lies in
$[0,1]$, whence
\begin{equation}\label{en:eq:pointwise-weight}
 0\le A_C(n)\le2^{\omega(n)}.
\end{equation}

\subsection{Upper bound for \texorpdfstring{$T_C$}{T\_C}}

We first estimate the complete character-twisted signed sum, and then impose
the congruence condition. Put
$\varepsilon=(1/2-\eta-\sigma)/2>0$. The short-divisibility result of
Drappeau--Maynard \cite[Proposition 3.1 and Corollary 3.2]{en:dm} gives a
number $\delta_1>0$, depending only on $\varepsilon$, such that for
$M\le C^{\delta_1}$, every character $\chi$ modulo $M$, and every fixed
$\phi\in C_c^\infty((0,\infty))$,
\begin{equation}\label{en:eq:dm-input}
 \sum_{q\le C^{1/2-\varepsilon}}
 \left|\sum_{q\mid n}\mu^2(n)\phi(n/C)\chi(n)K(n)\right|
 \ll_{\varepsilon,\phi}C^{1-\delta_1}.
\end{equation}
The bound is uniform in $M$ and $\chi$, including the principal character;
the factor $\mu^2(n)$ is already included in the result. The implied constant is depended
 on the fixed test $\phi$.

\paragraph{Expanding the signed sum.}
For an ordered marked tuple $p_1<\cdots<p_k$, write
\[
 w_C(P)=F\!\left(\frac{\log p_1}L,\ldots,
                 \frac{\log p_k}L\right)
                 \prod_{j=1}^k\mathbf1_{p_j\in\mathcal P_{j,C}},
 \qquad P=p_1\cdots p_k.
\]
Thus $0\le w_C(P)\le1$ and $H_C(n)=\sum_{P\mid n}w_C(P)$, where $P$
runs over products of $k$ distinct primes in increasing order. When $k=0$,
the only product is $P=1$, with $w_C(1)=1$. Since
$2^{\omega(n)}x_n=K(n)$, the character-twisted signed sum before imposing
any congruence condition is
\[
 \mathcal S(\chi)=\sum_n\mu^2(n)g(n/C)A_C(n)H_C(n)\chi(n)K(n).
\]
We now substitute \eqref{en:eq:convolution-identity}. On the support of
$g(n/C)$ we have $C<n<2C$, and
\[
 \frac{\log(n/e)}{\log n}
 =1-\frac{\log e}{\log n}
 =1-\frac{\log e}{L}\frac{L}{L+\log(n/C)}.
\]
For every nonzero term, all prime factors of $e$ are at most $z$ due to the support of $\xi_e$.
If $k\ge1$, every marked prime exceeds $C^\alpha>z$, so $(e,P)=1$;
the same conclusion is immediate when $k=0$. Hence the two conditions
$e\mid n$ and $P\mid n$ are equivalent to $eP\mid n$. Expanding the finite
sums gives the exact identity
\begin{equation}\label{en:eq:signed-exact-expansion}
 \begin{aligned}
 \mathcal S(\chi)
 &=\sum_{e,P}\xi_e w_C(P)
   \sum_{\substack{eP\mid n\\ C<n<2C}}
       \mu^2(n)g(n/C)\chi(n)K(n)\\
 &\quad-\sum_{e,P}\xi_e\frac{\log e}{L}w_C(P)
   \sum_{\substack{eP\mid n\\ C<n<2C}}\mu^2(n)
        \frac{Lg(n/C)}{L+\log(n/C)}\chi(n)K(n).
 \end{aligned}
\end{equation}

To put each term into the form required by \eqref{en:eq:dm-input}, set
$q=eP$. The product of the prime-power factors of $q$ belonging to primes
at most $z$ recovers $e$; the remaining prime factors, in increasing order,
recover the marked tuple. Thus each $q$ has at most one representation
contributing to either sum. For $k=0$, this simply says $q=e$.
Furthermore,
\[
 q\le C^{\eta+\sigma}=C^{1/2-2\varepsilon}
 \le C^{1/2-\varepsilon}.
\]
After reindexing by $q$, the coefficients in the first and second sums are,
respectively,
\[
 b_q^{(1)}=\xi_e w_C(P),\qquad
 b_q^{(2)}=\xi_e\frac{\log e}{L}w_C(P),
\]
with both coefficients set to zero if there is no such representation.
Since $|\xi_e|\le1$ and $0\le\log e/L\le\eta<1$, both have absolute
value at most $1$.

The two corresponding test functions are
\[
 g(t),\qquad G_C(t)=\frac{Lg(t)}{L+\log t}.
\]
Here $G_C$ is defined on $(1,2)$ and extended by zero outside that interval;
it is smooth because $g$ has compact support in $(1,2)$. The second test
depends on $C$, whereas \eqref{en:eq:dm-input} is stated for a fixed test.
We therefore need a uniform version of that input for this family of tests.

\paragraph{Uniformity for varying smooth tests.}
Fix a compact interval $I\Subset(1,2)$ whose interior contains
$\supp g$, and fix $C_0>1$. Let
\[
 \mathcal E_I=\{\phi\in C^\infty(\mathbb R;\mathbb C):
                         \supp\phi\subseteq I\},\qquad
 p_r(\phi)=\max_{0\le j\le r}\|\phi^{(j)}\|_\infty
 \quad(r=0,1,2,\ldots).
\]
We give this space the topology defined by these increasing seminorms:
convergence means uniform convergence of every fixed-order derivative.
This topology is metrizable and locally convex. It is also complete:
if a sequence is Cauchy for every $p_r$, its derivatives of every order
have uniform limits; the fundamental theorem of calculus shows that these
limits are the derivatives of a smooth limit, which is still supported in
$I$. Thus $\mathcal E_I$ is a Fr\'echet space, so Baire's theorem applies.

Let $\mathfrak L$ be the family of linear functionals
\[
 \mathcal L(\phi)=C^{-1+\delta_1}
 \sum_{q\le C^{1/2-\varepsilon}}b_q
 \sum_{q\mid n}\mu^2(n)\phi(n/C)\chi(n)K(n),
\]
where the family ranges over all $C\ge C_0$, all integers
$M\le C^{\delta_1}$, all characters $\chi$ modulo $M$, and all complex
coefficients $|b_q|\le1$. Each fixed $\mathcal L$ is a finite linear
combination of evaluations $\phi(n/C)$, because $I$ is compact. In
particular, for some finite constant $A_{\mathcal L}$,
\[
 |\mathcal L(\phi)-\mathcal L(\psi)|
 \le A_{\mathcal L}p_0(\phi-\psi).
\]
This proves continuity of each functional; no uniform bound for
$A_{\mathcal L}$ is asserted here. On the other hand, the triangle
inequality and \eqref{en:eq:dm-input} give pointwise boundedness of the
whole family:
\[
 \sup_{\mathcal L\in\mathfrak L}|\mathcal L(\phi)|<\infty
 \qquad\text{for every fixed }\phi\in\mathcal E_I.
\]
If the quoted asymptotic estimate starts at a threshold depending on
$\phi$, the remaining bounded range of $C$ is controlled directly by
the finite sums. Indeed, on any fixed interval $C_0\le C\le C_1$, the
integers $n,q,M$ lie in fixed finite ranges, $|\phi(n/C)|\le p_0(\phi)$,
and $|\chi(n)|,|b_q|\le1$, while $C^{-1+\delta_1}$ is bounded. The bound
in this display may still depend on $\phi$.

For $m=1,2,\ldots$, define
\[
 E_m=\bigcap_{\mathcal L\in\mathfrak L}
       \{\phi\in\mathcal E_I:|\mathcal L(\phi)|\le m\}.
\]
Each $E_m$ is closed by continuity of the functionals, and pointwise
boundedness gives $\mathcal E_I=\bigcup_{m\ge1}E_m$. Baire's theorem
states that a complete metrizable space cannot be a countable union of
closed sets with empty interior. Consequently, some $E_m$ contains an open
neighborhood of a function $\phi_0$. A basic neighborhood is specified by
finitely many seminorm inequalities. Since the $p_r$ are increasing, we
may take the largest of their indices and the smallest of their radii to
obtain one finite integer $r$ and one number $a>0$ such that
\[
 p_r(h)<a\quad\Longrightarrow\quad\phi_0+h\in E_m.
\]
Here $h$ is a perturbation of the center $\phi_0$. The integer $r$ and
the constants $m,a$ are obtained once for the entire family $\mathfrak L$;
they are independent of $C,M,\chi,b_q$, and of the test subsequently
inserted. In particular, $r$ is a fixed finite order furnished by this
argument, rather than an order chosen separately for each test.

Taking $h=0$ shows that $\phi_0\in E_m$. For $p_r(h)<a$, linearity and
the bounds at $\phi_0+h$ and $\phi_0$ therefore give
\[
 |\mathcal L(h)|
 =|\mathcal L(\phi_0+h)-\mathcal L(\phi_0)|\le2m
 \qquad(\mathcal L\in\mathfrak L).
\]
This is now a bound on a neighborhood of the zero function. For any
$\phi\ne0$, we have $p_r(\phi)>0$ because $p_r$ includes $p_0$.
Set $h=a\phi/(2p_r(\phi))$, so that $p_r(h)=a/2<a$. Then
\[
 \frac{a}{2p_r(\phi)}|\mathcal L(\phi)|
 =|\mathcal L(h)|\le2m.
\]
The case $\phi=0$ is immediate. Thus, with the fixed constant $B=4m/a$,
we obtain the uniform estimate
\begin{equation}\label{en:eq:uniform-test}
 |\mathcal L(\phi)|\le Bp_r(\phi)
 \qquad(\phi\in\mathcal E_I,\ \mathcal L\in\mathfrak L).
\end{equation}
This proves that the fixed-test input is uniform on any family supported
in $I$ whose derivatives of each fixed order are uniformly bounded.
The derivative bound may depend on the order; one constant controlling
all derivative orders is not required.

\paragraph{The character-twisted bound.}
Both $g$ and $G_C$ belong to $\mathcal E_I$. On $I$, the function $\log t$
and all of its fixed-order derivatives are bounded, and
$L+\log t\asymp L$ uniformly for $C\ge C_0$. For example,
\[
 G_C'(t)=\frac{Lg'(t)}{L+\log t}
          -\frac{Lg(t)}{t(L+\log t)^2}.
\]
Repeated differentiation, or the same calculation applied to
$G_C(t)=g(t)/(1+\log t/L)$, gives
\[
 \sup_{C\ge C_0}p_s(G_C)<\infty
 \qquad\text{for every fixed integer }s\ge0.
\]
Apply \eqref{en:eq:uniform-test} at the particular fixed order $r$ obtained
above to the two coefficient sequences $b_q^{(1)},b_q^{(2)}$ and the two
tests in \eqref{en:eq:signed-exact-expansion}. Removing the normalization
$C^{-1+\delta_1}$ in the definition of $\mathcal L$ yields
\[
 |\mathcal S(\chi)|
 \le B C^{1-\delta_1}\bigl(p_r(g)+p_r(G_C)\bigr)
 \ll C^{1-\delta_1},
\]
uniformly for $M\le C^{\delta_1}$ and every character $\chi$ modulo $M$.
The coefficients may vary with the marking rules, but their absolute
values are at most $1$, so the same bound is uniform in those rules.
No power saving is lost in the uniformity argument, and we may set
$\delta_0=\delta_1$.

\paragraph{Imposing the congruence condition.}
For a reduced residue class $a$ modulo $M$, put
\[
 T_a=\sum_{n\equiv a\ ({\rm mod}\ M)}
       \mu^2(n)g(n/C)A_C(n)H_C(n)K(n).
\]
Character orthogonality gives, on all integers,
\[
 \mathbf1_{n\equiv a\ ({\rm mod}\ M)}
 =\frac1{\varphi(M)}\sum_{\chi\ ({\rm mod}\ M)}
                       \overline{\chi(a)}\chi(n).
\]
Indeed, on integers not coprime to $M$ both sides vanish, since $a$ is
reduced and Dirichlet characters are extended by zero there. Hence
\[
 T_a=\frac1{\varphi(M)}\sum_{\chi\ ({\rm mod}\ M)}
             \overline{\chi(a)}\mathcal S(\chi),
 \qquad |T_a|\ll C^{1-\delta_0}.
\]
Finally, the classes in $\mathcal R_C$ are disjoint, and
$T_C=\sum_{a\in\mathcal R_C}T_a$. Using
$\Delta_C=d(C)|\mathcal R_C|C/(\varphi(M)L)$ gives
\begin{equation}\label{en:unif:signed-union}
 \begin{aligned}
 |T_C|&\ll |\mathcal R_C|C^{1-\delta_0},\\
 \frac{|T_C|}{\Delta_C}
 &\ll\frac{\varphi(M)L}{d(C)C^{\delta_0}}
 \le C^{\delta-\delta_0}L^{M_0+1}=o(1).
 \end{aligned}
\end{equation}
Here we used $\varphi(M)\le M\le C^\delta$ and
$d(C)\ge L^{-M_0}$. Choosing $\delta<\delta_0/2$ ensures both the modulus
range required above and convergence of the last expression to zero.

\subsection{Upper bound of \texorpdfstring{$D_C$}{D\_C}}
\label{en:unif:positive-proof}
\paragraph{Positive mass without congruence conditions.}
Our first goal is the estimate
\begin{equation}\label{en:unif:base-positive}
 \sum_{n\asymp C}\mu^2(n)\kappa^{\omega(n)}A_C(n)\ll C/L,
 \qquad\kappa=2,4,
\end{equation}
where $n\asymp C$ denotes any fixed interval
$c_1C<n\le c_2C$, with $0<c_1<c_2<\infty$. The implied constant may
depend on this interval and on the fixed sieve parameters. For $\kappa=2$,
this is the upper bound underlying $D_C$ before introducing the marked
primes and the congruence condition; the case $\kappa=4$ will also be
needed for the strengthened positive moment. We prove the estimate by
the nonnegative sieve-weight argument of Drappeau--Maynard
\cite[Lemma 3.3]{en:dm}, giving the details for both values of $\kappa$.

Put $f_\kappa(n)=\mu^2(n)\kappa^{\omega(n)}$. Since
$A_C(n)=\Theta_C(n)/\log n$ and $\log n\asymp L$ on the desired interval,
we first estimate
\[
 \sum_{n\le X}f_\kappa(n)\Theta_C(n)\ll X
 \qquad(X\asymp C).
\]
Throughout the argument, $C$ is sufficiently large that $y<X$.
The identity $\Theta_C=\Lambda*\theta_C$ gives, for squarefree $n$,
\[
 \Theta_C(n)=\sum_{p\mid n}\log p\,\theta_C(n/p),
\]
because the only divisors of a squarefree integer on which $\Lambda$
is nonzero are primes. On writing $n=pm$, squarefreeness is equivalent
to $m$ being squarefree and $p\nmid m$, and then
$f_\kappa(pm)=\kappa f_\kappa(m)$. Hence
\[
 \sum_{n\le X}f_\kappa(n)\Theta_C(n)
 =\kappa\sum_{m\le X}f_\kappa(m)\theta_C(m)
             \sum_{\substack{p\le X/m\\p\nmid m}}\log p.
\]
All coefficients $f_\kappa(m)\theta_C(m)$ are nonnegative. We may
therefore drop $p\nmid m$ and apply Chebyshev's bound
$\sum_{p\le t}\log p\ll t$ to obtain
\begin{align*}
 \sum_{n\le X}f_\kappa(n)\Theta_C(n)
 &\ll_\kappa X\sum_{m\le X}
                       \frac{f_\kappa(m)\theta_C(m)}m\\
 &\ll_\kappa X\sum_{P^+(m)\le X}
                       \frac{f_\kappa(m)\theta_C(m)}m.
\end{align*}
The last step enlarges the summation set, again using nonnegativity.
Although the new set allows $m>X$, only squarefree products of primes
at most $X$ contribute, so there are still only finitely many nonzero
terms.

We next evaluate this enlarged sum. Expanding
$\theta_C(m)=\sum_{e\mid m}\xi_e$ and interchanging the finite sums gives
\begin{align*}
 \sum_{P^+(m)\le X}\frac{f_\kappa(m)\theta_C(m)}m
 &=\sum_e\xi_e
      \sum_{\substack{P^+(m)\le X\\e\mid m}}\frac{f_\kappa(m)}m\\
 &=\sum_e\xi_e\frac{f_\kappa(e)}e
               \prod_{\substack{p\le X\\p\nmid e}}(1+\kappa/p)\\
 &=\prod_{p\le X}(1+\kappa/p)
       \sum_e\xi_e\prod_{p\mid e}\frac{\kappa}{p+\kappa}.
\end{align*}

Now apply \eqref{en:eq:fundamental-input} with
$h(p)=\kappa/(p+\kappa)$. At every prime $p\le z$, the resulting
sieve factor cancels the corresponding Euler-product factor exactly:
\[
 (1+\kappa/p)\left(1-\frac{\kappa}{p+\kappa}\right)=1.
\]
It follows that
\[
 \sum_{P^+(m)\le X}\frac{f_\kappa(m)\theta_C(m)}m
 \ll \prod_{z<p\le X}(1+\kappa/p),
 \qquad
 \sum_{n\le X}f_\kappa(n)\Theta_C(n)
 \ll X\prod_{z<p\le X}(1+\kappa/p).
\]
To bound the remaining product, use $\log(1+u)\le u$ for $u\ge0$
and Mertens' formula:
\[
 \prod_{z<p\le X}(1+\kappa/p)
 \le\exp\!\left(\kappa\sum_{z<p\le X}\frac1p\right)
 \ll_\kappa\left(\frac{\log X}{\log z}\right)^\kappa
 \ll_{\beta,\kappa}1.
\]
The final bound uses $z=C^\beta$ and $X\asymp C$, so that
$\log X/\log z$ stays bounded. This proves the required estimate.

Finally, on $c_1C<n\le c_2C$ we have $\log n\asymp L$, and
$f_\kappa(n)\Theta_C(n)\ge0$. Taking $X=c_2C$, we therefore obtain
\[
 \sum_{c_1C<n\le c_2C}f_\kappa(n)A_C(n)
 =\sum_{c_1C<n\le c_2C}
       \frac{f_\kappa(n)\Theta_C(n)}{\log n}
 \ll\frac1L\sum_{n\le c_2C}f_\kappa(n)\Theta_C(n)
 \ll\frac CL,
\]
which is \eqref{en:unif:base-positive}. Any fixed bounded nonnegative
test $\psi(n/C)$ with $\supp\psi\subseteq[c_1,c_2]$ can be inserted by
multiplying the upper bound by $\|\psi\|_\infty$.
In particular, when $k=0$ and $H_C=1$,
this proves the desired estimate without congruence conditions.

\paragraph{Retaining the density under the prime restrictions.}
We now aim to prove
\[
 D_\kappa^{\mathrm{all}}
 :=\sum_n f_\kappa(n)g(n/C)A_C(n)H_C(n)
 \ll d(C)\frac CL,\qquad\kappa=2,4.
\]
Thus the upper bound retains the product $d(C)=\prod_jd_j(C)$ associated
with the restrictions $p_j\in\mathcal P_{j,C}$. The case $k=0$ follows
from \eqref{en:unif:base-positive}, so assume $k\ge1$.

We first expand the marked-prime weight. For
$\mathbf p=(p_1,\ldots,p_k)$ with $p_1<\cdots<p_k$, write
\[
 P=p_1\cdots p_k,\qquad
 F_C(\mathbf p)=F\!\left(\frac{\log p_1}L,\ldots,
                         \frac{\log p_k}L\right).
\]
Since $P$ is squarefree, writing $n=Pm$ gives
$f_\kappa(Pm)=\kappa^k f_\kappa(m)$ when $(P,m)=1$, and zero otherwise.
Consequently, the definition of $H_C$ gives the exact expansion
\[
 D_\kappa^{\mathrm{all}}
 =\kappa^k
   \sum_{\substack{p_1<\cdots<p_k\\p_j\in\mathcal P_{j,C}}}
       F_C(\mathbf p)
   \sum_{\substack{m\ge1\\(P,m)=1}}
       f_\kappa(m)g(Pm/C)A_C(Pm).
\]
For every nonzero term, the marked primes exceed $C^\alpha>z$.
Since all prime factors of a sieve divisor $e$ are at most $z$, we have
$e\mid Pm$ if and only if $e\mid m$. Expanding the two divisor sums yields
\begin{equation}\label{en:unif:positive-cofactor}
 \begin{aligned}
 \theta_C(Pm)&=\theta_C(m),\\
 \Theta_C(Pm)&=\sum_{e\mid m}\xi_e\log(Pm/e)
             =\Theta_C(m)+\log P\,\theta_C(m).
 \end{aligned}
\end{equation}
These identities hold even if $(P,m)>1$, provided all prime factors of
$P$ exceed $z$. They will allow us to compare the arithmetic weights
after replacing the marked primes.

Fix $\lambda>1$ sufficiently close to $1$ for
\eqref{en:unif:prime-colors} to apply. Partition the positive axis into
half-open intervals $[\lambda^u,\lambda^{u+1})$, $u\in\mathbb Z$.
For the marked coordinates, retain the Cartesian-product boxes
\[
 \mathcal J=J_1\times\cdots\times J_k,
 \qquad J_j=[Y_j,\lambda Y_j),\qquad P_*=\prod_{j=1}^kY_j,
\]
whose images in the coordinates $(\log p_1/L,\ldots,\log p_k/L)$ meet
$\mathcal K$. In particular, a box is a product of $k$ intervals for the
marked primes; the cofactor $m$ is kept separate. Compactness of
$\mathcal K$ inside the region
$\alpha<t_1<\cdots<t_k$, $\sum_jt_j<\sigma$ gives a fixed positive
margin in each of these inequalities. Each grid interval has logarithmic
width $\log\lambda/L$, which tends to zero. Thus, for sufficiently large
$C$, every tuple of primes in a retained box is strictly ordered and has
distinct entries, all exceeding $C^\alpha>z$. An arbitrary replacement
tuple need not lie in the support of $F$; below we use only $0\le F\le1$.

For each retained box, put
\[
 \mathcal M_{\mathcal J}
 =\left\{m\ge1:\frac{C}{\lambda^kP_*}<m<\frac{2C}{P_*}\right\}.
\]
If $\mathbf p\in\mathcal J$ and $g(Pm/C)\ne0$, then $C<Pm<2C$
and $P_*\le P<\lambda^kP_*$, so $m\in\mathcal M_{\mathcal J}$.
The partition therefore rewrites the expansion as
\[
 D_\kappa^{\mathrm{all}}
 =\sum_{\mathcal J}\sum_{m\in\mathcal M_{\mathcal J}}
       R_{\mathcal J}(m),
\]
where
\[
 R_{\mathcal J}(m)
 =\kappa^k f_\kappa(m)
   \sum_{\substack{p_j\in J_j\cap\mathcal P_{j,C}\\(P,m)=1}}
       F_C(\mathbf p)g(Pm/C)A_C(Pm).
\]
It thus suffices to compare this sum, for each fixed box and cofactor,
with the corresponding sum over all primes in the same intervals.

If $m\in\mathcal M_{\mathcal J}$ and $p'_j\in J_j$, write
$P'=\prod_jp'_j$. Every such replacement satisfies
\[
 \frac{C}{\lambda^k}<P'm<2\lambda^kC.
\]
On this fixed interval, $\log(P'm)\asymp L$.
Moreover, $\log P'=\log P_*+O_{k,\lambda}(1)$ and
$\log P_*\ge k\alpha L-O_{k,\lambda}(1)$. Since
$\Theta_C(m),\theta_C(m)\ge0$, \eqref{en:unif:positive-cofactor} gives
\[
 A_C(P'm)
 =\frac{\Theta_C(m)+\log P'\,\theta_C(m)}{\log(P'm)}
 \asymp V_{\mathcal J}(m)
 :=\frac{\Theta_C(m)+\log P_*\,\theta_C(m)}L.
\]
The constants are uniform over the boxes, cofactors, and replacement
tuples. If both terms in the numerator vanish, both sides are zero.
Thus all arithmetic weights within this fixed box and cofactor are
comparable, including those with $(P',m)>1$.

Let
\[
 N_j^{\rm res}=\#(J_j\cap\mathcal P_{j,C}),\qquad
 N_j=\#(J_j\cap\mathbb P).
\]
The upper bound in \eqref{en:unif:prime-colors} for primes satisfying the
specified restrictions gives $N_j^{\rm res}\ll d_j(C)Y_j/L$.
The prime number theorem gives $N_j\asymp Y_j/L$. Hence
\[
 \prod_{j=1}^kN_j^{\rm res}
 \ll d(C)\prod_{j=1}^kN_j.
\]
Using $F_C\le1$, $g\le\|g\|_\infty$, and nonnegativity to remove the
coprimality condition, we obtain
\begin{align*}
 R_{\mathcal J}(m)
 &\ll\kappa^k f_\kappa(m)V_{\mathcal J}(m)
                  \prod_jN_j^{\rm res}\\
 &\ll d(C)\kappa^k f_\kappa(m)V_{\mathcal J}(m)
                  \prod_jN_j\\
 &\ll d(C)\kappa^k f_\kappa(m)
                  \sum_{p'_j\in J_j\cap\mathbb P}A_C(P'm).
\end{align*}
This is the required replacement: the sum over primes satisfying the
restrictions is bounded by $d(C)$ times the sum over all primes in the
same intervals, with a fixed implied constant.

We now sum over all boxes and cofactors. The preceding bound yields
\[
 D_\kappa^{\mathrm{all}}
 \ll d(C)\sum_{\mathcal J}\sum_{m\in\mathcal M_{\mathcal J}}
       \kappa^k f_\kappa(m)
       \sum_{p'_j\in J_j\cap\mathbb P}A_C(P'm).
\]
To count how often an integer occurs in this full sum, put $n=P'm$.
Every marked prime exceeds $C^\alpha$, and $n<2\lambda^kC<C^2$ for
large $C$. If $n$ has $r$ distinct prime factors exceeding $C^\alpha$,
their product divides $n$, so $C^{\alpha r}<C^2$ and
$r\le\lfloor2/\alpha\rfloor$. The marked product $P'$ is determined
by a choice of $k$ of these distinct primes, their increasing order
determines the tuple, and then $m=n/P'$ is determined. Each tuple
belongs to exactly one box. Thus each $n$ has at most
$\binom{\lfloor2/\alpha\rfloor}{k}=O_{k,\alpha}(1)$ representations
in the sum.

The nonzero terms have squarefree $m$. When $(P',m)=1$, the integer
$n=P'm$ is squarefree and
$\kappa^k f_\kappa(m)A_C(P'm)=f_\kappa(n)A_C(n)$.
When $(P',m)>1$, some marked prime $p>C^\alpha$ divides both $P'$ and
$m$, so $p^2\mid n$. Denote the total weight of these latter terms,
over all boxes and cofactors in the sum, by $E_{\rm sq}$.
The bounded number of representations gives
\[
 D_\kappa^{\mathrm{all}}
 \ll d(C)\left(
       \sum_{C/\lambda^k<n<2\lambda^kC}f_\kappa(n)A_C(n)
       +E_{\rm sq}\right).
\]

We next show that the terms divisible by the square of a marked prime
are negligible at this scale. For every $n\ge2$, including non-squarefree
$n$, \eqref{en:eq:convolution-identity} and $|\xi_e|\le1$ give
\[
 0\le A_C(n)
 \le\sum_{e\mid n}|\xi_e|\frac{\log(n/e)}{\log n}
 \le\tau_2(n),
\]
where $\tau_2$ is the divisor function. Since $m\mid n$, each
formal weight is at most
$\kappa^k\kappa^{\omega(n)}\tau_2(n)\ll_{k,\kappa,\varepsilon}C^\varepsilon$.

Combining this weight bound with the bounded number of representations
of each integer, and putting $c=2\lambda^k$, we obtain
\begin{equation}\label{en:eq:marker-collision-bound}
 \begin{aligned}
 E_{\rm sq}
 &\ll_\varepsilon C^\varepsilon
       \sum_{p>C^\alpha}\#\{n\le cC:p^2\mid n\}\\
 &\le cC^{1+\varepsilon}\sum_{p>C^\alpha}p^{-2}
 \ll C^{1-\alpha+\varepsilon}.
 \end{aligned}
\end{equation}
For any fixed $0<\varepsilon<\alpha/2$,
\[
 \frac{E_{\rm sq}}{C/L}\ll L C^{-\alpha+\varepsilon}=o(1).
\]
The factor $d(C)$ has already been retained outside the complete sum,
so the negligible contribution is $d(C)E_{\rm sq}=o(d(C)C/L)$.

Finally, applying \eqref{en:unif:base-positive} to the fixed interval
$C/\lambda^k<n<2\lambda^kC$ gives the uniform bound
\begin{equation}\label{en:unif:colored-positive}
 D_\kappa^{\mathrm{all}}
 \ll d(C)\bigl(C/L+E_{\rm sq}\bigr)
 \ll d(C)C/L,\qquad\kappa=2,4.
\end{equation}
The implied constant depends only on the fixed parameters and the
uniform constants in the prime restrictions, not on $C$ or the particular
admissible sets $\mathcal P_{j,C}$.

This prepares the principal-character estimate used below. The principal
character $\chi_0$ modulo $M$ only deletes the nonnegative terms with
$(n,M)>1$, so its twisted positive moment satisfies
\[
 0\le\mathcal P_\kappa(\chi_0)
 :=\sum_n f_\kappa(n)g(n/C)A_C(n)H_C(n)
                       \mathbf1_{(n,M)=1}
 \le D_\kappa^{\mathrm{all}}.
\]
After character orthogonality, its contribution to one reduced residue
class is therefore at most
$D_\kappa^{\mathrm{all}}/\varphi(M)\ll d(C)C/(\varphi(M)L)$.

We now turn to the nonprincipal characters and first establish two
preliminary estimates.

\paragraph{Nonprincipal characters with coprimality conditions.}
Let the nonprincipal character $\chi\pmod M$ be induced by a primitive
character $\chi^*$ of conductor $M_*>1$. For $q\ge1$, put
$\chi_q(m)=\chi(m)\mathbf1_{(m,q)=1}$. This is still completely
multiplicative, and
\[
 \chi_q(m)=\chi^*(m)\mathbf1_{(m,Mq)=1}
 =\chi^*(m)\sum_{d\mid(m,Mq)}\mu(d).
\]
For any real interval $I$, interchanging the finite sums and writing
$m=dv$ gives
\[
 \sum_{m\in I}\chi_q(m)
 =\sum_{d\mid Mq}\mu(d)\chi^*(d)
                     \sum_{v\in I/d}\chi^*(v).
\]
P\'olya--Vinogradov \cite{en:ik} bounds each inner interval sum by
$O(\sqrt{M_*}\log(2M_*))$. Using
$\tau_2(Mq)\ll_\varepsilon(Mq)^\varepsilon$ and absorbing the logarithm
into an arbitrarily small power, we obtain
\begin{equation}\label{en:eq:pv-coprimality}
 \sup_I\left|\sum_{m\in I}\chi_q(m)\right|
 \ll_\varepsilon M^{1/2+\varepsilon}q^\varepsilon.
\end{equation}

\paragraph{The multidimensional hyperbola argument and the Euler remainder.}
The second preliminary estimate is
\eqref{en:unif:nonprincipal-arithmetic} below, for the character sum with
weight $\mu^2(m)\kappa^{\omega(m)}$. We first estimate a related
convolution sum and then turn to the desired weight by comparing Euler
factors. Fix an integer $\kappa\ge2$, and put $a_\kappa=1-1/\kappa$.
Let $h_{\kappa,q}$ be the $\kappa$-fold Dirichlet convolution of $\chi_q$;
explicitly,
\[
 h_{\kappa,q}(n)
 =\sum_{m_1\cdots m_\kappa=n}\prod_{i=1}^\kappa\chi_q(m_i).
\]
Its partial sum is therefore
\[
 H_{\kappa,q}(X)
 :=\sum_{n\le X}h_{\kappa,q}(n)
 =\sum_{m_1\cdots m_\kappa\le X}
                  \prod_{i=1}^\kappa\chi_q(m_i),\qquad X\ge1.
\]

Partition the tuples in this sum into $\kappa$ disjoint classes by
assigning each tuple to the first position $j$ at which its largest entry
occurs. The conditions for assignment to $j$ are precisely
\[
 m_i<m_j\quad(i<j),\qquad m_i\le m_j\quad(i>j).
\]
Thus every tuple is
assigned to exactly one position, even when its maximum occurs several
times.

For a fixed $j$, hold the other $\kappa-1$ coordinates fixed and set
\[
 v_j=\prod_{i\ne j}m_i,\qquad
 \ell_j=\max\bigl(\{1\}\cup\{m_i+1:i<j\}
                              \cup\{m_i:i>j\}\bigr).
\]
The inequalities defining the assigned position give the lower bound
$m_j\ge\ell_j$, and the product condition gives $m_j\le X/v_j$.
Consequently, the possible values of this last variable are exactly the
integers in
\[
 I_j=\left[\ell_j,\left\lfloor\frac{X}{v_j}\right\rfloor\right]
                      \cap\mathbb Z.
\]
The inner sum is understood to be zero if this interval is empty.
Whenever $I_j$ is nonempty, choose $m_j\in I_j$ and write
$N=m_1\cdots m_\kappa\le X$. Since $m_j$ is a largest coordinate,
$m_j\ge N^{1/\kappa}$, and hence
\[
 v_j=\frac{N}{m_j}
 \le N^{1-1/\kappa}\le X^{a_\kappa}.
\]
We may therefore restrict the product of the fixed coordinates to this
range. Combining the disjoint assignment with the interval description
gives the exact decomposition
\[
 H_{\kappa,q}(X)
 =\sum_{j=1}^\kappa
   \sum_{\substack{(m_i)_{i\ne j}\in\mathbb N^{\kappa-1}\\
                   v_j\le X^{a_\kappa}}}
      \left(\prod_{i\ne j}\chi_q(m_i)\right)
      \sum_{\ell_j\le t\le X/v_j}\chi_q(t).
\]
The condition $v_j\le X^{a_\kappa}$ is only necessary for
a nonempty interval; any extra choices in the outer sum contribute zero.

For fixed $j$ and $v$, the number of choices of the remaining coordinates
with product $v$ is $\tau_{\kappa-1}(v)$. Some of these choices
may give an empty interval, but counting all of them is valid for an
upper bound. Since $|\chi_q(m_i)|\le1$, the triangle inequality and \eqref{en:eq:pv-coprimality} yield
\[
 |H_{\kappa,q}(X)|
 \ll_\varepsilon
 \kappa M^{1/2+\varepsilon}q^\varepsilon
       \sum_{v\le X^{a_\kappa}}\tau_{\kappa-1}(v).
\]

Let $Y\ge1$. For an integer $s\ge2$, first fix the first $s-1$ factors and count the
possible last factor:
\begin{align*}
 \sum_{v\le Y}\tau_s(v)
 &=\sum_{u_1\cdots u_{s-1}\le Y}
       \left\lfloor\frac{Y}{u_1\cdots u_{s-1}}\right\rfloor\\
 &\le Y\sum_{u_1\cdots u_{s-1}\le Y}
                   \frac1{u_1\cdots u_{s-1}}\\
 &\le Y\left(\sum_{1\le u\le Y}\frac1u\right)^{s-1}
 \le Y(1+\log Y)^{s-1}.
\end{align*}
For $s=1$, the corresponding bound
is simply $\sum_{v\le Y}\tau_1(v)=\lfloor Y\rfloor\le Y$.
Taking $s=\kappa-1$ and $Y=X^{a_\kappa}$ therefore gives
\[
 \sum_{v\le X^{a_\kappa}}\tau_{\kappa-1}(v)
 \ll_\kappa X^{a_\kappa}(\log(2X))^{\kappa-2}.
\]
Substituting this count into $H_{\kappa,q}$,
and absorbing the fixed number $\kappa$ of assigned positions into the
implied constant, proves
\begin{equation}\label{en:eq:hyperbola-character}
 |H_{\kappa,q}(X)|
 \ll_{\kappa,\varepsilon}
 M^{1/2+\varepsilon}q^\varepsilon
 X^{a_\kappa}(\log(2X))^{\kappa-2}.
\end{equation}
In particular, when $\kappa=2$ there is only one fixed coordinate,
$\tau_1(v)=1$, and this step incurs no log factor.

From $h_{\kappa,q}$ to
$\mu^2(m)\kappa^{\omega(m)}\chi_q(m)$, we compare their Dirichlet series.
For $\Re s>1$, comparison of Euler factors yields
\begin{equation}\label{en:eq:arithmetic-euler-factor}
 \begin{aligned}
 \sum_{m\ge1}\frac{\mu^2(m)\kappa^{\omega(m)}\chi_q(m)}{m^s}
 &=L(s,\chi_q)^\kappa G_{\kappa,q}(s),\\
 L(s,\chi_q)&=\prod_p(1-\chi_q(p)p^{-s})^{-1}.
 \end{aligned}
\end{equation}
Write $G_{\kappa,q}(s)=\sum_r g_{\kappa,q}(r)r^{-s}$. Its local polynomial is
\[
 (1+\kappa u)(1-u)^\kappa
 =1+\sum_{\nu=2}^{\kappa+1}c_{\kappa,\nu}u^\nu,
 \qquad u=\chi_q(p)p^{-s}.
\]
For every fixed $b>1/2$,
\begin{equation}\label{en:eq:euler-absolute-convergence}
 \sum_{r\ge1}\frac{|g_{\kappa,q}(r)|}{r^b}
 \le\prod_p\left(1+\sum_{\nu=2}^{\kappa+1}
                     |c_{\kappa,\nu}|p^{-\nu b}\right)
 \ll_{\kappa,b}1,
\end{equation}
because $\sum_p p^{-2b}<\infty$. The bound is uniform in $M,q,\chi$.
The coefficients of $L(s,\chi_q)^\kappa$ are $h_{\kappa,q}$.
Thus comparison of coefficients in \eqref{en:eq:arithmetic-euler-factor}
gives the convolution identity
\[
 \mu^2(m)\kappa^{\omega(m)}\chi_q(m)
 =\sum_{r\mid m}g_{\kappa,q}(r)h_{\kappa,q}(m/r).
\]
Summing over $m\le X$, writing $m=ru$, and using the definition of
$H_{\kappa,q}$ gives
\begin{align*}
 \sum_{m\le X}\mu^2(m)\kappa^{\omega(m)}\chi_q(m)
 &=\sum_{r\le X}g_{\kappa,q}(r)
                     \sum_{u\le X/r}h_{\kappa,q}(u)\\
 &=\sum_{r\le X}g_{\kappa,q}(r)H_{\kappa,q}(X/r).
\end{align*}
Given $\varepsilon>0$, apply \eqref{en:eq:hyperbola-character} with
$\varepsilon/2$ and use
$(\log(2Y))^{\kappa-2}\ll_{\kappa,\varepsilon}Y^{\varepsilon/2}$
for $Y\ge1$. The last identity then implies
\begin{align*}
 \left|\sum_{m\le X}\mu^2(m)\kappa^{\omega(m)}\chi_q(m)\right|
 &\ll_{\kappa,\varepsilon}
 M^{1/2+\varepsilon/2}q^{\varepsilon/2}X^{a_\kappa+\varepsilon/2}
 \sum_{r\le X}\frac{|g_{\kappa,q}(r)|}{r^{a_\kappa+\varepsilon/2}}\\
 &\ll_{\kappa,\varepsilon}
 M^{1/2+\varepsilon/2}q^{\varepsilon/2}X^{a_\kappa+\varepsilon/2}.
\end{align*}
For the final step, use \eqref{en:eq:euler-absolute-convergence} at
$b=a_\kappa+\varepsilon/2>1/2$. Enlarging the positive powers
from $\varepsilon/2$ to $\varepsilon$ proves the second preliminary estimate:
\begin{equation}\label{en:unif:nonprincipal-arithmetic}
 \sum_{m\le X}\mu^2(m)\kappa^{\omega(m)}\chi(m)\mathbf1_{(m,q)=1}
 \ll_{\kappa,\varepsilon}
 X^{1-1/\kappa+\varepsilon}M^{1/2+\varepsilon}q^\varepsilon.
\end{equation}
For $\kappa=2$, absolute convergence cannot be asserted at $b=1/2$ itself;
the arbitrarily small positive power above pays for precisely this endpoint.
This character argument in fact applies to every fixed integer $\kappa\ge2$;
the same derivation is available when $7^{\omega(n)}$ is used below.

\paragraph{Summing all divisibility levels after expansion.}
We now apply \eqref{en:unif:nonprincipal-arithmetic} to the complete
twisted positive moment
\[
 \mathcal P_\kappa(\chi)
 :=\sum_n f_\kappa(n)g(n/C)A_C(n)H_C(n)\chi(n),
 \qquad\kappa=2,4.
\]
Our goal in this paragraph is to prove, for every nonprincipal character
$\chi$ modulo $M\le C^\delta$ and every $\varepsilon>0$,
\[
 |\mathcal P_\kappa(\chi)|
 \ll_{\kappa,\varepsilon}
 C^{1-(1-\eta-\sigma)/\kappa+\delta/2+\varepsilon}.
\]
This will control the nonprincipal terms when we recover the residue
classes in the next paragraph.

Recall the coefficients $w_C(P)$ from the expansion in
\eqref{en:eq:signed-exact-expansion}, so that
$H_C(n)=\sum_{P\mid n}w_C(P)$ and $0\le w_C(P)\le1$.
Substituting \eqref{en:eq:convolution-identity} and expanding $H_C$ gives
\[
 \mathcal P_\kappa(\chi)
 =\sum_{e,P}\xi_e w_C(P)
   \sum_{\substack{eP\mid n\\C<n<2C}}
      f_\kappa(n)g(n/C)
      \frac{\log(n/e)}{\log n}\chi(n).
\]
For every nonzero term, $e$ and $P$ are squarefree and coprime, since
the sieve primes are at most $z$ and the marked primes exceed $z$.
Put $q=eP$ and $n=qm$. Then
\[
 f_\kappa(qm)
 =\mu^2(q)\kappa^{\omega(q)}
       f_\kappa(m)\mathbf1_{(m,q)=1},
 \qquad \chi(qm)=\chi(q)\chi(m).
\]
Writing
\[
 X_q=C/q,\qquad
 \psi_{e,C}(t)=g(t)\left(1-\frac{\log e}{L+\log t}\right),
\]
we have $qm/C=m/X_q$ and
$\log(qm/e)/\log(qm)=1-\log e/(L+\log(m/X_q))$.
The test $\psi_{e,C}$ is defined on $(1,2)$ and extended by zero outside
that interval. With $\chi_q(m)=\chi(m)\mathbf1_{(m,q)=1}$, the expansion
becomes the exact identity
\[
 \mathcal P_\kappa(\chi)
 =\sum_{e,P}\xi_e w_C(P)\mu^2(q)\kappa^{\omega(q)}\chi(q)
       I_{e,q},
 \qquad
 I_{e,q}:=\sum_m f_\kappa(m)\chi_q(m)
                          \psi_{e,C}(m/X_q).
\]
The product of the prime factors of $q$ at most $z$ recovers $e$, and
the remaining prime factors, in increasing order, recover $P$ and its
marked tuple. Thus each $q$ occurs at most once, and
$q\le C^{\eta+\sigma}$. When $k=0$, take $P=1$ and $\sigma=0$;
the same identities and uniqueness statement hold.

We next estimate the inner sum $I_{e,q}$. Since $e\le C^\eta$,
$0\le\log e\le\eta L$. Differentiation on the support of $g$ gives
\[
 \psi_{e,C}'(t)
 =g'(t)\left(1-\frac{\log e}{L+\log t}\right)
       +\frac{g(t)\log e}{t(L+\log t)^2}.
\]
Since $L+\log t\asymp L$ on this fixed support,
\[
 \|\psi_{e,C}\|_\infty+
          \int_1^2|\psi_{e,C}'(t)|\,dt\ll_{g,\eta}1
\]
uniformly in $e$ and $C$. Define the partial sum
\[
 S_q(u)=\sum_{m\le u}f_\kappa(m)\chi_q(m).
\]
The test vanishes near both endpoints of $(1,2)$, so summation by parts
followed by the change of variable $u=X_qt$ gives
\[
 I_{e,q}=-\int_1^2S_q(X_qt)\psi_{e,C}'(t)\,dt.
\]
Use an auxiliary number $\nu>0$ in
\eqref{en:unif:nonprincipal-arithmetic}. Since
$X_q\ge C^{1-\eta-\sigma}\to\infty$, that estimate and the derivative
bound imply
\begin{align*}
 |I_{e,q}|
 &\le\sup_{X_q\le u\le2X_q}|S_q(u)|
                          \int_1^2|\psi_{e,C}'(t)|\,dt\\
 &\ll_{\kappa,\nu}
       (C/q)^{a_\kappa+\nu}M^{1/2+\nu}q^\nu,
 \qquad a_\kappa=1-1/\kappa.
\end{align*}

The coefficient $\xi_e w_C(P)\chi(q)$ has absolute value at most $1$,
whereas the remaining factor satisfies
$\mu^2(q)\kappa^{\omega(q)}\le\kappa^{\omega(q)}\ll_{\kappa,\nu}q^\nu$.
Reindexing the complete outer sum by its distinct values $q=eP$ therefore
gives
\begin{align*}
 |\mathcal P_\kappa(\chi)|
 &\ll_{\kappa,\nu}
 M^{1/2+\nu}\sum_{q\le C^{\eta+\sigma}}
       q^\nu(C/q)^{a_\kappa+\nu}q^\nu\\
 &\ll_{\kappa,\nu} C^{a_\kappa+\nu}M^{1/2+\nu}
           \sum_{q\le C^{\eta+\sigma}}q^{-a_\kappa+\nu}\\
 &\ll_{\kappa,\nu}
 C^{a_\kappa+\nu}M^{1/2+\nu}
       (C^{\eta+\sigma})^{1-a_\kappa+\nu}.
\end{align*}
Using $M\le C^\delta$, the exponent of $C$ in the last expression is
\[
 a_\kappa+(\eta+\sigma)(1-a_\kappa)+\frac\delta2
       +\nu(1+\eta+\sigma+\delta)
 =1-\frac{1-\eta-\sigma}{\kappa}+\frac\delta2
       +\nu(1+\eta+\sigma+\delta).
\]
For the prescribed $\varepsilon>0$, choose
$0<\nu\le\varepsilon/(1+\eta+\sigma+\delta)$. We obtain the desired
bound, uniformly over the nonprincipal characters:
\begin{equation}\label{en:unif:positive-character-error}
 |\mathcal P_\kappa(\chi)|
 \ll_{\kappa,\varepsilon}
 C^{1-(1-\eta-\sigma)/\kappa+\delta/2+\varepsilon},
 \qquad\chi\ne\chi_0,\quad\kappa=2,4.
\end{equation}

\paragraph{Recovery of the class set and completion of the proposition.}
We return to the upper bound that motivated this subsection:
\[
 D_C=\sum_n\mu^2(n)g(n/C)2^{\omega(n)}A_C(n)H_C(n)
                  \mathbf1_{n\bmod M\in\mathcal R_C}
 \ll\Delta_C,
 \qquad
 \Delta_C=\frac{d(C)|\mathcal R_C|C}{\varphi(M)L}.
\]
We prove the bound simultaneously for $D_{2,C}=D_C$ and $D_{4,C}$.
For a reduced residue class $a$ modulo $M$, write
\[
 D_{\kappa,a}
 :=\sum_{n\equiv a\ ({\rm mod}\ M)}
           f_\kappa(n)g(n/C)A_C(n)H_C(n),
 \qquad
 D_{\kappa,C}=\sum_{a\in\mathcal R_C}D_{\kappa,a}.
\]
Character orthogonality separates the principal and nonprincipal terms:
\[
 D_{\kappa,a}
 =\frac{\mathcal P_\kappa(\chi_0)}{\varphi(M)}
  +\frac1{\varphi(M)}
     \sum_{\substack{\chi\ ({\rm mod}\ M)\\\chi\ne\chi_0}}
                  \overline{\chi(a)}\mathcal P_\kappa(\chi).
\]
For the principal character, all summands are nonnegative and
$\chi_0(n)=\mathbf1_{(n,M)=1}$. Thus the estimate retaining the prime
restrictions, \eqref{en:unif:colored-positive}, gives
\[
 0\le\frac{\mathcal P_\kappa(\chi_0)}{\varphi(M)}
 \le\frac{D_\kappa^{\mathrm{all}}}{\varphi(M)}
 \ll\frac{d(C)C}{\varphi(M)L}.
\]
For the nonprincipal characters, apply
\eqref{en:unif:positive-character-error}. Since $|\chi(a)|=1$ and there
are $\varphi(M)-1$ such characters, the triangle inequality yields
\begin{align*}
 \left|\frac1{\varphi(M)}
     \sum_{\substack{\chi\ ({\rm mod}\ M)\\\chi\ne\chi_0}}
                  \overline{\chi(a)}\mathcal P_\kappa(\chi)\right|
 &\le\frac1{\varphi(M)}
     \sum_{\substack{\chi\ ({\rm mod}\ M)\\\chi\ne\chi_0}}
                  |\mathcal P_\kappa(\chi)|\\
 &\ll C^{1-(1-\eta-\sigma)/\kappa+\delta/2+\varepsilon}.
\end{align*}
The factor $(\varphi(M)-1)/\varphi(M)$ is at most $1$; if there are
no nonprincipal characters, this contribution is zero. Combining the
two bounds and then summing over the disjoint classes in $\mathcal R_C$
gives
\begin{equation}\label{en:unif:positive-union}
 D_{\kappa,C}\ll\Delta_C+
 |\mathcal R_C|C^{1-(1-\eta-\sigma)/\kappa+\delta/2+\varepsilon}.
\end{equation}
Here $\Delta_C=d(C)|\mathcal R_C|C/(\varphi(M)L)$.

It remains to compare the second term with $\Delta_C$. Using
$\varphi(M)\le M\le C^\delta$ and $d(C)\ge L^{-M_0}$, we have
\begin{align*}
 \frac{|\mathcal R_C|
       C^{1-(1-\eta-\sigma)/\kappa+\delta/2+\varepsilon}}{\Delta_C}
 &=\frac{\varphi(M)L}{d(C)}
       C^{-(1-\eta-\sigma)/\kappa+\delta/2+\varepsilon}\\
 &\le C^{-(1-\eta-\sigma)/\kappa+3\delta/2+\varepsilon}
                 L^{M_0+1}.
\end{align*}
For $\kappa=2,4$, the largest exponent of $C$ occurs at $\kappa=4$.
Choose
\[
 \delta<\frac{1-\eta-\sigma}{12},\qquad
 0<\varepsilon<\frac{1-\eta-\sigma}{16}.
\]
Then
\[
 -\frac{1-\eta-\sigma}{\kappa}+\frac{3\delta}{2}+\varepsilon
 \le-\frac{1-\eta-\sigma}{4}+\frac{3\delta}{2}+\varepsilon
 <-\frac{1-\eta-\sigma}{16}.
\]
Consequently, the nonprincipal contribution in
\eqref{en:unif:positive-union} is $o(\Delta_C)$, and we have proved
\[
 D_C=D_{2,C}\ll\Delta_C,\qquad D_{4,C}\ll\Delta_C.
\]

To complete Proposition~\ref{en:thm:unified}, we also use the lower bound
$B_C\gg\Delta_C$ from Lemma~\ref{en:lem:unified-absolute}, which will
be proved independently in the following section. Let $\delta_{\rm abs}>0$ denote an
allowable exponent for $M\le C^{\delta_{\rm abs}}$ supplied by that
proof. Its construction there includes
\[
 \delta_{\rm abs}<\frac1{100},\qquad
 \delta_{\rm abs}<\frac\gamma2,\qquad
 \delta_{\rm abs}<\frac\alpha2\quad(k\ge1),
\]
where $\gamma$ is defined later in \eqref{en:unif-abs:gamma}.
The same construction makes $\delta_{\rm abs}$ sufficiently small for
the bounds on normalized prime character sums in
\eqref{en:unif-abs:prime-bias}, including the conductor range required
by the uniform prime number theorem used there. These restrictions on
$\delta_{\rm abs}$ are obtained in that later proof.

The choices can therefore be made in the following order. First fix the
sieve and marking parameters. The argument leading to
\eqref{en:unif:signed-union} then supplies $\delta_0>0$,
and the proof of Lemma~\ref{en:lem:unified-absolute} supplies
$\delta_{\rm abs}>0$ using those fixed parameters. Choose
\[
 0<\delta<\min\left\{\frac{\delta_0}2,
        \frac{1-\eta-\sigma}{12},\delta_{\rm abs}\right\},
\]
then choose $\varepsilon$ in the range above, and finally choose $\nu$
as required for \eqref{en:unif:positive-character-error}.
Lemma~\ref{en:lem:unified-absolute} gives $B_C\gg\Delta_C$ for this
choice of $\delta$. The squarefree Weil bound gives $B_C\le D_C$,
so the upper bounds just proved imply
$B_C\asymp D_C\asymp\Delta_C$. Together with
\eqref{en:unif:signed-union}, which gives $T_C=o(\Delta_C)$, and the
bound for $D_{4,C}$, this proves Proposition~\ref{en:thm:unified}.
Finally, $\mathcal R_C\ne\varnothing$, $d(C)\ge L^{-M_0}$, and
$\varphi(M)\le M\le C^\delta$ give
\[
 \Delta_C\ge C^{1-\delta}L^{-M_0-1}\longrightarrow\infty.
\]

\section{Lower bound of \texorpdfstring{$B_a(C)$}{B\_a(C)}}
\label{en:unif:absolute-proof}

\begin{lemma}[Absolute mass in each reduced residue class]\label{en:lem:unified-absolute}
Under the hypotheses of the general marked construction, there is a sufficiently
small fixed $\delta>0$ such that, uniformly for $M\le C^\delta$ and $(a,M)=1$,
\[
 B_a(C):=\sum_{n\equiv a\ ({\rm mod}\ M)}
 \mu^2(n)g(n/C)|K(n)|A_C(n)H_C(n)
 \gg \frac{d(C)C}{\varphi(M)L}.
\]
The implied constant may depend on the fixed sieve parameters, marking regions,
prime-counting constants, and $g$, but not on $a,M,C$.
Summing over any prescribed set $\mathcal R_C$ of reduced residue classes gives
$B_C\gg\Delta_C$.
\end{lemma}

\subsection{Exceptional real zeros and products of three or four primes}
To obtain the lower bound for $B_a(C)$, we will construct integers
$n=Pq_1\cdots q_r$, where $P$ is a fixed marked product and free prime indicator
$r\in\{3,4\}$. When $(P,M)=1$, the condition $n\equiv a\pmod M$
is equivalent to
\[
 q_1\cdots q_r\equiv b\pmod M,\qquad b=aP^{-1}\pmod M.
\]
For $k=0$, take $P=1$. The construction in the next subsection ensures
$(P,M)=1$ when $k\ge1$. Our first task is to count sufficiently many
prime tuples satisfying this congruence. Later parts of the proof will
show that a fixed positive proportion of the resulting candidates have
$|K(n)|$ bounded away from zero.

Fix $\lambda>1$ sufficiently close to $1$, and consider intervals
$J=[Y,\lambda Y)$ with $C^{0.12}\le Y\ll C$. Write
\[
 N_J=\#(J\cap\mathbb P),\qquad
 a_J(\chi)=\frac1{N_J}\sum_{p\in J}\chi(p),\qquad
 \rho_J(u)=\frac1{N_J}\#\{p\in J:p\equiv u\pmod M\}.
\]
Here $N_J>0$ for sufficiently large $C$. Initially require
$M\le C^\delta$ with $\delta<0.01$; we will decrease $\delta$ below.
Every prime in these intervals then exceeds $M$, so it is coprime to
$M$. In particular, $a_J(\chi_0)=1$ and
$\sum_{u\in(\mathbb Z/M\mathbb Z)^\times}\rho_J(u)=1$.

For $r\in\{3,4\}$ and intervals $J_1,\ldots,J_r$ in this range, put
$N=\prod_{i=1}^rN_{J_i}$. For a reduced class $b$ modulo $M$, define
the proportion
\[
 R_r(b)=\frac1N
 \#\left\{(q_1,\ldots,q_r)\in\prod_{i=1}^r(J_i\cap\mathbb P):
                   q_1\cdots q_r\equiv b\pmod M\right\}.
\]
We aim to choose $r=3$ or $r=4$ so that, uniformly in the reduced target
class $b$ and the permitted intervals,
\[
 R_r(b)\ge\frac1{2\varphi(M)}.
\]
Character orthogonality gives
\begin{align*}
 R_r(b)
 &=\frac1{N\varphi(M)}
   \sum_{\chi\ ({\rm mod}\ M)}\overline{\chi(b)}
          \prod_{i=1}^r\sum_{q_i\in J_i\cap\mathbb P}\chi(q_i)\\
 &=\frac1{\varphi(M)}
   \left(1+\sum_{\substack{\chi\ ({\rm mod}\ M)\\\chi\ne\chi_0}}
           \overline{\chi(b)}\prod_{i=1}^r a_{J_i}(\chi)\right).
\end{align*}
Thus we need to control the nonprincipal part. We will use a square-sum estimate over all characters
together with the individual estimates in
\eqref{en:unif-abs:prime-bias}, treating a possible exceptional real
character separately.

First, the Brun--Titchmarsh inequality
\cite[Theorem 2, (1.10)]{en:bt} and the prime number theorem give
\[
 \rho_J(u)\le\frac{\pi(\lambda Y;M,u)}{N_J}
 \ll_\lambda\frac1{\varphi(M)}.
\]
Indeed, $\log(\lambda Y/M)\asymp\log Y$ in the stated ranges and
$N_J\asymp_\lambda Y/\log Y$. Choose a fixed $K_\lambda\ge1$ such
that $\rho_J(u)\le K_\lambda/\varphi(M)$ uniformly in $C,M,J,u$.
Since $a_J(\chi)=\sum_u\rho_J(u)\chi(u)$, another application of
character orthogonality gives
\begin{equation}\label{en:unif-abs:l2}
 \sum_{\chi\ ({\rm mod}\ M)}|a_J(\chi)|^2
 =\varphi(M)\sum_u\rho_J(u)^2
 \le\varphi(M)\max_u\rho_J(u)\sum_u\rho_J(u)
 \le K_\lambda.
\end{equation}
The sums over $u$ are over reduced residue classes. This bound is
independent of the number $\varphi(M)$ of characters.

We next establish the estimates needed in $R_r(b)$. Given any fixed $0<\varepsilon_1<1$, we may
choose $\delta>0$ sufficiently small so that, for each admissible $M$,
there is at most one exceptional nonprincipal character $\chi_*\pmod M$
and
\begin{equation}\label{en:unif-abs:prime-bias}
 |a_J(\chi)|\le\varepsilon_1\quad(\chi\ne\chi_0,\chi_*),
 \qquad a_J(\chi_*)=-b_J+O_\lambda(\varepsilon_1),
 \quad 0\le b_J\le1.
\end{equation}
The exceptional character, when present, is induced by a real primitive
quadratic character; otherwise its term is omitted. The same character
$\chi_*$ is used for all intervals $J$ at a given modulus $M$.
We prove \eqref{en:unif-abs:prime-bias} without assuming that this
exceptional character is absent.

Let $\chi'$ be the nonprincipal primitive character inducing $\chi$,
with conductor $q\mid M$. The $\mathrm{GL}_1$ case of
Kaneko--Thorner \cite[Theorem 2.4]{en:kt} supplies absolute constants
$A,c>0$ such that, for $t\ge(4q)^A$,
\[
 \sum_{n\le t}\chi'(n)\Lambda(n)
 =-\frac{t^{\beta_{\chi'}}}{\beta_{\chi'}}
 +O\!\left(t\exp\left\{-\frac{c\log t}
                 {\log(4q)+\sqrt{\log t}}\right\}\right).
\]
Here $\beta_{\chi'}$ is the largest real zero exceeding $3/4$, with
the convention $\beta_{\chi'}=1/2$ if there is no such zero. The smaller
error factor in the cited theorem has been bounded above by $t$.
To ensure the conductor range throughout our intervals, require
$A\delta<0.06$. Since $q\le M\le C^\delta$ and $Y\ge C^{0.12}$,
we then have $(4q)^A\le4^AC^{A\delta}\le Y$ for all sufficiently
large $C$.

Put $\vartheta_{\chi'}(t)=\sum_{p\le t}\chi'(p)\log p$.
Removing the terms supported on higher prime powers gives
\[
 \left|\sum_{n\le t}\chi'(n)\Lambda(n)-\vartheta_{\chi'}(t)\right|
 \le\sum_{h\ge2}\sum_{p\le t^{1/h}}\log p
 \ll\sqrt t(\log(2t))^2.
\]
For this bound, separate $h=2$ and bound the at most $\log_2t$
remaining exponents trivially. Thus, uniformly for $Y\le t\le\lambda Y$,
\[
 \vartheta_{\chi'}(t)=-t^{\beta_{\chi'}}/\beta_{\chi'}+E(t),
 \qquad
 \sup_{Y\le t\le\lambda Y}|E(t)|
 \ll_\lambda Y\mathcal E(Y,M)+\sqrt Y(\log Y)^2,
\]
where
\[
 \mathcal E(Y,M)=
 \exp\left\{-\frac{c\log Y}{\log(4M)+\sqrt{\log Y}}\right\}.
\]
This follows by using $q\le M$, $t\asymp_\lambda Y$, and decreasing
the absolute constant $c$ if necessary, so that one bound holds on the
whole interval.

Partial summation removes the log weight:
\[
 \sum_{p\in J}\chi'(p)
 =\int_Y^{\lambda Y}\frac{d\vartheta_{\chi'}(t)}{\log t}+O(1).
\]
The $O(1)$ term comes from the half-open endpoints. The derivative of
the main term is $-t^{\beta_{\chi'}-1}$, while the error contributes
\[
 \left[\frac{E(t)}{\log t}\right]_Y^{\lambda Y}
       +\int_Y^{\lambda Y}\frac{E(t)}{t(\log t)^2}\,dt
 \ll_\lambda\frac{\sup_{Y\le t\le\lambda Y}|E(t)|}{\log Y}.
\]
Consequently,
\begin{align*}
 \sum_{p\in J}\chi(p)
 &=-\int_Y^{\lambda Y}\frac{t^{\beta_{\chi'}-1}}{\log t}\,dt\\
 &\quad+O_\lambda\!\left(
       \frac Y{\log Y}\mathcal E(Y,M)+\sqrt Y\log Y\right).
\end{align*}
Here $\chi$ and its inducing character $\chi'$ agree at every prime in
$J$, because those primes exceed $M$. Dividing by
$N_J\asymp_\lambda Y/\log Y$ gives
\[
 a_J(\chi)
 =-\frac1{N_J}\int_Y^{\lambda Y}
                \frac{t^{\beta_{\chi'}-1}}{\log t}\,dt
   +O_\lambda\!\left(\mathcal E(Y,M)+Y^{-1/2}(\log Y)^2\right).
\]

The zero-free region in \cite[Theorem 12.3]{en:koukoulopoulos} ensures
that a primitive character with a real zero larger than
$1-c_0/\log(4M)$ must be real, for a sufficiently small absolute
$c_0>0$. The Landau--Page theorem \cite{en:page}, in the form of
\cite[Corollary 12.7]{en:koukoulopoulos}, permits at most one such
primitive character of conductor at most $M$. If its conductor divides $M$, denote its
induced character modulo $M$ by $\chi_*$; otherwise no exceptional
character occurs among the characters being summed. For every remaining
nonprincipal character, the normalized integral is bounded by
\[
 \frac1{N_J}\int_Y^{\lambda Y}
                \frac{t^{\beta_{\chi'}-1}}{\log t}\,dt
 \ll_\lambda Y^{\beta_{\chi'}-1}
 \le\exp\left\{-\frac{c_0\log Y}{\log(4M)}\right\},
\]
after decreasing $c_0$ if needed. The convention
$\beta_{\chi'}=1/2$ is also covered by this upper bound.
Since $M\le C^\delta$ and $C^{0.12}\le Y\ll C$, the two ratios in the
exponentials satisfy
\[
 \frac{\log Y}{\log(4M)+\sqrt{\log Y}}
 \ge\frac{0.12L}{\delta L+O(1)+O(\sqrt L)},
 \qquad
 \frac{\log Y}{\log(4M)}
 \ge\frac{0.12L}{\delta L+O(1)}.
\]
Choosing the fixed $\delta$ sufficiently small, and then $C$ sufficiently
large, makes both exponential upper bounds as small as required in
terms of $\varepsilon_1$. The term $Y^{-1/2}(\log Y)^2$ also tends
to zero uniformly in the permitted range. This proves
$|a_J(\chi)|\le\varepsilon_1$ for $\chi\ne\chi_0,\chi_*$.

For the possible exceptional character, the normalized integral is
nonnegative and at most $1+o(1)$: use $\beta_{\chi'}\le1$ and
$N_J\sim\int_Y^{\lambda Y}dt/\log t$. Define
\[
 b_J=\min\left\{1,\frac1{N_J}\int_Y^{\lambda Y}
                   \frac{t^{\beta_{\chi'}-1}}{\log t}\,dt\right\}.
\]
This changes the normalized integral by $o(1)$ and ensures $0\le b_J\le1$.
The same error estimate therefore gives
$a_J(\chi_*)=-b_J+O_\lambda(\varepsilon_1)$, completing the proof of
\eqref{en:unif-abs:prime-bias}.

We now insert these estimates into the character expansion of $R_r(b)$.
For all characters other than $\chi_0$ and the possible $\chi_*$, use
\eqref{en:unif-abs:prime-bias} on $r-2$ factors and then Cauchy's
inequality on the remaining two:
\begin{align*}
 \sum_{\chi\ne\chi_0,\chi_*}
          \left|\prod_{i=1}^r a_{J_i}(\chi)\right|
 &\le\varepsilon_1^{r-2}
          \sum_{\chi\ ({\rm mod}\ M)}
                |a_{J_1}(\chi)a_{J_2}(\chi)|\\
 &\le\varepsilon_1^{r-2}
       \left(\sum_\chi|a_{J_1}(\chi)|^2\right)^{1/2}
       \left(\sum_\chi|a_{J_2}(\chi)|^2\right)^{1/2}\\
 &\le K_\lambda\varepsilon_1^{r-2},
\end{align*}
where the last line uses \eqref{en:unif-abs:l2}.

For the exceptional real character, \eqref{en:unif-abs:prime-bias} gives
\[
 \prod_{i=1}^r a_{J_i}(\chi_*)
 =(-1)^r\prod_{i=1}^r b_{J_i}+O_\lambda(r\varepsilon_1).
\]
Indeed, a telescoping expansion of the difference of these two products
has $r$ terms, and all other factors have absolute value at most $1$.
Also $\overline{\chi_*(b)}=\chi_*(b)\in\{-1,1\}$. Combining this term,
the principal term, and the bound for the remaining characters proves
the desired proportion formula
\begin{equation}\label{en:unif-abs:parity}
 R_r(b)=\frac1{\varphi(M)}
 \left(1+\chi_*(b)(-1)^r\prod_{j=1}^r b_{J_j}
       +O_\lambda(K_\lambda\varepsilon_1^{r-2}+r\varepsilon_1)\right).
\end{equation}
The exceptional term is omitted when there is no $\chi_*$.

Finally choose the number of free primes as follows. If there is no
exceptional character, take $r=3$. Otherwise take $r=3$ when
$\chi_*(b)=-1$, and $r=4$ when $\chi_*(b)=1$. In both exceptional
cases, $\chi_*(b)(-1)^r=1$, so the additional main term in
\eqref{en:unif-abs:parity} is nonnegative. First choose
$\varepsilon_1>0$ sufficiently small that the absolute value of the
error is at most $1/2$ for both $r=3$ and $r=4$.
Then choose $\delta$ as in the proof of \eqref{en:unif-abs:prime-bias},
and take $C$ sufficiently large. We obtain the required lower bound
\[
 R_r(b)\ge\frac1{2\varphi(M)},\qquad
 \#\{\text{prime tuples in the target class}\}
       \ge\frac{N}{2\varphi(M)}.
\]
This choice of $r$ is uniform over the prime intervals for fixed $M,b$.

\begin{remark}[Effect of the new zero-free results]\label{en:rem:new-zero-input}
The preprint of OpenAI\cite[Theorem 1]{en:oaisiegel26} gives the uniform real-zero
gap $(1-\beta)\log q\ge c>0$. This alone extends the small upper bound
for $|a_J(\chi)|$ to the possible exceptional character; after decreasing
$\delta$, one can use $r=3$ throughout. However, for $q\asymp C^\delta$
and $Y\asymp C^a$, the resulting bound
$Y^{\beta-1}\le\exp(-ca/\delta+o(1))$ is only a fixed small quantity and need
not tend to zero. The other preprint\cite[Theorem 1.1]{en:oaizero26} gives the
stronger zero-free half-plane $\Re s>7/8$ for every Dirichlet $L$-function.
Section~\ref{en:num:main} derives an upper bound for $|a_J(\chi)|$ that
tends to zero uniformly over the nonprincipal characters, using that
input, and extends the explicit family of moduli.
\end{remark}

\subsection{Prime intervals for a fixed marked tuple}
\label{en:unif-abs:prime-boxes}
We now apply \eqref{en:unif-abs:parity} to the integers in
$B_a(C)$. The number $k$ of marked primes remains fixed throughout the
proof. For each fixed marked tuple, we will find
$\gg C/(P\varphi(M)L)$ candidates in the target class. The next two
subsections will show that a fixed positive proportion of these
candidates have $|K(n)|$ bounded away from zero.

Choose $1<u<v<2$ and $g_0>0$ such that $g\ge g_0$ on $[u,v]$.
Use a fixed $\lambda>1$ sufficiently close to $1$ that
$u\lambda^4<v$ and \eqref{en:unif:prime-colors} applies. We use
the former subsection with this choice of $\lambda$, and choose
$\delta$ accordingly.

For $k\ge1$, fix $(p_1,\ldots,p_k)$ with
$p_j\in\mathcal P_{j,C}$ and
$(\log p_1/L,\ldots,\log p_k/L)$ in the open box $Q$, and set
\[
 P=p_1\cdots p_k<C^\sigma,\qquad
 X=C/P>C^{1-\sigma}>C^{1/2}.
\]
For $k=0$, use $P=1$ and $X=C$. When $k\ge1$, also require
$\delta<\alpha/2$. Then $p_j>C^\alpha>M$, so $(P,M)=1$.
Set $b=aP^{-1}\pmod M$, and choose $r\in\{3,4\}$ by the rule following
\eqref{en:unif-abs:parity}. This choice is fixed once the marked tuple
and target class are fixed.

We next choose the intervals for the free primes. For $r=3$, take
\[
 \frac{\log Q_1}{\log X}\in(0.320,0.322),\qquad
 \frac{\log Q_2}{\log X}\in(0.330,0.332),\qquad
 Q_3=\frac{uX}{Q_1Q_2}.
\]
For $r=4$, take
\[
 \frac{\log Q_1}{\log X}\in(0.240,0.241),\quad
 \frac{\log Q_2}{\log X}\in(0.247,0.248),\quad
 \frac{\log Q_3}{\log X}\in(0.253,0.254),\quad
 Q_4=\frac{uX}{Q_1Q_2Q_3}.
\]
In each case, put $J_i=[Q_i,\lambda Q_i)$. For the first $r-1$
coordinates, take $Q_i$ from the grid $\{\lambda^m:m\in\mathbb Z\}$
and retain only intervals lying wholly inside the stated exponent
interval. The last interval is determined by the first $r-1$; its lower
endpoint need not be a grid point. Its exponent lies in
$(0.346,0.350)+O(1/L)$ when $r=3$, and in
$(0.257,0.260)+O(1/L)$ when $r=4$.

These exponent ranges are separated. Hence, for sufficiently large
$C$, the intervals are ordered and disjoint, and their free primes
are distinct. Their lower endpoints exceed $X^{0.240}$, and thus
exceed $C^{0.12}$ because $X>C^{1/2}$. All free primes therefore
exceed both $M$ and $z$. The ranges also give
\[
 Q_1\cdots Q_r=uX,\qquad
 uC\le Pq_1\cdots q_r<u\lambda^rC<vC
 \quad(q_i\in J_i).
\]
These exponent intervals will also make the local moment bounds in
Section~\ref{en:unif-abs:conditional} small enough for the next step.

We count these boxes before removing any primes. If $(s_i,t_i)$ is
one of the first $r-1$ exponent intervals, the number of complete grid
intervals in that range is
\[
 \frac{(t_i-s_i)\log X}{\log\lambda}+O(1)\asymp L.
\]
Here $(1-\sigma)L\le\log X\le L$. There are therefore
$\asymp L^{r-1}$ boxes. In each box, continue to write
$N=\prod_{i=1}^rN_{J_i}$ for the number of free tuples before deletion.
The prime number theorem gives, uniformly over the marked tuples and
these boxes,
\[
 N_{J_i}=(1+o(1))\frac{(\lambda-1)Q_i}{\log Q_i}
             \asymp\frac{Q_i}{L},\qquad
 N\asymp\frac{X}{L^r}.
\]
It follows that
\[
 \sum_{\text{free-prime boxes}}N\asymp\frac XL.
\]
Different boxes cannot contain the same free tuple: each of its first
$r-1$ coordinates belongs to exactly one half-open grid interval,
which then determines the last interval.

We must also ensure that no free prime equals a marked prime. In each
box, define
\[
 \begin{split}
 \mathcal V&=\prod_{i=1}^r
       \bigl((J_i\cap\mathbb P)\setminus\{p_1,\ldots,p_k\}\bigr),
       \qquad n(\mathbf q)=P\prod_{i=1}^r q_i,\\
 \mathcal V_a&=\{\mathbf q\in\mathcal V:
                        n(\mathbf q)\equiv a\pmod M\},
       \qquad N_a=|\mathcal V_a|.
 \end{split}
\]
For $k=0$, the set of marked primes being removed is empty. We will
use this notation in the following subsections. In particular, $N$
always denotes the size of the original box, before deletion.

Equation \eqref{en:unif-abs:parity}, with the chosen $r$, gives at
least $N/(2\varphi(M))$ tuples in the target class.
At most $k$ primes are removed from each coordinate. The total number
of tuples removed, including those outside the target class, is at
most
\[
 \sum_{i=1}^r k\prod_{h\ne i}N_{J_h}
 =N\sum_{i=1}^r\frac{k}{N_{J_i}}
 \ll_k N\frac L{Q_{\min}},\qquad Q_{\min}=\min_iQ_i.
\]
Since $Q_{\min}>C^{.12}$ and $\varphi(M)\le M\le C^\delta$ with
$\delta<.01$, we have $\varphi(M)L/Q_{\min}=o(1)$ uniformly.
Thus $|\mathcal V|=(1+o(1))N$, and, for sufficiently large $C$,
\[
 N_a\ge\frac{N}{2\varphi(M)}-O_k\!\left(\frac{NL}{Q_{\min}}\right)
       \ge\frac{N}{3\varphi(M)}.
\]

For every $\mathbf q\in\mathcal V_a$, the integer
$n=n(\mathbf q)$ has exactly $k+r$ distinct prime factors, all
larger than $z$. Consequently,
\[
 \mu^2(n)=1,\qquad A_C(n)=1,\qquad g(n/C)\ge g_0.
\]
For the fixed marked tuple, summing over the free-prime boxes gives
\[
 \sum_{\text{free-prime boxes}}N_a
 \gg\frac X{\varphi(M)L}
 =\frac C{P\varphi(M)L}.
\]
These tuples give distinct integers when the marked tuple is fixed:
unique factorization and the ordered free-prime intervals determine
each $q_i$ from $n/P$. This is the required candidate count. It remains
to retain a fixed positive proportion for which $|K(n)|$ is not small.

\subsection{A bilinear bound for the local Kloosterman factors}
\label{en:unif-abs:fourier}
Section~\ref{en:unif-abs:prime-boxes} provides enough prime tuples in
the our required residue class. We still need to show that many of them have
large local Kloosterman factors. The next subsection will do this using
Chebyshev moments. Here we prove the bilinear bound
\eqref{en:unif-abs:bilinear} needed to control those moments, including
when the variables are larger than the local prime modulus.

All prime factors in our construction tend to infinity, so we may take
the local modulus $\ell$ to be odd. For $t\in\mathbb F_\ell^\times$,
the group of nonzero residue classes modulo $\ell$, define
\[
 \kappa_\ell(t)=\ell^{-1/2}S(1,t^{-2};\ell),\qquad
 F_{\ell,j}(t)=U_j(\kappa_\ell(t)/2),\quad j\ge1.
\]
Here $U_j$ is the Chebyshev polynomial of the second kind. The Weil
bound gives $\kappa_\ell(t)\in[-2,2]$. And we have the estimation
\begin{equation}\label{en:unif-abs:mellin}
 \max_{\psi\ ({\rm mod}\ \ell)}
 \left|\sum_{t\in\mathbb F_\ell^\times}
       F_{\ell,j}(ct)\overline{\psi(t)}\right|
 \ll_j\sqrt\ell,\qquad c\in\mathbb F_\ell^\times.
\end{equation}
The maximum includes the principal character, and the constant is
uniform in $\ell$, $c$, and $\psi$ for each fixed $j$.

We recall why the geometric results apply to this particular function.
Katz's monodromy theorem \cite[Theorem 11.1]{en:katz} gives the group
$\mathrm{SL}_2$ for the rank-two Kloosterman sheaf. After normalization,
pulling this sheaf back by $t\mapsto t^{-2}$ gives trace
$-\kappa_\ell(t)$. Its $j$th symmetric power has trace
$(-1)^jU_j(\kappa_\ell(t)/2)$; the fixed sign does not affect the bound.
The inverse-square map is a finite \'etale covering when $\ell$ is odd.
Its geometric fundamental-group image has finite index, so the
connected monodromy group remains $\mathrm{SL}_2$. The symmetric power
is therefore geometrically irreducible of rank $j+1>1$.

Multiplication of the trace by $\overline{\psi(t)}$ corresponds to a
rank-one Kummer twist. This preserves irreducibility, so the resulting
sheaf has no constant rank-one component. Its conductor is $O_j(1)$,
uniformly in the character and the scaling $c$. These facts are the
reason for the square-root bound: in the trace formula, the degree-two
term vanishes, while the degree-one term has only $O_j(1)$ eigenvalues.
More precisely, $H_c^2=0$ and the Grothendieck--Ogg--Shafarevich formula
gives $\dim H_c^1=O_j(1)$. Deligne's weight theorem
\cite[Theorem 3.3.1]{en:deligne} bounds each normalized eigenvalue by
$\sqrt\ell$, proving \eqref{en:unif-abs:mellin}. The same
inverse-power pullback and its conductor bound are treated in
\cite[\S5.2]{en:fkmprime}.

We change from the complete sum to sums over integers. Let
$\mathcal A\subset[A,c_1A]$ and $\mathcal B\subset[B,c_1B]$, where
$A,B\ge1$ and $c_1>1$ is fixed. Assume all members are coprime to
$\ell$, and let $|\alpha_s|,|\beta_t|\le1$. Put
\[
 \mathcal S=\sum_{s\in\mathcal A}\sum_{t\in\mathcal B}
               \alpha_s\beta_tF_{\ell,j}(cst).
\]
For a multiplicative character $\psi$ modulo $\ell$, define
\[
 \widehat F_c(\psi)=\sum_{u\in\mathbb F_\ell^\times}
                 F_{\ell,j}(cu)\overline{\psi(u)},\qquad
 A(\psi)=\sum_{s\in\mathcal A}\alpha_s\psi(s),\qquad
 B(\psi)=\sum_{t\in\mathcal B}\beta_t\psi(t).
\]
Character orthogonality gives the inversion formula
\[
 F_{\ell,j}(cu)=\frac1{\ell-1}
                     \sum_\psi\widehat F_c(\psi)\psi(u).
\]
Substituting $u=st$ and using $\psi(st)=\psi(s)\psi(t)$ separates
the two variables:
\[
 \mathcal S=\frac1{\ell-1}
                 \sum_\psi\widehat F_c(\psi)A(\psi)B(\psi).
\]

The sets $\mathcal A$ and $\mathcal B$ consist of actual integers.
Several of them may have the same residue modulo $\ell$. For each
nonzero residue $v$, let
\[
 a_v=\sum_{\substack{s\in\mathcal A\\s\equiv v\ ({\rm mod}\ \ell)}}
                      \alpha_s,\qquad
 r_v=\#\{s\in\mathcal A:s\equiv v\pmod\ell\}.
\]
The interval containing $\mathcal A$ has length $(c_1-1)A$, so
\[
 r_v\le1+\frac{(c_1-1)A}{\ell}\ll_{c_1}1+A/\ell.
\]
Orthogonality first sums over the residue classes. Cauchy's inequality
then bounds the sum within each class:
\begin{align*}
 \sum_\psi|A(\psi)|^2
 &=(\ell-1)\sum_{v\in\mathbb F_\ell^\times}|a_v|^2\\
 &\le(\ell-1)\sum_v r_v
             \sum_{\substack{s\in\mathcal A\\s\equiv v\ ({\rm mod}\ \ell)}}
                         |\alpha_s|^2\\
 &\ll_{c_1}(\ell-1)(1+A/\ell)|\mathcal A|.
\end{align*}
The same argument gives
\[
 \sum_\psi|B(\psi)|^2
 \ll_{c_1}(\ell-1)(1+B/\ell)|\mathcal B|.
\]
Using \eqref{en:unif-abs:mellin} in $\mathcal S$
and applying Cauchy's inequality, we obtain
\[
 |\mathcal S|
 \ll_j\frac{\sqrt\ell}{\ell-1}
       \left(\sum_\psi|A(\psi)|^2\right)^{1/2}
       \left(\sum_\psi|B(\psi)|^2\right)^{1/2}.
\]
Each square root supplies a factor $(\ell-1)^{1/2}$. Together these
cancel the denominator from inversion. We have proved
\begin{equation}\label{en:unif-abs:bilinear}
 \left|\sum_{s\in\mathcal A,t\in\mathcal B}
        \alpha_s\beta_tF_{\ell,j}(cst)\right|
 \ll_{j,c_1}\sqrt{\ell(1+A/\ell)(1+B/\ell)
                         |\mathcal A|\,|\mathcal B|}.
\end{equation}
The factors $1+A/\ell$ and $1+B/\ell$ count the 
integers in each residue class. Thus the estimate remains valid when
one or both intervals extend beyond $\ell$.

In Section~\ref{en:unif-abs:conditional}, each variable will be a prime
or a product of primes from disjoint ordered intervals. Unique
factorization makes the latter representation unique, so these products
form sets of integers to which \eqref{en:unif-abs:bilinear} applies.
Any character modulo $M$ can be included in $\alpha_s$ and $\beta_t$;
their absolute values remain at most $1$.

\subsection{Large local factors in the residue class}
\label{en:unif-abs:conditional}
Fix a marked tuple and one free-prime box from
Section~\ref{en:unif-abs:prime-boxes}. We use the sets
$\mathcal V$, $\mathcal V_a$ and the counts $N$, $N_a$ defined there.
Our goal is to obtain at least half of $\mathcal V_a$ while making
every local factor of $K(n(\mathbf q))$ larger than a fixed positive
constant in absolute value. The constant will be independent of the
marked tuple, the box, $M$, and the class $a$.

Every $n=n(\mathbf q)$ with $\mathbf q\in\mathcal V$ is squarefree.
The Chinese remainder theorem gives:
\[
 K(n)=\prod_{\ell\mid n}\kappa_\ell(n/\ell).
\]
Each cofactor $n/\ell$ is coprime to $\ell$. Thus there are
$k+r\le k+4$ factors, each in $[-2,2]$.
We first use the bilinear estimate from
Section~\ref{en:unif-abs:fourier} to control the Chebyshev moments of
each one, including after restriction to the residue class.

\paragraph{A marked prime as the local modulus.}
Suppose $k\ge1$ and fix $\ell=p_i$. In a three-prime box, group the
free variables as $q_1$ and $q_2q_3$, and set
$A=Q_1$, $B=Q_2Q_3$. In a four-prime box, group them as
$q_1q_2$ and $q_3q_4$, and set $A=Q_1Q_2$, $B=Q_3Q_4$.
Let $\mathcal A$ and $\mathcal B$ be the corresponding sets of
integer products after the deletions already made in $\mathcal V$.
If the first group contains $h$ primes, then
\[
 |\mathcal A|\asymp\frac{A}{L^h},\qquad
 |\mathcal B|\asymp\frac{B}{L^{r-h}},\qquad
 |\mathcal A||\mathcal B|=|\mathcal V|\asymp N.
\]
These comparisons from the prime counts in
Section~\ref{en:unif-abs:prime-boxes} and the unique grouped representations
noted at the end of Section~\ref{en:unif-abs:fourier}. Also
$AB=uX\asymp X$, and the chosen exponent intervals give
$A,B\ge X^{0.24}$.
Each group has at most two primes, so its products lie in
$[A,\lambda^2A]$ or $[B,\lambda^2B]$. We may therefore use
$c_1=\lambda^2$ in \eqref{en:unif-abs:bilinear}.

For a fixed $j\ge1$ and any character $\chi$ modulo $M$, put
\[
 S_{\ell,j}(\chi)
 :=\sum_{\mathbf q\in\mathcal V}\chi(n(\mathbf q))
       U_j\!\left(\frac{\kappa_\ell(n(\mathbf q)/\ell)}2\right)
 =\chi(P)\sum_{s\in\mathcal A,\ t\in\mathcal B}
       \chi(s)\chi(t)F_{\ell,j}((P/\ell)st).
\]
All the variables and $P/\ell$ are coprime to $\ell$.
The coefficients have absolute value at most $1$, so
\eqref{en:unif-abs:bilinear} gives
\begin{align*}
 \frac{|S_{\ell,j}(\chi)|}{|\mathcal A||\mathcal B|}
 &\ll_j L^{r/2}
       \sqrt{\frac{\ell(1+A/\ell)(1+B/\ell)}{AB}}\\
 &\ll_j L^{r/2}
       \sqrt{\ell/X+1/A+1/B+1/\ell}.
\end{align*}
Since $C^\alpha<\ell<C^\sigma$ and
$X>C^{1-\sigma}$, the square roots of the four terms are at most
a fixed constant times
\[
 C^{-(1-2\sigma)/2},\qquad
 C^{-0.12(1-\sigma)},\qquad
 C^{-0.12(1-\sigma)},\qquad C^{-\alpha/2},
\]
respectively. Each upper bound tends to zero as a fixed negative
power of $C$. The choice of $\gamma$ below also absorbs $L^{r/2}$.

\paragraph{A free prime as the local modulus.}
Now consider the position occupied by $\ell=q_i$. First fix the
actual prime value of $q_i$. For $r=3$, vary the other two primes.
For $r=4$, choose one other position once, fix the prime value at that
position as well, and vary the remaining two primes. Denote these fixed values by
$\mathbf v$, their corresponding slice of $\mathcal V$ by
$\mathcal V_{\mathbf v}$, and its size by $N_{\mathbf v}$.
The two varying primes range independently over their intervals after we drooped
the marked primes. Their interval parameters will be denoted
by $A$ and $B$.

Write the sum on a slice as
\[
 S_{\mathbf v,j}(\chi)
 :=\sum_{\mathbf q\in\mathcal V_{\mathbf v}}\chi(n(\mathbf q))
       U_j\!\left(\frac{\kappa_\ell(n(\mathbf q)/\ell)}2\right).
\]
On this slice, $n/\ell$ is a nonzero fixed constant times the product
of the two varying primes. The fixed part of $\chi(n)$ has absolute value
at most $1$, and the remaining character values enter the two
coefficients. We may therefore apply
\eqref{en:unif-abs:bilinear} again. Write
$A=X^a$, $B=X^b$, and $\ell=X^t$. Since the two varying prime sets
have sizes $\asymp A/L$ and $\asymp B/L$, division by
$N_{\mathbf v}$ gives the bound
\[
 \frac{|S_{\mathbf v,j}(\chi)|}{N_{\mathbf v}}
 \ll_j L\left(
 X^{-(a+b-t)/2}+X^{-a/2}+X^{-b/2}+X^{-t/2}\right).
\]
The exponent intervals in Section~\ref{en:unif-abs:prime-boxes}
give, uniformly in the fixed values,
\[
 a+b-t\ge
 \begin{cases}
 0.320+0.330-0.350-O(1/L)=0.300-O(1/L),&r=3,\\
 0.240+0.247-0.260-O(1/L)=0.227-O(1/L),&r=4,
 \end{cases}
 \qquad a,b,t\ge 0.240-O(1/L).
\]
Thus the upper bound is $O_j(LX^{-0.11})$. As $X>C^{1/2}$,
it is also $O_j(C^{-1/40})$. In particular, for every slice the
unnormalized sum is bounded by a fixed constant times
$N_{\mathbf v}C^{-1/40}$.

Finally sum over $\mathbf v$. These are fixed prime values, not
residue classes modulo $M$; distinct prime values are kept distinct
even when they have the same residue. The slices partition
$\mathcal V$, so
\[
 \sum_{\mathbf v}N_{\mathbf v}=|\mathcal V|\le N.
\]
Consequently, summing the slice estimates gives
$O_j(NC^{-1/40})$ for the full box.

\paragraph{Restriction to the residue class.}
For any one marked or free position, write
$Z(\mathbf q)=\kappa_\ell(n(\mathbf q)/\ell)$, where $\ell$ is
the prime at that position. At a free position, it varies with
$\mathbf q$. The two arguments above show that we may take
\begin{equation}\label{en:unif-abs:gamma}
 \gamma=\begin{cases}
 1/40,&k=0,\\
 \min\{(1/2-\sigma)/2,\alpha/4,1/40\},&k\ge1.
 \end{cases}
\end{equation}
Indeed, in the marked case the first and fourth exponents above are
at least $2\gamma$, and $0.12(1-\sigma)>0.06>1/40$; hence all fixed
logarithmic factors can be absorbed. For each fixed $j\ge1$ there
is a constant $C_j$ such that
\[
 \left|\sum_{\mathbf q\in\mathcal V}\chi(n(\mathbf q))
                   U_j(Z(\mathbf q)/2)\right|
 \le C_jNC^{-\gamma}
\]
for every character $\chi$ modulo $M$ and every local position.
The constant is independent of the marked tuple, the box, $M$, and
$\chi$. The reference count is the original box size $N$.

Character orthogonality now gives
\[
 \sum_{\mathbf q\in\mathcal V_a}U_j(Z(\mathbf q)/2)
 =\frac1{\varphi(M)}\sum_{\chi\ ({\rm mod}\ M)}
       \overline{\chi(a)}
       \sum_{\mathbf q\in\mathcal V}
             \chi(n(\mathbf q))U_j(Z(\mathbf q)/2).
\]
There are $\varphi(M)$ terms, each with a coefficient of absolute
value $1/\varphi(M)$, so the absolute value of this sum is at most
$C_jNC^{-\gamma}$. Dividing by the bound
$N_a\ge N/(3\varphi(M))$ from
Section~\ref{en:unif-abs:prime-boxes}, we obtain
\[
 \left|\frac1{N_a}\sum_{\mathbf q\in\mathcal V_a}
                   U_j(Z(\mathbf q)/2)\right|
 \le3C_j\varphi(M)C^{-\gamma}
 \le3C_jC^{\delta-\gamma}=o(1),
\]
after requiring $\delta<\gamma/2$. It also applies to a free
position because the slice argument already summed over every value
of its varying modulus.

\paragraph{All local factors large at once.}
We now turn these moment estimates into a count. The measure
$\nu_{\rm ST}$ is defined in the notation. Put
$m_0=k+4$. Choose a sufficiently small fixed $0<b_0<1$ and a
continuous function $f:[-2,2]\to[0,1]$ such that $f=1$ on
$[-b_0,b_0]$ and $f=0$ when $|t|\ge2b_0$. Since the density of
$\nu_{\rm ST}$ is at most $1/\pi$,
\[
 \int f\,d\nu_{\rm ST}\le\frac{4b_0}{\pi}.
\]
Choose a fixed $\varepsilon_0>0$ such that
$4b_0/\pi+2\varepsilon_0<1/(4m_0)$. By the Weierstrass
approximation theorem, there is a real polynomial $Q_0$ satisfying
$\|Q_0-f\|_\infty\le\varepsilon_0$ on $[-2,2]$.
The polynomial $R=Q_0+\varepsilon_0$ then satisfies
\[
 R\ge f\ge\mathbf1_{[-b_0,b_0]},\qquad
 R\ge0\quad\hbox{on }[-2,2],\qquad
 \int R\,d\nu_{\rm ST}<\frac1{4m_0}.
\]
Here the last inequality follows from $R\le f+2\varepsilon_0$.

Expand this fixed polynomial in the Chebyshev basis:
\[
 R(t)=a_0+\sum_{j=1}^{J_0}a_jU_j(t/2).
\]
For $j\ge1$, the substitution $t=2\cos\theta$ gives
\[
 \int U_j(t/2)\,d\nu_{\rm ST}(t)
 =\frac2\pi\int_0^\pi
            \sin((j+1)\theta)\sin\theta\,d\theta=0.
\]
The constant polynomial has integral $1$, so
$a_0=\int R\,d\nu_{\rm ST}<1/(4m_0)$.
Let $\mathbb P_a$ and $\mathbb E_a$ denote proportion and average
on the finite set $\mathcal V_a$. For every local position, the
majorant and the conditional moment bounds give
\[
 \mathbb P_a(|Z|\le b_0)
 \le\mathbb E_a R(Z)
 \le a_0+3\varphi(M)C^{-\gamma}
                    \sum_{j=1}^{J_0}|a_j|C_j.
\]
The degree, coefficients, and constants in this finite sum are all
fixed. Its second term is therefore eventually
smaller than $1/(4m_0)$, uniformly in all the choices under
consideration. Each event $|Z|\le b_0$ consequently occupies less
than $1/(2m_0)$ of $\mathcal V_a$.

There are at most $m_0$ positions. The union bound shows that fewer
than half the tuples have a small factor at any position. Hence at
least $N_a/2\ge N/(6\varphi(M))$ tuples satisfy
$|\kappa_\ell(n/\ell)|>b_0$ for every $\ell\mid n$. For these tuples,
\[
 |K(n)|=\prod_{\ell\mid n}|\kappa_\ell(n/\ell)|
       >b_0^{k+r}\ge b_0^{k+4}.
\]
This proves the desired conclusion in each fixed box and
class. The next subsection sums these contributions over the marked tuples.

\subsection{Summing over the marked primes}
\label{en:unif-abs:marked-sum}
We can now complete the lower bound
\[
 B_a(C)\gg\frac{d(C)C}{\varphi(M)L}.
\]
For each fixed marked tuple and each free-prime box,
Section~\ref{en:unif-abs:conditional} retains at least $N_a/2$ tuples
in $\mathcal V_a$ for which $|K(n(\mathbf q))|\ge b_0^{k+4}$.
The box count in Section~\ref{en:unif-abs:prime-boxes} therefore gives
at least
\[
 c\frac C{P\varphi(M)L}
\]
such tuples after summing over all free-prime boxes, where $c>0$ is
independent of the marked tuple, $M$, and $a$. All resulting integers
are squarefree and satisfy $A_C(n)=1$ and $g(n/C)\ge g_0$.
The chosen marked tuple lies in the given box, so $H_C(n)\ge c_F$.
When $k=0$, we use the convention $c_F=1$, since $H_C=F=1$.
Thus each distinct integer obtained in this way contributes at least
\[
 g_0c_Fb_0^{k+4}
\]
to $B_a(C)$.

For $k\ge1$, we next sum the factor $1/P$ over the marked tuples.
This step uses a grid for the marked primes, separate from the grid
for the free primes in Section~\ref{en:unif-abs:prime-boxes}. In
coordinate $j$ of the existing box $Q$, retain the complete
intervals $[Y,\lambda Y)$ from the same half-open grid that lie in
$(C^{a_j},C^{b_j})$. On each such interval,
\[
 \frac{\#(\mathcal P_{j,C}\cap[Y,\lambda Y))}{\lambda Y}
 \le\sum_{\substack{p\in\mathcal P_{j,C}\\Y\le p<\lambda Y}}\frac1p
 \le\frac{\#(\mathcal P_{j,C}\cap[Y,\lambda Y))}{Y}.
\]
The prime-counting hypothesis
\eqref{en:unif:prime-colors} makes this sum
$\asymp d_j(C)/L$. There are $\asymp L$ complete intervals in each
fixed exponent interval. We drop the at most two incomplete
intervals at its endpoints; they are not needed for a lower bound.
If $\mathcal W_j$ denotes the primes of $\mathcal P_{j,C}$ in the
retained intervals, it follows that
\[
 \sum_{p\in\mathcal W_j}\frac1p\asymp d_j(C).
\]
The coordinate intervals of $Q$ are ordered and disjoint, so every
choice $p_j\in\mathcal W_j$ is an admissible marked tuple. Hence
\[
 \sum_{p_1\in\mathcal W_1,\ldots,p_k\in\mathcal W_k}
          \frac1{p_1\cdots p_k}
 =\prod_{j=1}^k\sum_{p\in\mathcal W_j}\frac1p
 \asymp\prod_{j=1}^k d_j(C)=d(C).
\]
The number of retained construction representations is consequently
at least
\[
 c\frac C{\varphi(M)L}
   \sum_{p_1\in\mathcal W_1,\ldots,p_k\in\mathcal W_k}\frac1P
 \gg\frac{d(C)C}{\varphi(M)L}.
\]

It remains to bound how often one integer is represented. The number
$k$ was fixed before the construction, and every retained integer
$n$ has $m=k+r$ distinct prime factors, with $r=3$ or $4$. Thus
$m=\omega(n)$ determines $r=m-k$; the two possible values of $r$ do
not give an extra choice for a fixed $n$. Every representation assigns
the $m$ prime factors of $n$ to the ordered positions
\[
 (p_1,\ldots,p_k,q_1,\ldots,q_r).
\]
There are at most $m!\le(k+4)!$ such assignments. Each assignment
determines $P$ and $X=C/P$. It also determines at most one free-prime
box, since each of the first $r-1$ free primes lies in a unique
half-open grid interval and the last interval is fixed by those
intervals and $X$. The grid used to compute the marked-prime sum
also assigns each marked prime to a unique interval; it does not
create additional representations. Therefore each integer has at
most $(k+4)!$ construction representations.

Dividing the representation count by $(k+4)!$ and multiplying by
the fixed contribution $g_0c_Fb_0^{k+4}$ proves the asserted lower
bound for $B_a(C)$. For $k=0$, there is just the empty marked tuple,
so $\sum 1/P=1=d(C)$, and the same counting and
representation argument applies. Finally, summing over
$a\in\mathcal R_C$ gives $B_C\gg\Delta_C$ and completes the proof.
\hfill$\square$

\section{Two signs of fixed amplitude and small Hecke eigenvalues}
\label{en:app:main}
\subsection{Proof of the main theorem and consequences of the weighted estimates}
\label{en:app:main-proof}

To prove Theorem~\ref{en:thm:main-hecke}, we first verify the prime-counting
assumption \eqref{en:unif:prime-colors} for the sets in
\eqref{en:def:hecke-colors}, with $d_j(C)=m_i(C)$ for $j\in J_i$.
Proposition~\ref{en:thm:unified} will then give the weighted estimates
\eqref{en:eq:hecke-main}. Hecke multiplicativity will give
\eqref{en:eq:normalized-hecke-small} on the squarefree integers with
$H_C(n)>0$.

\paragraph{The effective prime count.}
For a fixed non-CM holomorphic newform $f$ of trivial nebentypus, the effective
Sato--Tate estimate \cite[Theorem 1.1]{en:thornerST} is, uniformly over all
intervals $I\subseteq[0,\pi]$,
\begin{equation}\label{en:eq:effective-st}
 \#\{p\le X:p\nmid N_f,\ \vartheta_f(p)\in I\}
 =\mu_{\rm ST}(I)\pi(X)
 +O_f\!\left(\pi(X)\frac{\log(w_fN_f\log X)}{\sqrt{\log X}}\right).
\end{equation}
One may also obtain this estimate for any fixed level from
Hoey--Iskander--Jin--Trejos Su\'arez
\cite[Remark 2.1 and Theorem 3.1]{en:hijtsST}.
The bound for the arithmetic conductor recalled in their Remark 2.1
gives a constant $Q_f\ge1$, depending only on $f$, such that
$q_{\operatorname{Sym}^m f}\le Q_f^{m+1}$ for every $m\ge1$.
Their Theorem 3.1 then gives
\[
 \sum_{\substack{p\le X\,,\ p\nmid N_f\\\vartheta_f(p)\in I}}\log p
 =\mu_{\rm ST}(I)X+
 O_f\!\left(X\frac{\log\log X}{\sqrt{\log X}}\right),
\]
uniformly in $I$. Partial summation divides the endpoint error by $\log X$
and adds the error integral
\[
 \int_3^X O_f\!\left(\frac{\log\log t}{(\log t)^{5/2}}\right)dt
 =O_f\!\left(\frac{X\log\log X}{(\log X)^{5/2}}\right).
\]
It follows that the unweighted prime count is
\[
 \mu_{\rm ST}(I)\operatorname{Li}(X)
 +O_f\!\left(\frac{X\log\log X}{(\log X)^{3/2}}\right).
\]
The prime number theorem, in the form
$\pi(X)-\operatorname{Li}(X)\ll X\exp(-c\sqrt{\log X})$, now gives
\eqref{en:eq:effective-st}. All of $w_f,N_f,Q_f$ are fixed, and no
uniformity over varying forms or levels is needed.

\paragraph{Verification of the prime restrictions.}
Put $\varepsilon_i=L^{-A_i/k_i}$. By the relation
$\lambda_{f_i}(p)=2\cos\vartheta_{f_i}(p)$, the condition in
\eqref{en:def:hecke-colors} corresponds to the interval
\[
 I_{i,C}=[\arccos(\varepsilon_i/2),\pi-\arccos(\varepsilon_i/2)].
\]
Its Sato--Tate mass is the quantity $m_i(C)$ defined in Section~\ref{en:unif:definition}.
The substitution $t=2\cos\vartheta$ and expansion at zero give
\begin{align*}
 m_i(C)&=\frac1{2\pi}\int_{-\varepsilon_i}^{\varepsilon_i}
                      \sqrt{4-t^2}\,dt
       =\frac1\pi\int_0^{\varepsilon_i}
                    \left(2-\frac{t^2}{4}+O(t^4)\right)dt\\
       &=\frac2\pi\varepsilon_i-\frac{\varepsilon_i^3}{12\pi}
                    +O(\varepsilon_i^5).
\end{align*}
This proves the first estimate in \eqref{en:eq:hecke-density}.

For a fixed $\lambda>1$ close to $1$, write
$N(Y)=\pi(\lambda Y)-\pi(Y)$ and $E_C=(\log L)/\sqrt L$.
In the range $C^\alpha\le Y\le C^\sigma$, we have $\log Y\asymp L$.
Applying \eqref{en:eq:effective-st} at $Y$ and $\lambda Y$ therefore gives
\[
 \#(\mathcal P_{j,C}\cap[Y,\lambda Y))
 =m_i(C)N(Y)+O(E_CY/L)\qquad(j\in J_i),
\]
uniformly in $Y$. The possible endpoint terms are $O(1)$ and are absorbed
by this error. Here the constant may depend on $f_i,\lambda,\alpha,\sigma$.
The prime number theorem gives
$N(Y)\sim(\lambda-1)Y/\log Y\asymp_\lambda Y/L$.
Thus the relative error is $O(E_C/m_i(C))$. Since $k_i>2A_i$,
\[
 \frac{m_i(C)}{E_C}
 \asymp\frac{L^{1/2-A_i/k_i}}{\log L}\longrightarrow\infty.
\]
Consequently, the prime-counting assumption \eqref{en:unif:prime-colors}
holds with $d_j(C)=m_i(C)$, and in fact the prime count has relative error
$o(1)$.

The numbers $k_i$ and their sum $k$ are fixed as in
Theorem~\ref{en:thm:main-hecke}. Multiplying the densities over its $k$
marked positions gives
\[
 d(C)=\prod_i\left\{\frac2\pi L^{-A_i/k_i}
                      (1+O(L^{-2A_i/k_i}))\right\}^{k_i}
       =\left((2/\pi)^k+o(1)\right)L^{-A}.
\]
This proves the second estimate in \eqref{en:eq:hecke-density}. It also gives
$d(C)\ge L^{-M_0}$ for any fixed $M_0>A$ and all sufficiently large $C$.
All the hypotheses of Proposition~\ref{en:thm:unified} are now satisfied,
so \eqref{en:unif:main} gives \eqref{en:eq:hecke-main}.

For squarefree $n$ with $H_C(n)>0$, choose a marked tuple contributing to
\eqref{en:unif:marker}. Its $k_i$ primes assigned to $f_i$ each contribute
at most $L^{-A_i/k_i}$ to $|\lambda_{f_i}(n)|$. Every other prime factor
contributes at most $2$, by the local bounds in the notation.
Hecke multiplicativity therefore gives
\[
 |\lambda_{f_i}(n)|\le2^{\omega(n)-k_i}
       (L^{-A_i/k_i})^{k_i}=2^{\omega(n)-k_i}L^{-A_i}.
\]
This proves \eqref{en:eq:normalized-hecke-small}
and completes the proof of Theorem~\ref{en:thm:main-hecke}.

The next corollary uses only \eqref{en:unif:main}. It separates the absolute
mass into the two signs.

\begin{corollary}[Sign masses and a bound on the number of prime factors]\label{en:cor:unified-applications}
Under the hypotheses of Proposition~\ref{en:thm:unified}, define
\[
 B_C^\pm=\sum_{\pm K(n)>0}\mu^2(n)g(n/C)|x_n|W_C(n).
\]
Then the following conclusions hold.
\begin{enumerate}[label=\textup{(\roman*)}]
\item The two signs have asymptotically equal absolute mass:
$B_C^\pm/B_C\to1/2$.
\item Uniformly for all integers $R\ge0$,
\begin{equation}\label{en:eq:unified-omega-tail}
 \sum_{\omega(n)>R}\mu^2(n)g(n/C)W_C(n)\ll2^{-R}\Delta_C.
\end{equation}
For every integer-valued function $R(C)\to\infty$, restricting
all three sums to $\omega(n)\le R(C)$ preserves
$D_C\asymp B_C\asymp\Delta_C$ and $T_C=o(\Delta_C)$.
\item There is a fixed integer $R_0$ such that each sign has at least
$\gg\Delta_C$ integer witnesses satisfying
\[
 C<n<2C,\quad\mu^2(n)=1,\quad\omega(n)\le R_0,\quad
 H_C(n)>0,\quad n\bmod M\in\mathcal R_C,\quad \pm K(n)>0.
\]
\end{enumerate}
The uniformity and dependence on fixed parameters are as in the proposition.
\end{corollary}
\begin{proof}
Write $w_n=\mu^2(n)g(n/C)W_C(n)\ge0$. Since $x_n$ has the same sign as
$K(n)$, separating its positive and negative values gives
\[
 B_C=B_C^++B_C^-,\qquad T_C=B_C^+-B_C^-,
 \qquad \frac{B_C^\pm}{B_C}=\frac12\pm\frac{T_C}{2B_C}.
\]
The last term tends to zero by $B_C\gg\Delta_C$ and $T_C=o(\Delta_C)$
in \eqref{en:unif:main}. This proves (i).

For (ii), on $\omega(n)>R$ we have $2^{\omega(n)}\le2^{-R}4^{\omega(n)}$.
All factors of the weight are nonnegative, so
\[
 E_R:=\sum_{\omega(n)>R}w_n
 \le2^{-R}\sum_n\mu^2(n)g(n/C)4^{\omega(n)}A_C(n)H_C(n)
                  \mathbf1_{n\bmod M\in\mathcal R_C}
 =2^{-R}D_{4,C}\ll2^{-R}\Delta_C,
\]
where the last bound is again \eqref{en:unif:main}.
Let $D_C^{(R)},T_C^{(R)},B_C^{(R)}$ denote the sums restricted to
$\omega(n)\le R$. Since $|x_n|\le1$ on squarefree integers,
\[
 0\le D_C-D_C^{(R)}\le E_R,\quad
 |T_C-T_C^{(R)}|\le E_R,\quad
 0\le B_C-B_C^{(R)}\le E_R.
\]
If $R=R(C)\to\infty$, then $E_R=o(\Delta_C)$, which proves the stated
estimates for the restricted sums.

For (iii), conclusion (i) and \eqref{en:unif:main} give a fixed $c>0$
such that $B_C^\pm\ge c\Delta_C$ for sufficiently large $C$.
Choose a fixed $R_0$ so that the constant in
\eqref{en:eq:unified-omega-tail} multiplied by $2^{-R_0}$ is at most $c/2$.
The remaining absolute mass of each sign is then at least $c\Delta_C/2$.

To turn this into a count, use the bound for $H_C$ stated after
\eqref{en:unif:marker}. More explicitly, for $C>2$, an integer $n<2C$
has at most $\lfloor2/\alpha\rfloor$ distinct prime factors exceeding
$C^\alpha$. Since the primes in each marked tuple are ordered and $0\le F\le1$,
\[
 H_C(n)\le H_0:=\binom{\lfloor2/\alpha\rfloor}{k}
 \quad(k\ge1),\qquad H_0=1\quad(k=0).
\]
Together with \eqref{en:eq:pointwise-weight}, this gives, for every
retained integer,
\[
 w_n|x_n|\le\|g\|_\infty H_0\,2^{\omega(n)}A_C(n)
             \le\|g\|_\infty H_0\,4^{R_0}.
\]
Dividing the remaining mass by this fixed bound proves (iii).
Positive weight implies all the listed support conditions.

The distinction between (ii) and (iii) is that $R_0$ in (iii) is fixed.
The comparison gives only
\[
 |T_C^{(R_0)}|\le |T_C|+E_{R_0}
       \le o(\Delta_C)+O(2^{-R_0}\Delta_C).
\]
This bound does not prove $T_C^{(R_0)}=o(\Delta_C)$. Conclusion (iii)
asserts only the counting lower bounds at a fixed $R_0$.
\end{proof}

\subsection{Fixed amplitudes and small Hecke coefficients}
\label{en:app:separation}

Corollary~\ref{en:cor:unified-applications} gives many integers of
each sign with a fixed bound on the number of prime factors. We now
also require $|x_n|$ to exceed a fixed positive constant. This will
give unweighted counts while retaining the marked-prime conditions
and the residue classes.

\begin{lemma}[Two signs of fixed amplitude]\label{en:app:amplitude}
Under the hypotheses of Proposition~\ref{en:thm:unified}, there are a
fixed $b>0$ and an integer $R_0$ such that each sign has
$\gg\Delta_C$ witnesses satisfying
\[
 C<n<2C,\quad \mu^2(n)=1,\quad \omega(n)\leq R_0,\quad
 H_C(n)>0,\quad n\bmod M\in\mathcal R_C,\quad \pm x_n\geq b.
\]
The dependencies of constants and the uniformity in $M,\mathcal R_C$
are as in the proposition.
\end{lemma}
\begin{proof}
We retain a fixed positive amount of the absolute mass of each sign,
then use a pointwise upper bound to count the remaining integers.
By Proposition~\ref{en:thm:unified}, choose fixed
$c_B,C_D,C_4>0$ such that, for sufficiently large $C$,
\[
 B_C\ge c_B\Delta_C,\qquad D_C\le C_D\Delta_C,
 \qquad D_{4,C}\le C_4\Delta_C.
\]
Since $T_C=o(\Delta_C)$, 
so that $|T_C|\le c_B\Delta_C/3$. The identities in the proof of
Corollary~\ref{en:cor:unified-applications} then give
\[
 B_C^\pm=\frac{B_C\pm T_C}{2}
       \ge\frac{c_B}{3}\Delta_C.
\]

Fix $0<b\le\min(1/2,c_B/(12C_D))$, and choose an integer $R_0$
large enough that $2^{-R_0}C_4\le c_B/12$. These choices are
independent of $C$, $M$, and $\mathcal R_C$. Using the nonnegative
weight $w_n$ from the proof of
Corollary~\ref{en:cor:unified-applications}, we have
\[
 \sum_{|x_n|<b}w_n|x_n|
       \le bD_C\le\frac{c_B}{12}\Delta_C,
 \qquad
 \sum_{\omega(n)>R_0}w_n|x_n|
       \le2^{-R_0}D_{4,C}\le\frac{c_B}{12}\Delta_C.
\]
The second inequality uses the same estimate as
\eqref{en:eq:unified-omega-tail}, together with $|x_n|\le1$.
After imposing both $|x_n|\ge b$ and
$\omega(n)\le R_0$, the absolute mass of each sign is therefore at
least
\[
 \left(\frac{c_B}{3}-\frac{c_B}{12}-\frac{c_B}{12}\right)\Delta_C
       =\frac{c_B}{6}\Delta_C.
\]

The pointwise bound in the proof of
Corollary~\ref{en:cor:unified-applications} shows that each retained
integer contributes at most $C_{\rm pt}:=\|g\|_\infty H_0\,4^{R_0}$.
If $N_\pm$ counts the retained integers of the given sign with
positive weight, then
\[
 C_{\rm pt}N_\pm\ge\frac{c_B}{6}\Delta_C.
\]
Thus $N_\pm\gg\Delta_C$.
\end{proof}

\begin{corollary}[Small eigenvalues in the two-marker example]\label{en:app:double}
Use the two-marker weight in \eqref{en:def:double-marker}. There are
fixed $b>0$ and $R_0$ such that, in every reduced class $a\pmod M$
in the range of Theorem~\ref{en:thm:main-hecke}, each sign has
\[
 \gg\frac{C}{\varphi(M)L^{3/2}}
\]
squarefree integers $n\in(C,2C)$ satisfying the product, ratio, and
two prime-eigenvalue conditions in \eqref{en:def:double-marker},
together with
\[
 \omega(n)\leq R_0,\qquad \pm x_n\geq b,\qquad
 \left|\frac{\tau(n)}{n^{11/2}}\right|
       \leq2^{R_0-2}L^{-1/2}.
\]
In particular, with $M_C=2^{\lfloor\log_2 L\rfloor}$ and $a=1$,
the number of integers of each sign is $\gg C/L^{5/2}$.
Thus the Hecke eigenvalue tends to zero, whereas
$\pm K(n)\geq b2^{\omega(n)}$.
\end{corollary}
\begin{proof}
Take $\mathcal R_C=\{a\}$. For the weight in
\eqref{en:def:double-marker}, the density calculation in
Section~\ref{en:app:main-proof} gives $d(C)\asymp L^{-1/2}$.
Hence $\Delta_C\asymp C/(\varphi(M)L^{3/2})$, and
Lemma~\ref{en:app:amplitude} gives the claimed count with a fixed
amplitude and a fixed bound on $\omega(n)$.

Here $k_1=2$ and $A_1=1/2$. For every retained integer,
\eqref{en:eq:normalized-hecke-small} and $\omega(n)\le R_0$ give
\[
 |\lambda_\Delta(n)|=|\tau(n)n^{-11/2}|
       \le2^{\omega(n)-2}L^{-1/2}
       \le2^{R_0-2}L^{-1/2}.
\]
Finally, $M_C\asymp L$ and $\varphi(M_C)\asymp L$, which give
the count $\gg C/L^{5/2}$. The relation between $K(n)$ and $x_n$
gives the claimed signed amplitude.
\end{proof}

Theorem~\ref{en:thm:main-hecke} already allows finitely many forms
$f_1,\ldots,f_s$. The next corollary combines that theorem with
Lemma~\ref{en:app:amplitude} to obtain unweighted counts with fixed
amplitude and a fixed bound on $\omega(n)$.

\begin{corollary}[Finitely many forms and arbitrary fixed logarithmic decay]\label{en:app:finite-forms}
Fix non-CM normalized holomorphic newforms $f_1,\ldots,f_s$ with
trivial nebentypus, and fix $A_1,\ldots,A_s>0$.
There exist fixed constants $\delta,b,c_1,\ldots,c_s>0$ and an
integer $R_0$ such that, for every sufficiently large $C$, every
$M\leq C^\delta$, and every $(a,M)=1$, each sign has
\begin{equation}\label{en:app:finite-count}
 \gg\frac{C}{\varphi(M)L^{1+A_1+\cdots+A_s}}
\end{equation}
integers satisfying
\begin{gather*}
 C<n<2C,\quad \mu^2(n)=1,\quad n\equiv a\pmod M,\quad
 \omega(n)\leq R_0,\quad \pm x_n\geq b,\\
 |\lambda_{f_i}(n)|\leq c_iL^{-A_i}\qquad(1\leq i\leq s).
\end{gather*}
Any fixed restrictions on the marked-prime positions given by an
admissible $F$ in \eqref{en:unif:marker} may also be retained.
The constants may depend on the forms, exponents, and fixed parameters
of Proposition~\ref{en:thm:unified}.
\end{corollary}
\begin{proof}
Choose fixed integers $k_i>2A_i$, and use the marked-prime conditions
and admissible regions described in
Section~\ref{en:unif:definition}. That construction gives suitable
fixed sieve parameters for the total number $k=\sum_i k_i$ of
marked primes. The verification in
Section~\ref{en:app:main-proof} applies to these choices, and
\eqref{en:eq:hecke-density} gives
\[
 d(C)\asymp L^{-A_1-\cdots-A_s}.
\]
For $\mathcal R_C=\{a\}$, Theorem~\ref{en:thm:main-hecke} and
Lemma~\ref{en:app:amplitude} give \eqref{en:app:finite-count}.
On every retained integer, $H_C(n)>0$, so the bound
\eqref{en:eq:normalized-hecke-small}, proved on this support in
Section~\ref{en:app:main-proof}, gives
\[
 |\lambda_{f_i}(n)|
 \le2^{\omega(n)-k_i}L^{-A_i}
 \le2^{R_0-k_i}L^{-A_i}.
\]
We may therefore take $c_i=2^{R_0-k_i}$. The marked-prime and residue
conditions are retained throughout the application of the lemma.
\end{proof}

\section{Explicit constants for two Hecke markers}
\label{en:num:main}
Fix $0\le g\le1$ with $\operatorname{supp}g\subset(1,2)$ and
$g=1$ on $[5/4,7/4]$. Define
\[
 H_C(n)=\sum_{\substack{p<q,\ pq\mid n\\
 C^{0.36}<pq<C^{0.38}\\C^{0.02}<q/p<C^{0.04}}}
 \mathbf1_{|\lambda_\Delta(p)|\le L^{-1/4}}
 \mathbf1_{|\lambda_\Delta(q)|\le L^{-1/4}},\qquad L=\log C,
\]
where $\lambda_\Delta(n)=\tau(n)/n^{11/2}$ is the normalized eigenvalue
of the Ramanujan form. Put
\[
 \mu_C=\mu_{\rm ST}\{|2\cos\theta|\le L^{-1/4}\},\qquad
 d_C=\mu_C^2\asymp L^{-1/2},\qquad
 \Delta_a=d_C\frac C{\varphi(M)L}.
\]
For a reduced residue class $a\pmod M$, let $D_a,T_a,B_a$ denote the
three masses in the main theorem. When the class is fixed, we also write
$D_C,T_C,B_C,\Delta_C$ for these quantities.
We use a seven-dimensional upper-bound sieve of total level
$y=C^\eta$, with $z=C^{\eta/66}$, preceded by an exact sieve at the
fixed primes below $e^{10}$. The coefficients are specified below.
As before, set
\[
 \theta=\xi*1\ge0,\qquad \Theta=\Lambda*\theta=\xi*\log,
 \qquad A_C=\Theta/\log,
 \qquad W_C=2^\omega A_CH_C\mathbf1_{n\equiv a\ (M)}.
\]
The coefficients satisfy $|\xi_d|\le1$ and vanish for $d>y$.
Consequently, the complete divisibility expansion has length at most
$C^{\eta+0.38}$.

\begin{theorem}[Explicit parameters using the zero-free half-plane]\label{en:num:theorem}
Using the zero-free half-plane for Dirichlet $L$-functions in the OpenAI
preprint \cite[Theorem 1.1]{en:oaizero26}, choose
\[
 \eta=\frac1{2000},\qquad \beta=\frac1{132000},\qquad
 R_0=83,\qquad b=10^{-21}.
\]
Uniformly over all integers $1\le M\le C^{11/200}$ and every reduced class
$a\pmod M$, for sufficiently large $C$ one has
\[
 D_a\asymp B_a\asymp\Delta_a,\qquad T_a=o(\Delta_a).
\]
For each sign, at least $3\cdot10^{-8}\Delta_a/(384\cdot4^{83})$
integers satisfy simultaneously
\[
 C<n<2C,\quad \mu^2(n)=1,\quad n\equiv a\pmod M,\quad
 H_C(n)>0,\quad\omega(n)\le83,\quad\pm x_n\ge10^{-21}.
\]
They also satisfy $|\lambda_\Delta(n)|\le2^{81}L^{-1/2}$.
The uniform starting point may depend on $g$, is independent of $M$,
and is not evaluated numerically.
\end{theorem}

\begin{remark}[Using only reference~\cite{en:oaisiegel26}]
If the zero-free half-plane is replaced by the Landau--Siegel zero gap in
\cite{en:oaisiegel26}, the range \(M\le C^{11/200}\) must be replaced by
\(M\le C^\delta\) for a sufficiently small, non-explicit \(\delta>0\).
The error \(O(C^{-1/100})\) in the three-prime estimate must also be
replaced by a fixed small error. With \(\delta\) chosen sufficiently
small, the constants \(B_*=3\cdot10^{-8}\), \(b=10^{-21}\), and
\(R_0=83\) remain unchanged.
\end{remark}

\subsection{The signed estimate and square divisors in a residue class}
\label{en:num:signed}
We prove the signed estimate $T_a=o(\Delta_a)$ in
Theorem~\ref{en:num:theorem}. Since $d_C\asymp L^{-1/2}$, for
$M\le C^d$ the mass to which we compare the sum satisfies
\[
 \Delta_a\asymp\frac{C}{\varphi(M)L^{3/2}}
           \ge c\frac{C^{1-d}}{L^{3/2}}
\]
for a fixed $c>0$. We will prove
$T_a\ll_\varepsilon C^{1-s+\varepsilon}$ with a fixed $s>d$,
and then choose $\varepsilon<s-d$. We obtain this estimate by
expanding the squarefree condition, counting the large square
divisors in the residue class, and applying a spectral estimate
to the remaining terms.

\paragraph{Reduction to short divisor sums.}
Put
\[
 u=\eta+0.38,\qquad a_0=1/2-u=0.12-\eta,\qquad Q=C^u.
\]
For coefficients $\boldsymbol b=(b_q)_{q\le Q}$ with $|b_q|\le1$
and a smooth test function $G$ supported in a fixed compact
subinterval of $(1,2)$, define
\[
 \mathcal S_a(G,\boldsymbol b)
 =\sum_{n\equiv a\ ({\rm mod}\ M)}
       \mu^2(n)G(n/C)K(n)\sum_{q\mid n,\ q\le Q}b_q.
\]
We need this estimate for exactly the two sums arising from
\eqref{en:eq:signed-exact-expansion}. In the notation of that
expansion, $H_C(n)=\sum_{P\mid n}w_C(P)$, now with two marked
primes and $P<C^{0.38}$. The coefficients constructed in
Section~\ref{en:r3:presieve} have $e\le C^\eta$ and all prime
factors at most $z$. Both marked primes exceed $C^{0.16}>z$.
Thus $(e,P)=1$, and $q=eP$ has at most one such representation:
its prime-power factors at primes at most $z$ recover $e$, and
the remaining prime factors recover $P$.

Define
\[
 b_q^{(1)}=\xi_e w_C(P),\qquad
 b_q^{(2)}=\xi_e\frac{\log e}{L}w_C(P)
 \quad(q=eP),
 \qquad G_C(t)=\frac{Lg(t)}{L+\log t},
\]
with the coefficients equal to zero when there is no representation.
Here $G_C$ is extended by zero outside $(1,2)$. We have $q\le Q$ and
$|b_q^{(1)}|,|b_q^{(2)}|\le1$, because
$0\le w_C(P)\le1$ and $0\le\log e/L\le\eta<1$.
The same expansion, with the residue condition retained, gives
the exact identity
\[
 \begin{split}
 T_a
 &=\sum_{n\equiv a\ ({\rm mod}\ M)}
       \mu^2(n)g(n/C)K(n)A_C(n)H_C(n)\\
 &=\mathcal S_a(g,\boldsymbol b^{(1)})
       -\mathcal S_a(G_C,\boldsymbol b^{(2)}).
 \end{split}
\]
Both test functions have the same fixed compact support, and every
derivative is bounded uniformly in $C$. We therefore estimate
$\mathcal S_a(G,\boldsymbol b)$ uniformly for this family of tests
and all the stated coefficients.

\paragraph{The spectral estimate and the required divisor sums.}
The pointwise estimates in the proof of Drappeau--Maynard
\cite[\S3.1, (3.2)--(3.3)]{en:dm} give, for every character
$\chi$ modulo $M$ and every $\ell$ with
$\ell\sqrt M\le C^{1/2-\nu}$, where $\nu>0$ is fixed,
\begin{equation}\label{en:num:dm-pointwise}
 \left|\sum_{\ell\mid n}G(n/C)\chi(n)K(n)\right|
 \ll_{\varepsilon,\nu} C^{1/2+\varepsilon}
 \left\{M^{1/2}+
 \left(\frac{\sqrt{CM}}\ell\right)^{7/16}
 \left(1+\frac{\sqrt M}\ell\right)\right\}.
\end{equation}
This uses the admissible spectral exponent $\vartheta=7/64$ and
retains the displayed dependence on $M$. We use their combined
estimate directly; the underlying spectral bounds are
\cite[Lemmas 2.7--2.9]{en:topacogullari}. If $(\ell,M)>1$, every
character value in the sum is zero.
For fixed $\varepsilon,\nu$, divide these linear forms in $G$
by the positive expression on the right of
\eqref{en:num:dm-pointwise}. Applying the finite-seminorm argument
used for \eqref{en:eq:uniform-test} to this family of normalized
forms makes the estimate uniform over $g$ and $G_C$. Their fixed
support and bounded derivatives give the required uniform seminorm
bound.

The identity $\mu^2(n)=\sum_{r^2\mid n}\mu(r)$ will replace the
conditions $q\mid n$ and $r^2\mid n$ by $[q,r^2]\mid n$.
To sum \eqref{en:num:dm-pointwise} over these pairs, we use
\begin{equation}\label{en:num:lcm}
 \sum_{r\le R,q\le Q}[q,r^2]^{-v}
 \ll_{v,\varepsilon}Q^{1-v}R^{1-2v+\varepsilon}
 \quad(0<v<1/2),\qquad
 \sum_{r\le R,q\le Q}[q,r^2]^{-w}\ll_w1
 \quad(w>1).
\end{equation}
Here $[q,r^2]$ denotes the least common multiple. To prove the
first bound, use $[q,r^2]=qr^2/(q,r^2)$ and
$(q,r^2)^v\le\sum_{h\mid(q,r^2)}h^v$. For fixed $r$ this gives
\begin{align*}
 \sum_{q\le Q}[q,r^2]^{-v}
 &\le r^{-2v}\sum_{h\mid r^2}h^v
                  \sum_{\substack{q\le Q\\h\mid q}}q^{-v}\\
 &\ll_v r^{-2v}Q^{1-v}\sum_{h\mid r^2}h^{v-1}
 \ll_{v,\varepsilon}r^{-2v+\varepsilon}Q^{1-v}.
\end{align*}
Summing over $r\le R$ proves the claim. For the second bound,
extend both sums to all positive integers:
\[
 \begin{split}
 \sum_{r,q\ge1}[q,r^2]^{-w}
 &\le\sum_{r\ge1}r^{-2w}\sum_{h\mid r^2}h^w
                     \sum_{\substack{q\ge1\\h\mid q}}q^{-w}\\
 &=\zeta(w)\sum_{r\ge1}\frac{\tau_2(r^2)}{r^{2w}}<\infty.
 \end{split}
\]
The last series converges since $w>1$ and the usual divisor
function satisfies $\tau_2(r^2)\ll_\varepsilon r^\varepsilon$.

\paragraph{Large square divisors.}
Suppose $M\le C^d$, choose $s>d$, and put $R=C^s/M$.
In particular, $R\ge C^{s-d}\to\infty$ uniformly in $M$.
Expanding the squarefree condition gives
\[
 \mathcal S_a(G,\boldsymbol b)
 =\sum_{r\le\sqrt{2C}}\mu(r)
       \sum_{\substack{n\equiv a\ ({\rm mod}\ M)\\r^2\mid n}}
            G(n/C)K(n)\sum_{q\mid n,\ q\le Q}b_q.
\]
Write $\mathcal S_{a,>R}$ and $\mathcal S_{a,\le R}$ for the two
parts of this sum. The support of $G$ restricts $n$ to $(C,2C)$.
For the part with $r>R$, we claim
\begin{equation}\label{en:num:square-tail}
 |\mathcal S_{a,>R}|
 \ll_\varepsilon C^\varepsilon
       \left(\frac{C}{MR}+\sqrt C\right)
 =O_\varepsilon(C^{1-s+\varepsilon}+C^{1/2+\varepsilon}).
\end{equation}
Indeed, the divisor bound and the Weil--Estermann bound for all
moduli \cite[Corollary 11.12]{en:ik} give
\[
 \left|\sum_{q\mid n,\ q\le Q}b_q\right|\le\tau_2(n),\qquad
 |K(n)|\le\tau_2(n),\qquad
 \tau_2(n)^2\ll_\varepsilon C^\varepsilon.
\]
If $r^2\mid n$ and $n\equiv a\pmod M$, then $(a,M)=1$ forces
$(r,M)=1$. The Chinese remainder theorem specifies one residue
class modulo $Mr^2$, containing $O(C/(Mr^2)+1)$ integers in
$(C,2C)$. Therefore
\[
 |\mathcal S_{a,>R}|
 \ll_\varepsilon C^\varepsilon
       \sum_{R<r\le\sqrt{2C}}
                   \left(\frac C{Mr^2}+1\right)
 \ll_\varepsilon C^\varepsilon
                   \left(\frac C{MR}+\sqrt C\right).
\]
This proves \eqref{en:num:square-tail}.

\paragraph{Small square divisors.}
For $r\le R$, character orthogonality gives the exact expansion
\[
 \mathcal S_{a,\le R}
 =\frac1{\varphi(M)}\sum_{\chi\ ({\rm mod}\ M)}
       \overline{\chi(a)}
       \sum_{r\le R}\mu(r)\sum_{q\le Q}b_q
       \sum_{[q,r^2]\mid n}G(n/C)\chi(n)K(n).
\]
If $2s<a_0$, then for every pair $q,r$ in this expression,
\[
 [q,r^2]\sqrt M\le QR^2\sqrt M
       =C^{u+2s}M^{-3/2}
       \le C^{1/2-(a_0-2s)}.
\]
Thus \eqref{en:num:dm-pointwise} applies with the fixed choice
$\nu=a_0-2s>0$. The character coefficients have total absolute
value $\varphi(M)^{-1}\sum_\chi1=1$.

Expand the right side of \eqref{en:num:dm-pointwise} into three
terms. Apart from an arbitrarily small power of $C$, the first
is $C^{1/2}M^{1/2}$ for each of at most $QR$ pairs, giving
\[
 C^{1/2}M^{1/2}QR
       =C^{1/2+u+s}M^{-1/2}
       =C^{1-a_0+s}M^{-1/2}.
\]
For the second term, use \eqref{en:num:lcm} with $v=7/16$:
\[
 \begin{split}
 C^{23/32}M^{7/32}Q^{9/16}R^{1/8}
 &=C^{23/32+9u/16+s/8}M^{7/32-1/8}\\
 &=C^{1-9a_0/16+s/8}M^{3/32}.
 \end{split}
\]
The third term is $C^{23/32}M^{23/32}[q,r^2]^{-23/16}$.
Its double sum is bounded by the second estimate in
\eqref{en:num:lcm}, with $w=23/16>1$.
Apply the preliminary estimates with smaller positive error
exponents, so that all such factors combine into one $C^\varepsilon$.
We obtain
\begin{equation}\label{en:num:spectral-budget}
 |\mathcal S_{a,\le R}|\ll_\varepsilon C^\varepsilon
 \left(C^{1-a_0+s}M^{-1/2}
       +C^{1-9a_0/16+s/8}M^{3/32}
       +C^{23/32}M^{23/32}\right).
\end{equation}

\paragraph{Choice of parameters.}
The conditions
\begin{equation}\label{en:num:feasible}
 d<s,\qquad 2s<a_0,\qquad 12s+d<6a_0,\qquad
 s+\frac{23}{32}d<\frac9{32}
\end{equation}
make every term above $O_\varepsilon(C^{1-s+\varepsilon})$.
For the first term, use $M^{-1/2}\le1$; the difference of its
exponent from $1-s$ is $-a_0+2s<0$.
For the second, $M\le C^d$ gives the difference
\[
 -\frac{9a_0}{16}+\frac{9s}{8}+\frac{3d}{32}<0,
\]
which is equivalent to $12s+d<6a_0$. The difference for the
third term is $-9/32+s+23d/32<0$. Also $2s<a_0<1/2$
implies $s<1/2$, so the term $C^{1/2+\varepsilon}$ in
\eqref{en:num:square-tail} is smaller than
$C^{1-s+\varepsilon}$. Combining the two parts proves
\[
 \mathcal S_a(G,\boldsymbol b)\ll_\varepsilon C^{1-s+\varepsilon}.
\]
The exact two-sum identity at the start of this subsection then
gives the same bound for $T_a$.

For Theorem~\ref{en:num:theorem}, take
\[
 \eta=\frac1{2000},\qquad d=\frac{11}{200},\qquad
 s=\frac{551}{10000},\qquad
 u=\frac{761}{2000},\qquad a_0=\frac{239}{2000}.
\]
The four required inequalities have the following positive bounds:
\begin{equation}\label{en:num:single-margins}
 \begin{aligned}
 a_0-2s&=\frac{93}{10000},&
 \frac{9a_0}{16}-\frac{9s}{8}-\frac{3d}{32}&=\frac3{40000},\\
 \frac9{32}-s-\frac{23d}{32}&=\frac{29859}{160000},&
 s-d&=\frac1{10000}.
 \end{aligned}
\end{equation}
Choose a fixed $0<\varepsilon<3/80000$. This is smaller than
half the least of these bounds. Since
$\Delta_a\asymp C/(\varphi(M)L^{3/2})$ and
$\varphi(M)\le M\le C^d$, the resulting bound gives
\[
 \frac{|T_a|}{\Delta_a}
 \ll C^{d-s+\varepsilon}L^{3/2}
 =C^{-1/10000+\varepsilon}L^{3/2}\longrightarrow0.
\]
This proves $T_a=o(\Delta_a)$ uniformly in the modulus and the
reduced class in Theorem~\ref{en:num:theorem}.

\subsection{Prime products in each reduced residue class}
\label{en:num:prime-products}

We now prove the congruence estimate needed for the absolute lower bound:
each reduced class contains a proportion
\[
 \frac{1+O(C^{-1/100})}{\varphi(M)}
\]
of every three-prime box in Section~\ref{en:num:free-boxes}, uniformly
for $M\le C^{11/200}$.
This will allow us to use three free primes for every permitted modulus.
Fix a marked pair $p,q$, put $P=pq$, and suppose $(P,M)=1$.
Then $Pq_1q_2q_3\equiv a\pmod M$ is equivalent to
$q_1q_2q_3\equiv b\pmod M$, where $b=aP^{-1}\pmod M$.
The construction in Section~\ref{en:num:free-boxes} has
$p,q>C^{0.165}$ and free-prime intervals
$J_i=[Q_i,\lambda Q_i)$ with $Q_i\ge C^{0.171+o(1)}$.
Thus all these primes exceed $M$. We prove the estimate for any such
intervals in a fixed power range, with fixed $\lambda>1$.

\paragraph{Character expansion of the required proportion.}
Use the notation $N_i=N_{J_i}$ and $a_i(\chi)=a_{J_i}(\chi)$ from
Section~4.1. Thus $N_i$ counts the primes in $J_i$, and $a_i(\chi)$ is
their average character value. Put
\[
 N_b=\#\{(q_1,q_2,q_3)\in\textstyle\prod_{i=1}^3(J_i\cap\mathbb P):
                    q_1q_2q_3\equiv b\pmod M\}.
\]
Character orthogonality gives exactly
\[
 \frac{N_b}{N_1N_2N_3}
 =\frac1{\varphi(M)}
   \left(1+\sum_{\substack{\chi\ ({\rm mod}\ M)\\\chi\ne\chi_0}}
          \overline{\chi(b)}a_1(\chi)a_2(\chi)a_3(\chi)\right).
\]
The Brun--Titchmarsh and orthogonality argument giving
\eqref{en:unif-abs:l2} also applies here: the required comparison is
$\log(Q_i/M)\asymp\log Q_i$, which holds in the given ranges. It gives
\[
 \sum_{\chi\ ({\rm mod}\ M)}|a_i(\chi)|^2\ll_\lambda1.
\]
Consequently, Cauchy--Schwarz bounds the nonprincipal part by
\[
 \max_{\chi\ne\chi_0}|a_3(\chi)|
 \left(\sum_\chi|a_1(\chi)|^2\right)^{1/2}
 \left(\sum_\chi|a_2(\chi)|^2\right)^{1/2}
 \ll_\lambda\max_{\chi\ne\chi_0}|a_3(\chi)|.
\]
It remains to bound each nonprincipal character average. For $M=1$,
there are no nonprincipal characters.

\paragraph{The required character average.}
Fix $0<d<a_*<a^*$. We will prove
\begin{equation}\label{en:num:prime-bias}
 \max_{\substack{\chi\ ({\rm mod}\ M)\\\chi\ne\chi_0}}
 \frac{\left|\sum_{Y\le r<\lambda Y\,,\ r\in\mathbb P}\chi(r)\right|}
      {\#([Y,\lambda Y)\cap\mathbb P)}
 \ll Y^{-1/8}L^2\ll C^{-a_*/16},
 \qquad M\le C^d,\quad C^{a_*}\le Y\le C^{a^*}.
\end{equation}
Here we use the zero-free statement asserted in
OpenAI \cite[Theorem 1.1]{en:oaizero26}: Dirichlet $L$-functions have no zeros in
$\Re s>7/8$. This statement covers nonreal characters and all heights.
For a primitive nonprincipal character $\chi'$, the functional equation
reflects a nontrivial zero $\rho$ to $1-\rho$ for the conjugate character.
Applying the same input to that character gives
\[
 \frac18\le\Re\rho\le\frac78.
\]

Let $\chi'$ induce $\chi$, and write $q_*\mid M$ for its conductor.
Since $\chi$ is nonprincipal, $q_*>1$.
Let $N_{\chi'}(T)$ count the nontrivial zeros with
$|\Im\rho|\le T$, including multiplicities. The standard zero-counting
estimate \cite{en:ik} is
\[
 N_{\chi'}(T)\ll(T+1)\log(q_*(T+2)).
\]
For the zeros with $|\Im\rho|\le1$, the inequality $|\rho|\ge1/8$
therefore gives $O(\log(3q_*))$. For the other zeros,
partial summation gives
\[
 \sum_{1<|\Im\rho|\le2Y}\frac1{|\rho|}
 \le\frac{N_{\chi'}(2Y)}{2Y}
       +\int_1^{2Y}\frac{N_{\chi'}(t)}{t^2}\,dt
 \ll\log^2(q_*Y).
\]
Combining these two ranges, we obtain
$\sum_{|\Im\rho|\le2Y}|\rho|^{-1}\ll\log^2(q_*Y)$,
with a constant independent of $q_*$.

Define $\psi(t,\chi')=\sum_{n\le t}\chi'(n)\Lambda(n)$.
The truncated explicit formula \cite{en:ik}, at a height
$T_1\in[Y,2Y]$ chosen away from zero ordinates, gives, for half-integers
$t\asymp_\lambda Y$,
\[
 \psi(t,\chi')=-\sum_{|\Im\rho|\le T_1}\frac{t^\rho}{\rho}
 +O_\lambda\!\left(\frac{t\log^2(q_*tY)}Y+\log^2(q_*t)\right).
\]
The zero sum is at most
$t^{7/8}\sum_{|\Im\rho|\le2Y}|\rho|^{-1}$ in absolute value.
Since $\log(q_*Y)\ll L$, this proves
\[
 |\psi(t,\chi')|\ll_\lambda Y^{7/8}L^2.
\]
This bound also holds for every real $Y\le t\le\lambda Y$.
Indeed, $\psi(t,\chi')$ is constant between consecutive integers;
when $t$ passes an integer $n$, its change is $\chi'(n)\Lambda(n)$,
of absolute value at most $\log n$. Replacing $t$ by a neighboring
half-integer therefore changes the sum by only $O(\log(2t))$.

Removing the terms $n=r^j$ with $j\ge2$ costs
$O(\sqrt Y\log^2Y)$, which is smaller than the bound above.
Hence, with
$\vartheta_{\chi'}(t)=\sum_{r\le t\,,\ r\in\mathbb P}\chi'(r)\log r$,
we have $\vartheta_{\chi'}(t)\ll_\lambda Y^{7/8}L^2$ on the same interval.
Since $Y>M$, the characters $\chi$ and $\chi'$ agree at every prime
in $[Y,\lambda Y)$. Partial summation now gives
\[
 \begin{aligned}
 \sum_{Y\le r<\lambda Y\,,\ r\in\mathbb P}\chi(r)
 &=\frac{\vartheta_{\chi'}(\lambda Y)}{\log(\lambda Y)}
  -\frac{\vartheta_{\chi'}(Y)}{\log Y}
  +\int_Y^{\lambda Y}\frac{\vartheta_{\chi'}(t)}{t\log^2t}\,dt
  +O(1)\\
 &\ll_\lambda\frac{Y^{7/8}L^2}{\log Y}.
 \end{aligned}
\]
The $O(1)$ term accounts for the two endpoints of the half-open interval.
The prime number theorem gives
$\#([Y,\lambda Y)\cap\mathbb P)\asymp_\lambda Y/\log Y$.
Dividing by this value proves the first bound in
\eqref{en:num:prime-bias}. Finally,
$Y^{-1/8}\le C^{-a_*/8}$ and $L^2\le C^{a_*/16}$ for sufficiently
large $C$, which proves the second bound. The constants depend only on
the fixed power ranges and $\lambda$, not on the prime factors of $M$.

\paragraph{The resulting three-prime count.}
Return to the intervals $J_i$ above. Take $d=11/200$ and $a_*=0.17$
in \eqref{en:num:prime-bias}. Since $Q_i\ge C^{0.171+o(1)}$ and
$0.17/16>1/100$, we obtain
$\max_{\chi\ne\chi_0}|a_i(\chi)|\ll C^{-1/100}$.
Substituting this into the character expansion proves
\begin{equation}\label{en:num:three-prime-class}
 \frac{N_b}{N_1N_2N_3}
 =\frac1{\varphi(M)}\left(1+O(C^{-1/100})\right).
\end{equation}
In particular, the proportion is at least
$9/(10\varphi(M))$ for sufficiently large $C$, uniformly in the boxes
and the reduced class $b$. Three free primes therefore suffice, without
the parity choice used in Section~4.1 for an exceptional real character.

For $M=p_0^m$ with $p_0$ fixed, Banks--Shparlinski
\cite[Theorem 3.2]{en:banksShparlinskiSmooth}, combined with the
classical zero-free region \cite[Theorem 12.3]{en:koukoulopoulos}
and separate treatment of the finitely many primitive real characters,
gives \eqref{en:num:three-prime-class} with relative error
$O(\exp(-c(L/ \log L)^{1/3}))$.
Thus the numerical lower bounds below also hold for this family without
the half-plane input; the extension to arbitrary moduli uses
\cite[Theorem 1.1]{en:oaizero26}.

\subsection{A fixed pre-sieve in dimension seven}
\label{en:r3:presieve}
We now construct the sieve coefficients announced at the start of
this section. The next subsection needs one choice of coefficients
that works for both $j=2$ and $j=7$. Our goal is
\begin{equation}\label{en:r3:fundamental-sharp}
 \sum_d\xi_d\prod_{p\mid d}\frac{j}{p+j}
 \le F_*\prod_{p<z}(1+j/p)^{-1},\qquad j=2,7,
 \qquad F_*=1+e^{-2.29}.
\end{equation}
The same coefficients will also satisfy the support and positivity
conditions needed for the signed estimate and the absolute lower bound.

Let $w_0=e^{10}$, $Q_0=\prod_{p<w_0}p$, and write
$P(z)=\prod_{p<z}p$. With $y=C^\eta$ and $z=C^\beta$ as at the
start of the section, we have $\beta=\eta/66$. Put
\[
 D=\frac y{Q_0},\qquad
 s_{\mathrm{hi}}=\frac{\log D}{\log z}
     =66-\frac{\log Q_0}{\beta\log C}.
\]
Since $Q_0$ is fixed, we may take $C$ sufficiently large that
$z>w_0$ and $s_{\mathrm{hi}}\ge65.99>64=9\cdot7+1$.

Apply the fundamental lemma
\cite[Lemma 9.1]{en:matomakiTeravainen} in dimension $7$, with
parameters $(D,z)$. It gives upper-bound coefficients
$\lambda_d^+$, chosen independently of the density function, such that
\[
 |\lambda_d^+|\le1,\qquad
 \lambda_d^+=0\quad\hbox{unless }d\le D\hbox{ and }d\mid P(z),
 \qquad
 \sum_{d\mid n}\lambda_d^+\ge\mathbf1_{(n,P(z))=1}.
\]
Taking $n=1$ shows that $\lambda_1^+=1$.
Define $\lambda_b^{\mathrm{hi}}=\lambda_b^+$ when $(b,Q_0)=1$,
and $\lambda_b^{\mathrm{hi}}=0$ otherwise. Now set
\[
 \xi_{ab}=\mu(a)\lambda_b^{\mathrm{hi}}
 \quad\hbox{for }a\mid Q_0,\ (b,Q_0)=1,
\]
and set all other coefficients to zero. The decomposition into the
part supported on primes below $w_0$ and the remaining part is unique.
Thus $\xi_1=1$, $|\xi_d|\le1$, and a nonzero coefficient has
$d\mid P(z)$ and $d\le Q_0D=y$. In particular, the coefficients
are supported on squarefree integers and have total level at most $y$.

For the function $\theta$ already defined at the start of the section,
this construction gives
\[
 \theta(n)
 =\left(\sum_{a\mid(n,Q_0)}\mu(a)\right)
       \sum_{b\mid n}\lambda_b^{\mathrm{hi}}
 =\mathbf1_{(n,Q_0)=1}\sum_{b\mid n}\lambda_b^{\mathrm{hi}}
 \ge\mathbf1_{(n,P(z))=1}\ge0.
\]
The convolution identity
\eqref{en:eq:convolution-identity} now gives $A_C(n)\ge0$.
If $P^-(n)>z$, only $d=1$ contributes to that identity, so
$A_C(n)=1$. 

It remains to prove \eqref{en:r3:fundamental-sharp} with its
constant. For $j=2,7$, let $h_j$ and $h_j^{\mathrm{hi}}$ be
multiplicative functions with prime values
\[
 h_j(p)=\frac{j}{p+j},\qquad
 h_j^{\mathrm{hi}}(p)=
 \begin{cases}
 0,&p<w_0,\\
 j/(p+j),&p\ge w_0.
 \end{cases}
\]
Only their values on squarefree integers will be used. We check the
product hypothesis of the fundamental lemma for both functions
$h_j^{\mathrm{hi}}$ with the same dimension and constant.

The two bounds in Rosser--Schoenfeld
\cite[Theorem 5, (3.17)--(3.18)]{en:rosserSchoenfeld}, applied at
the endpoints, give, for $v\ge u\ge w_0$,
\[
 \sum_{u\le p<v}\frac1p
 \le\log\frac{\log v}{\log u}
       +\frac1{2\log^2v}+\frac1{2\log^2u}
 \le\log\frac{\log v}{\log u}+\frac1{\log^2u}.
\]
Since $(1-h_j^{\mathrm{hi}}(p))^{-1}=1+j/p$ for $p\ge w_0$,
we obtain
\[
 \begin{aligned}
 \log\prod_{u\le p<v}(1-h_j^{\mathrm{hi}}(p))^{-1}
 &\le j\sum_{u\le p<v}\frac1p\\
 &\le j\log\frac{\log v}{\log u}+\frac j{\log^2u}
 \le7\log\frac{\log v}{\log u}+\frac7{100}.
 \end{aligned}
\]
Here we used $j\le7$, $\log u\ge10$, and $\log(1+t)\le t$.
If $u<w_0\le v$, apply this bound from $w_0$ to $v$; if
$v\le w_0$, the product is $1$. Hence, for every $v\ge u\ge2$,
\[
 \prod_{u\le p<v}(1-h_j^{\mathrm{hi}}(p))^{-1}
 \le K\left(\frac{\log v}{\log u}\right)^7,
 \qquad K=e^{7/100}.
\]

Part (iii) of the same fundamental lemma therefore applies to the
already chosen coefficients $\lambda_d^+$, for each $j=2,7$.
It gives
\[
 \begin{aligned}
 \sum_b\lambda_b^{\mathrm{hi}}h_j(b)
 &=\sum_d\lambda_d^+h_j^{\mathrm{hi}}(d)\\
 &\le\bigl(1+e^{63-s_{\mathrm{hi}}}K^{10}\bigr)
                \prod_{w_0\le p<z}(1-h_j(p)).
 \end{aligned}
\]
The upper factor satisfies
\[
 1+e^{63-s_{\mathrm{hi}}}K^{10}
 \le1+\exp(63-65.99+7/10)
 =1+e^{-2.29}=F_*.
\]

Finally, multiplicativity and the coprimality of $a$ and $b$ give
the exact factorization
\[
 \sum_d\xi_dh_j(d)
 =\left(\sum_{a\mid Q_0}\mu(a)h_j(a)\right)
       \left(\sum_b\lambda_b^{\mathrm{hi}}h_j(b)\right).
\]
The first factor is the positive product
$\prod_{p<w_0}(1-h_j(p))$. Multiplying the upper bound
by this product yields
\[
 \sum_d\xi_dh_j(d)
 \le F_*\prod_{p<z}(1-h_j(p))
 =F_*\prod_{p<z}(1+j/p)^{-1}.
\]
This proves \eqref{en:r3:fundamental-sharp} for $j=2$ and $j=7$
with one set of coefficients.

\subsection{Explicit upper bounds for the positive sums}
\label{en:r3:positive}
We prove \eqref{en:r3:positive-final}, namely
\[
 D_{j,a}\le14j\beta^{-j}\Delta_a,\qquad j=2,7,
 \qquad M\le C^{11/200},\quad(a,M)=1,
\]
where
\[
 f_j(n)=\mu^2(n)j^{\omega(n)},\qquad
 D_{j,a}=\sum_{n\equiv a\ ({\rm mod}\ M)}
              f_j(n)g(n/C)A_C(n)H_C(n).
\]
The case $j=2$ bounds $D_a$. The case $j=7$ will bound the contribution
from integers with many prime factors in Section~\ref{en:r3:tail}.
We first use \eqref{en:r3:fundamental-sharp} without the marked-prime or
congruence conditions. We then restore the prime restrictions, count the
possible marked pairs for each integer, and impose the residue class.

\paragraph{The bound before imposing prime restrictions.}
For $X>y$, the squarefree convolution expansion from
Section~\ref{en:unif:positive-proof} gives
\[
 \sum_{n\le X}f_j(n)\Theta(n)
 =j\sum_{m\le X}f_j(m)\theta(m)
       \sum_{\substack{r\le X/m\\r\text{ prime},\ r\nmid m}}\log r.
\]
All outer weights are nonnegative. Delete $r\nmid m$ and use the bound
$\sum_{r\le t,\ r\text{ prime}}\log r<1.01624t$ from
Rosser--Schoenfeld \cite[Theorem 9]{en:rosserSchoenfeld}. Enlarging the
$m$-sum by nonnegativity gives
\[
 \sum_{n\le X}f_j(n)\Theta(n)
 \le1.01624jX\sum_{P^+(m)\le X}\frac{f_j(m)\theta(m)}m.
\]
The finite Euler-product identity proved in
Section~\ref{en:unif:positive-proof}, followed by
\eqref{en:r3:fundamental-sharp}, gives
\[
 \begin{aligned}
 \sum_{P^+(m)\le X}\frac{f_j(m)\theta(m)}m
 &=\prod_{p\le X}(1+j/p)\sum_d\xi_dh_j(d)\\
 &\le F_*\prod_{z\le p\le X}(1+j/p).
 \end{aligned}
\]
Here $h_j$ is the multiplicative function from
Section~\ref{en:r3:presieve}, and $d\le y<X$.

Fix $\lambda=1.01$ and put $X=2\lambda^2C$. Since $z=C^\beta$,
Mertens' formula yields
\[
 \begin{aligned}
 \log\prod_{z\le p\le X}(1+j/p)
 &=j\sum_{z\le p\le X}\frac1p
       +O_j\!\left(\sum_{p\ge z}p^{-2}\right)\\
 &=j\log\frac{\log X}{\log z}+o(1).
 \end{aligned}
\]
Thus the product is $(1+o(1))\beta^{-j}$. On
$C/\lambda^2<n<2\lambda^2C$, we have
$\log n\ge L-2\log\lambda=(1+o(1))L$. Dividing the cumulative bound
by this lower bound for $\log n$ proves
\begin{equation}\label{en:r3:base-positive}
 \sum_{C/\lambda^2<n<2\lambda^2C}f_j(n)A_C(n)
 \le\bigl(2\lambda^2\cdot1.01624F_*+o(1)\bigr)
       j\beta^{-j}\frac CL.
\end{equation}

\paragraph{Restoring the two prime restrictions.}
Let $D_j^{\rm all}$ denote $D_{j,a}$ with the congruence condition omitted.
For this comparison, partition the marked-prime ranges into half-open
intervals $[\lambda^u,\lambda^{u+1})$ and retain the boxes
$\mathcal J=I_1\times I_2$ meeting the product and ratio conditions
in $H_C$. Write $I_i=[Y_i,\lambda Y_i)$ and
\[
 \mathcal M_{\mathcal J}
 =\left\{m\ge1:\frac{C}{\lambda^2Y_1Y_2}<m<\frac{2C}{Y_1Y_2}\right\}.
\]
If $p\in I_1$, $q\in I_2$, and $g(pqm/C)\ne0$, then
$m\in\mathcal M_{\mathcal J}$. For every such $m$, replacing $p,q$
by any $p'\in I_1$, $q'\in I_2$ gives
\[
 \frac C{\lambda^2}<p'q'm<2\lambda^2C.
\]
This is the enlarged interval used in \eqref{en:r3:base-positive}.

The verification in Section~\ref{en:app:main-proof}, with
$L^{-1/4}$, gives uniformly in these intervals
\[
 \#\{r\in I_i:r\text{ prime},\ |\lambda_\Delta(r)|\le L^{-1/4}\}
 =(1+o(1))\mu_C\,\#(I_i\cap\mathbb P),\qquad i=1,2.
\]
Indeed, the relative error there is
$O((\log L)L^{-1/2}/\mu_C)=O((\log L)L^{-1/4})=o(1)$.

All marked primes exceed $z$. The identity
\eqref{en:unif:positive-cofactor} therefore gives
\[
 A_C(pqm)=\frac{\Theta(m)+\log(pq)\theta(m)}{\log(pqm)},
\]
even if $p$ or $q$ divides $m$. With the box and $m$ fixed, replacing
$p,q$ by $p',q'$ changes each logarithm by $O(1)$. The two numerator
terms are nonnegative, $\log(pq)\asymp L$, and
$\log(pqm)\asymp L$. Hence
\[
 A_C(p'q'm)=(1+O(1/L))A_C(pqm)
\]
uniformly over all these choices. If the numerator is zero, both sides
are zero.

Apply the box comparison from Section~\ref{en:unif:positive-proof}
using $g\le1$ and the indicator defining the marked-prime positions.
The two prime counts above give the factor
$(1+o(1))\mu_C^2=(1+o(1))d_C$. Removing the coprimality condition
with $m$ enlarges the nonnegative sum. Thus
\[
 D_j^{\rm all}\le(1+o(1))d_C
 \sum_{\mathcal J}\sum_{m\in\mathcal M_{\mathcal J}}j^2f_j(m)
       \sum_{\substack{p'\in I_1\cap\mathbb P\\q'\in I_2\cap\mathbb P}}
                      A_C(p'q'm).
\]

\paragraph{Counting representations and removing square divisors.}
The original product and ratio bounds imply
\[
 0.16<\frac{\log p}{L}<0.18,\qquad
 0.19<\frac{\log q}{L}<0.21.
\]
Thus every retained box lies in
\[
 C^{0.16}/\lambda\le p'\le\lambda C^{0.18},\qquad
 C^{0.19}/\lambda\le q'\le\lambda C^{0.21}.
\]
These ranges are disjoint for large $C$. For a fixed integer
$n=p'q'm<2\lambda^2C$, let $a$ and $b$ count its distinct prime factors
in the first and second range. Both counts are bounded, since all these
primes are at least $C^{0.16}/\lambda$. Their product divides $n$, so
\[
 (C^{0.16}/\lambda)^a(C^{0.19}/\lambda)^b\le n<2\lambda^2C,
 \qquad 0.16a+0.19b\le1+O(1/L).
\]
If $ab\ge7$, one of the pairs
$(1,7),(2,4),(3,3),(4,2),(7,1)$ is coordinatewise no larger than
$(a,b)$. On these five pairs, the minimum of $0.16a+0.19b$ is $1.02$.
For large $C$ this contradicts the upper bound. Hence $ab\le6$.

There are at most $ab$ choices for $(p',q')$; each fixes $m=n/(p'q')$
and its unique half-open grid box. Thus every $n$ has at most six
representations in the complete sum. This argument also applies when
$n$ is not squarefree. In particular, it proves $H_C(n)\le6$ on the
original interval $C<n<2C$.

If $(p'q',m)=1$, then
$j^2f_j(m)=f_j(n)$. If $(p'q',m)>1$ and $f_j(m)\ne0$, some prime
$r\ge C^{0.16}/\lambda$ divides both $p'q'$ and $m$, so $r^2\mid n$.
Let $E_j$ be the sum of the formal weights $j^2f_j(m)A_C(n)$ over
all such additional terms. The divisor bound used in
\eqref{en:eq:marker-collision-bound} gives each formal weight
$O_\varepsilon(C^\varepsilon)$, even for nonsquarefree $n$.
Enlarging the sum over $r$ to all integers in the stated range,
the representation count and the number of multiples of $r^2$ give
\[
 \begin{aligned}
 E_j
 &\ll_\varepsilon C^\varepsilon
   \sum_{r\ge C^{0.16}/\lambda}
        \#\{n\le2\lambda^2C:r^2\mid n\}\\
 &\ll_\varepsilon C^{1+\varepsilon}
               \sum_{r\ge C^{0.16}/\lambda}r^{-2}
 \ll_\varepsilon C^{0.84+\varepsilon}=o(C/L),
 \end{aligned}
\]
where the last "=" holds for a fixed $\varepsilon<0.16$.
We have bounded all the additional terms together, over all boxes.
Consequently,
\[
 \begin{aligned}
 D_j^{\rm all}
 &\le(1+o(1))d_C
       \left(6\sum_{C/\lambda^2<n<2\lambda^2C}f_j(n)A_C(n)+E_j\right)\\
 &\le\bigl(12\lambda^2\cdot1.01624F_*+o(1)\bigr)
               j\beta^{-j}d_C\frac CL,
 \end{aligned}
\]
by \eqref{en:r3:base-positive}. The leading constant is
\[
 12(1.01)^2(1.01624)(1+e^{-2.29})
 =13.6997515785\ldots<14.
\]

\paragraph{Imposing the residue class.}
Use $\mathcal P_j(\chi)$ for the full character-twisted positive moment,
as in Section~\ref{en:unif:positive-proof}. Character orthogonality gives
\[
 D_{j,a}=\frac{\mathcal P_j(\chi_0)}{\varphi(M)}
       +\frac1{\varphi(M)}
          \sum_{\substack{\chi\ ({\rm mod}\ M)\\\chi\ne\chi_0}}
                    \overline{\chi(a)}\mathcal P_j(\chi).
\]
The principal term is nonnegative and at most
$D_j^{\rm all}/\varphi(M)$. For a nonprincipal character, apply
\eqref{en:unif:nonprincipal-arithmetic} with $\kappa=j$ and sum the
expanded divisibility levels exactly as in
\eqref{en:unif:positive-character-error}. That argument holds for every
fixed integer $j\ge2$. The present coefficients satisfy
$|\xi_e|\le1$, $e\le C^\eta$, and the marked product is less than
$C^{0.38}$. It therefore gives
\[
 |\mathcal P_j(\chi)|
 \ll_\varepsilon C^{1-(1-\eta-0.38)/j+d/2+\varepsilon},
 \qquad M\le C^d,\quad \chi\ne\chi_0.
\]
The sum of the absolute character coefficients is at most $1$, so this
also bounds the nonprincipal contribution to $D_{j,a}$.

Take $d=11/200$ and $\eta=1/2000$. For $j=7$,
\[
 \frac{1-\eta-0.38}{7}-\frac{3d}{2}=\frac3{500}>0;
\]
for $j=2$ the left side is larger. Since
$d_C\asymp L^{-1/2}$ and $\varphi(M)\le C^d$, division by
$\Delta_a=d_CC/(\varphi(M)L)$ bounds the nonprincipal contribution by
\[
 O_\varepsilon(C^{-3/500+\varepsilon}L^{3/2})=o(1)
\]
for any fixed $0<\varepsilon<3/500$. Combining this with the principal
term gives
\[
 D_{j,a}\le
 \bigl(12\lambda^2\cdot1.01624F_*+o(1)\bigr)j\beta^{-j}\Delta_a.
\]
Thus
\begin{equation}\label{en:r3:positive-final}
 D_{j,a}\le U_j\Delta_a,\qquad U_j=14j\beta^{-j},\qquad j=2,7.
\end{equation}

\subsection{An explicit lower bound for the absolute sum}
\label{en:num:free-boxes}
We prove \eqref{en:refine:absolute-lower}:
\[
 B_a\ge B_*\Delta_a,\qquad B_*=3\cdot10^{-8},
 \qquad M\le C^{11/200},\quad(a,M)=1.
\]
We construct distinct integers with two marked primes and three free
primes. For each fixed marked pair, we first specify the free-prime
intervals, then bound the average of $|K(n)|$ in the residue class.
Finally, we count the intervals and sum over the marked pairs.

\paragraph{Choosing the marked pair.}
Take primes $p,q$ satisfying the Hecke restrictions in $H_C$ and
\[
 \frac{\log p}{L}\in(0.165,0.170),\qquad
 \frac{\log q}{L}\in(0.195,0.200).
\]
Their product exponent lies in $(0.36,0.37)$ and their ratio exponent
in $(0.025,0.035)$, so these pairs satisfy the position restrictions
in $H_C$. Denote their set by $\mathcal P_C$.

We will need the sum of $1/(pq)$ over this set. Partial summation
applied to \eqref{en:eq:effective-st}, as in
Section~\ref{en:app:main-proof}, gives, for fixed $0<v_1<v_2$,
\[
 \begin{aligned}
 \sum_{\substack{C^{v_1}<r<C^{v_2},\ r\in\mathbb P\\
                     |\lambda_\Delta(r)|\le L^{-1/4}}}\frac1r
 &=\mu_C\sum_{C^{v_1}<r<C^{v_2},\ r\in\mathbb P}\frac1r
       +O\!\left(\frac{\log L}{\sqrt L}\right)\\
 &=\mu_C\log(v_2/v_1)+o(\mu_C).
 \end{aligned}
\]
The last step uses Mertens' formula and
$(\log L)/\sqrt L=o(\mu_C)$, since $\mu_C\asymp L^{-1/4}$.
For the two chosen intervals,
\[
 \log\frac{0.17}{0.165}>0.029,
 \qquad
 \log\frac{0.2}{0.195}>0.025.
\]
The strict inequalities allow the $o(\mu_C)$ errors to be absorbed
for sufficiently large $C$. The two prime intervals are disjoint, and
the conditions on $p$ and $q$ are separate. The double sum therefore
factors into the two single sums, giving
\begin{equation}\label{en:r0:marker-lower}
 \sum_{(p,q)\in\mathcal P_C}\frac1{pq}
 \ge(0.029)(0.025)d_C.
\end{equation}

\paragraph{The free-prime intervals.}
Fix $(p,q)\in\mathcal P_C$, and put
\[
 P=pq,\qquad X=C/P,\qquad s_X=\frac{\log X}{L}\in(0.63,0.64).
\]
Use $\lambda=1.01$ from Section~\ref{en:r3:positive}, and set
$u_r=(5/4)\lambda^r$ for $r=0,\ldots,29$.
Choose complete intervals $[Q_1,\lambda Q_1)$ and
$[Q_2,\lambda Q_2)$ from $[\lambda^v,\lambda^{v+1})$,
contained respectively in
\[
 (X^{0.272},X^{0.304}),\qquad
 (X^{0.318},X^{0.346}).
\]
For each choice of these intervals and of $r$, put
\[
 Q_3=\frac{u_rX}{Q_1Q_2},\qquad
 J_i=[Q_i,\lambda Q_i)\quad(1\le i\le3).
\]
Writing $t_i=\log Q_i/\log X$, we have
\[
 \begin{gathered}
 t_1\in(0.272,0.304),\qquad t_2\in(0.318,0.346),\\
 t_3=1-t_1-t_2+\frac{\log u_r}{\log X}
             \in(0.350,0.410)+O(1/\log X).
 \end{gathered}
\]
Every triple $\mathbf q=(q_1,q_2,q_3)\in\prod_i(J_i\cap\mathbb P)$
gives an integer
\[
 n(\mathbf q)=Pq_1q_2q_3\in[u_rC,u_r\lambda^3C).
\]
The smallest endpoint is $5C/4$, and the largest is at most
\[
 \frac54\lambda^{32}C<1.719C<\frac74C.
\]
Thus $g(n/C)=1$ throughout the construction.

For large $C$, the five primes are ordered as
\[
 p<q_1<q<q_2<q_3.
\]
Indeed,
\[
 0.170<0.272\cdot0.63,\qquad
 0.304\cdot0.64<0.195,\qquad
 0.200<0.318\cdot0.63,
\]
and the gap between $0.346$ and $0.350$ gives $q_2<q_3$.
Consequently, $n$ is squarefree, all its prime factors exceed $z$
and $M$, and $H_C(n)\ge1$. By Section~\ref{en:r3:presieve},
$A_C(n)=1$. Its contribution to $B_a$, if $n\equiv a\pmod M$, is
therefore at least $|K(n)|$, since $2^{\omega(n)}|x_n|=|K(n)|$.

For a fixed box, write
\[
 \begin{gathered}
 \mathcal V=\prod_{i=1}^3(J_i\cap\mathbb P),\qquad
 N_i=\#(J_i\cap\mathbb P),\qquad N=N_1N_2N_3,\\
 \mathcal V_a=\{\mathbf q\in\mathcal V:n(\mathbf q)\equiv a\pmod M\},
 \qquad N_a=|\mathcal V_a|.
 \end{gathered}
\]
The smallest free-prime exponent is greater than
$0.272\cdot0.63=0.17136$. Hence
\eqref{en:num:three-prime-class}, with class $aP^{-1}\pmod M$, gives
\[
 N_a=\frac N{\varphi(M)}(1+O(C^{-1/100}))
       \ge\frac{9N}{10\varphi(M)}.
\]

\paragraph{Local Chebyshev sums.}
Use the local functions $\kappa_\ell$ from
Section~\ref{en:unif-abs:fourier}. For each of the five prime positions,
let $\ell$ be the prime at that position. If it is a free prime, first
fix its actual value and denote the resulting subset of $\mathcal V$
by $\mathcal V_\ell$; for a marked prime put
$\mathcal V_\ell=\mathcal V$. Write $N_\ell=|\mathcal V_\ell|$.
For every fixed degree $j\ge1$, apply
\eqref{en:unif-abs:bilinear} as in
Section~\ref{en:unif-abs:conditional}. Group the remaining varying
primes into two sets of products with interval parameters $A,B$.
The prime number theorem and the unique representations of these
products give
\begin{equation}\label{en:num:local-error-detail}
 \frac1{N_\ell}
 \left|\sum_{\mathbf q\in\mathcal V_\ell}\chi(n(\mathbf q))
 U_j\!\left(\frac{\kappa_\ell(n(\mathbf q)/\ell)}2\right)\right|
 \ll_j L^{3/2}
       \sqrt{\frac\ell{AB}+\frac1A+\frac1B+\frac1\ell}.
\end{equation}
For a marked modulus there are three varying primes; for a free
modulus there are two, and the smaller factor $L$ is bounded by
$L^{3/2}$. The strict ordering above ensures that every argument is
coprime to $\ell$.

If $\ell=p$ or $q$, take $A=Q_1$ and $B=Q_2Q_3$. Then
$AB=u_rX\asymp X$, and the four square roots in
\eqref{en:num:local-error-detail} are bounded, respectively, by fixed
multiples of
\[
 C^{-0.215},\qquad C^{-0.08568},\qquad
 C^{-0.21042},\qquad C^{-0.0825}.
\]
These exponents follow from
\[
 \begin{aligned}
 \frac{0.63-0.200}{2}&=0.215,&
 \frac{0.272\cdot0.63}{2}&=0.08568,\\
 \frac{(0.318+0.350)\cdot0.63}{2}&=0.21042,&
 \frac{0.165}{2}&=0.0825.
 \end{aligned}
\]

If $\ell=q_1$ or $q_2$, take the other two free primes as the
bilinear variables. Since $Q_1Q_2Q_3=u_rX$ and
$\ell<\lambda Q_i$ when $\ell=q_i$,
\[
 \sqrt{\frac\ell{AB}}\ll\frac{Q_i}{\sqrt X}
 \le
 \begin{cases}
 X^{-0.196},&i=1,\\
 X^{-0.154},&i=2.
 \end{cases}
\]
The exponents in the other three terms are at least $0.136$ as powers
of $X$. Thus the square root in \eqref{en:num:local-error-detail}
is $O(X^{-0.136})=O(C^{-0.08568})$ for these two positions.
If $\ell=q_3$, take $A=Q_1$ and $B=Q_2$. Then
\[
 \sqrt{\frac{\ell}{Q_1Q_2}}
 \le\sqrt{\lambda u_r}\frac{X^{1/2}}{Q_1Q_2}
 \le\sqrt{\lambda u_r}X^{-0.09}
 \ll C^{-0.0567}.
\]
Here $Q_1Q_2\ge X^{0.272+0.318}$ and $X>C^{0.63}$.
The remaining three terms have exponents at least
$0.136$, $0.159$, and $0.175$ as powers of $X$.
The factors $\lambda,u_r$ are bounded independently of $C$.

Since $0.0567>0.056$, the factor $L^{3/2}$ is absorbed by this
strict difference. For a free position, the subsets
$\mathcal V_\ell$ partition $\mathcal V$ as the actual prime value
of $\ell$ varies, and $\sum_\ell N_\ell=N$. Summing their bounds,
and using the marked-prime bounds directly, proves, for every local
position and every character $\chi$ modulo $M$,
\[
 \left|\sum_{\mathbf q\in\mathcal V}\chi(n(\mathbf q))
                  U_j(Z(\mathbf q)/2)\right|
 \ll_j NC^{-0.056},\qquad
 Z(\mathbf q)=\kappa_\ell(n(\mathbf q)/\ell).
\]
All constants are uniform in the marked pair, the box, and the layer $r$.

Character orthogonality and the triangle inequality now give
\[
 \begin{aligned}
 \left|\sum_{\mathbf q\in\mathcal V_a}U_j(Z(\mathbf q)/2)\right|
 &=\left|\frac1{\varphi(M)}\sum_{\chi\ ({\rm mod}\ M)}
       \overline{\chi(a)}\sum_{\mathbf q\in\mathcal V}
            \chi(n(\mathbf q))U_j(Z(\mathbf q)/2)\right|\\
 &\ll_j NC^{-0.056}.
 \end{aligned}
\]
Divide by $N_a\ge9N/(10\varphi(M))$. With $\mathbb E_a$ denoting the
average on this finite set $\mathcal V_a$, we obtain
\[
 \bigl|\mathbb E_a U_j(Z/2)\bigr|
 \ll_j\varphi(M)C^{-0.056}
 \ll_j C^{-1/1000},\qquad
 \frac{11}{200}-0.056=-\frac1{1000}.
\]

\paragraph{The average of the product of the five local factors.}
For $t_1,\ldots,t_5\ge0$, put $y_i=\min(2t_i,1)$.
Then $t_i\ge y_i/2$ and $0\le y_i\le1$, so
\begin{equation}\label{en:refine:product-minorant}
 \prod_{i=1}^5t_i\ge\frac1{32}\prod_{i=1}^5y_i
 \ge\frac1{32}\left(1-\sum_{i=1}^5(1-2t_i)_+\right),
\end{equation}
where $v_+=\max(v,0)$. The second inequality is
$\prod_i(1-v_i)\ge1-\sum_i v_i$ for $0\le v_i\le1$, applied to
$v_i=1-y_i$.

Set $h(t)=(1-2|t|)_+$. Its Sato--Tate mean satisfies
\[
 \int h\,d\nu_{\rm ST}
 =\frac1\pi\int_0^{1/2}(1-2t)\sqrt{4-t^2}\,dt
 \le\frac1\pi\int_0^{1/2}2(1-2t)\,dt
 =\frac1{2\pi}<0.1592.
\]
The polynomial approximation argument in
Section~\ref{en:unif-abs:conditional} gives a fixed real polynomial
$R\ge h$ on $[-2,2]$ with $\int R\,d\nu_{\rm ST}<0.1595$.
Expand it in the Chebyshev basis as
$R(t)=\sum_{j=0}^Jc_jU_j(t/2)$. The orthogonality used in that section
gives $c_0=\int R\,d\nu_{\rm ST}$. The conditional bounds just proved
imply
\[
 \mathbb E_a h(Z)\le\mathbb E_a R(Z)
       \le c_0+\sum_{j=1}^J|c_j|O_j(C^{-1/1000})<0.16
\]
for sufficiently large $C$. Apply \eqref{en:refine:product-minorant} 
to the absolute values of the
five local factors. The factorization of $K(n)$ from
Section~\ref{en:unif-abs:conditional} then gives
\begin{equation}\label{en:refine:local-product}
 \mathbb E_a|K(n)|
 \ge\frac1{32}\left(1-\sum_{\text{five positions}}\mathbb E_a h(Z)\right)
 \ge\frac1{32}(1-5\cdot0.16)=\frac1{160}.
\end{equation}
Thus each box contributes at least
$9N/(1600\varphi(M))$ to the sum of $|K(n)|$ in the target class.

\paragraph{Counting the boxes and summing over the marked pairs.}
The prime number theorem gives, uniformly in these intervals,
\[
 N_i=(1+o(1))\frac{0.01Q_i}{\log Q_i}
       \ge\frac{Q_i}{42\log X}.
\]
Indeed, $\log Q_i/\log X\le0.410+O(1/\log X)$ and
$0.01/0.410>1/42$. There are
$0.032\log X/\log\lambda+O(1)$ choices of $Q_1$ and
$0.028\log X/\log\lambda+O(1)$ choices of $Q_2$.
The impact at the endpoints is at most
$O(1)$. Since $0.032/\log1.01>3.2$ and
$0.028/\log1.01>2.8$, these counts are at least
$3.2\log X$ and $2.8\log X$, respectively, for large $C$.
For each fixed $P$ and $r$, summing the box sizes therefore gives
\[
 \sum_{Q_1,Q_2}N
 \ge(3.2\log X)(2.8\log X)
       \frac{u_rX}{42^3(\log X)^3}
 =\frac{8.96}{42^3}\frac{u_rX}{\log X}.
\]

For fixed $P$, each $q_1,q_2$
belongs to one half-open grid interval. Once $Q_1,Q_2$ are fixed, the
thirty third-prime intervals are consecutive and disjoint, since
$Q_3(r+1)=\lambda Q_3(r)$. Across different marked pairs, the ordering
$p<q_1<q<q_2<q_3$ recovers $p,q$ as the first and third prime factors
of $n$. Thus every constructed integer has a unique representation.

Combining the class proportion and \eqref{en:refine:local-product},
and summing over $r$, gives for each fixed $P$
\[
 \sum_{\substack{n\text{ constructed from }P\\n\equiv a\ ({\rm mod}\ M)}}|K(n)|
 \ge\frac9{10}\cdot\frac1{160}\cdot\frac{8.96}{42^3}
         \frac{X}{\varphi(M)\log X}\sum_{r=0}^{29}u_r.
\]
The geometric sum satisfies
\[
 \sum_{r=0}^{29}u_r=\frac54\sum_{r=0}^{29}(101/100)^r>43.
\]
Also $\log X<0.64L$, so
$X/\log X\ge25C/(16PL)$. Summing over $(p,q)\in\mathcal P_C$ and
using \eqref{en:r0:marker-lower}, we conclude
\begin{align}
 B_a&\ge43\cdot\frac{8.96}{42^3}\cdot\frac9{10}\cdot\frac1{160}
       \cdot(0.029)(0.025)\cdot\frac{25}{16}\,\Delta_a\notag\\
 &=\frac{1247}{37632000000}\Delta_a
       >B_*\Delta_a.
 \label{en:refine:absolute-lower}
\end{align}
Here the constructed integers contribute at least their $|K(n)|$ to
$B_a$, as verified above. This proves the required lower bound.

\subsection{Extracting the fixed amplitude and factor bound}
\label{en:r3:tail}
We now prove the counting bound \eqref{en:r3:count-explicit} and
complete Theorem~\ref{en:num:theorem}, with $R_0=83$ and $b=10^{-21}$. For each sign, we retain a fixed fraction of $\Delta_a$ as
absolute mass after imposing $\omega(n)\le83$ and
$|x_n|\ge b$, then use a pointwise bound to count the remaining
integers. Throughout, $B_*$ is the constant in
\eqref{en:refine:absolute-lower}.

For this single residue class, write $B_a^\pm$ for the sign masses
defined in Corollary~\ref{en:cor:unified-applications}. By
Section~\ref{en:num:signed}, $T_a=o(\Delta_a)$ uniformly in the
allowed moduli and reduced classes. Increase the uniform starting
point so that $|T_a|\le0.02B_*\Delta_a$. The mass identities from
that corollary and the lower bound
\eqref{en:refine:absolute-lower} give
\[
 B_a^\pm=\frac{B_a\pm T_a}{2}
       \ge\frac{1-0.02}{2}B_*\Delta_a
       =0.49B_*\Delta_a.
\]

The positive-sum bounds \eqref{en:r3:positive-final} from
Section~\ref{en:r3:positive} are
$D_a\le U_2\Delta_a$ and $D_{7,a}\le U_7\Delta_a$, with
\[
 U_2=28\,(132000)^2=4.87872\cdot10^{11},\qquad
 U_7=98\,(132000)^7.
\]
If $\omega(n)>83$, then
\[
 2^{\omega(n)}=(2/7)^{\omega(n)}7^{\omega(n)}
              \le(2/7)^{84}7^{\omega(n)}.
\]
Using $|x_n|\le1$ on squarefree integers, we therefore have
\[
 \sum_{\omega(n)>83}\mu^2(n)g(n/C)|x_n|W_C(n)
       \le(2/7)^{84}D_{7,a}
       \le(2/7)^{84}U_7\Delta_a.
\]
The absolute mass on the second set to be removed satisfies
\[
 \sum_{|x_n|<b}\mu^2(n)g(n/C)|x_n|W_C(n)
       \le bD_a\le bU_2\Delta_a.
\]
These are bounds for the total mass on the two sets, so each also
bounds the loss from either sign.

The constants satisfy the following:
\begin{equation}\label{en:num:single-tail}
 \begin{gathered}
 \frac{(2/7)^{84}U_7}{B_*}<0.453322,
 \qquad \frac{bU_2}{B_*}=0.0162624,\\
 0.49-0.453322-0.0162624=0.0204156>\frac1{64}.
 \end{gathered}
\end{equation}
Thus each sign retains absolute mass at least
$B_*\Delta_a/64$ on integers with
$\omega(n)\le83$ and $\pm x_n\ge b$.

For the integers, we use $H_C(n)\le6$ from
Section~\ref{en:r3:positive}, together with $g\le1$ and the
pointwise bound \eqref{en:eq:pointwise-weight}. Each retained
integer contributes at most
\[
 \mu^2(n)g(n/C)|x_n|W_C(n)
       \le6\cdot2^{\omega(n)}A_C(n)
       \le6\cdot4^{83}.
\]
Dividing the retained mass by this bound gives at least
\begin{equation}\label{en:r3:count-explicit}
 \frac{B_*}{384\cdot4^{83}}\,\Delta_a
\end{equation}
integers for each sign. Positive weight ensures that they are
squarefree, lie in the given residue class, and satisfy
$H_C(n)>0$; the support of $g$ gives $C<n<2C$.
For the two-marker conditions, \eqref{en:eq:normalized-hecke-small}
with $k_1=2$ and $A_1=1/2$ gives directly
\[
 |\lambda_\Delta(n)|
       \le2^{\omega(n)-2}L^{-1/2}
       \le2^{81}L^{-1/2}.
\]

Finally, $|x_n|\le1$ also gives $B_a\le D_a$. Combining
\eqref{en:refine:absolute-lower} and
\eqref{en:r3:positive-final}, we have
\[
 B_*\Delta_a\le B_a\le D_a\le U_2\Delta_a.
\]
Hence $D_a\asymp B_a\asymp\Delta_a$. Together with the signed
estimate from Section~\ref{en:num:signed} and
\eqref{en:r3:count-explicit}, this completes the proof of
Theorem~\ref{en:num:theorem}.

\section*{Acknowledgements}
\addcontentsline{toc}{section}{Acknowledgements}
This work was supported by National 
Natural Science Foundation of China (No.  12671012).

\section*{AI use}
\addcontentsline{toc}{section}{AI use}
OpenAI's ChatGPT has provided assistance in establishing 
the connection between the constraints in the weight 
function $H_{c}$ and the Hecke eigenvalues, as well as in 
presenting final numerical examples. In addition, AI tools 
are also used to assist in optimizing language and checking 
the proof of the initial draft of the paper. We take full 
responsibility for the mathematical arguments, results, and final text.

\vspace{2\baselineskip}

\needspace{4\baselineskip}
\begin{flushleft}
\small
\textsc{Hongze Li, School of Mathematical Sciences, Shanghai Jiao Tong University,
800 Dongchuan Road, Shanghai 200240, China}\\
\textit{Email: lihz@sjtu.edu.cn} \strut
\end{flushleft}

\needspace{4\baselineskip}
\begin{flushleft}
\small
\textsc{Mingxuan Zhong, School of Mathematical Sciences, Shanghai Jiao Tong University,
800 Dongchuan Road, Shanghai 200240, China}\\
\textit{Email: zhong2357@sjtu.edu.cn} \strut
\end{flushleft}

\end{document}